\documentstyle[amssymb]{article}
\newcommand{\R}{{\Bbb R}}
\newcommand{\p}{\frac{N+2}{N-2}}
\newcommand{\q}{\frac{2N}{N-2}}
\newcommand{\e}{\varepsilon}
\newcommand{\al}{\alpha}

\newcommand{\bR}{\bar{R}}
\newcommand{\bS}{\bar{S}}
\newcommand{\bal}{\bar{\alpha}}
\newcommand{\baf}{\bar{f}}
\newcommand{\be}{\bar{\e}}

\newcommand{\ba}{\bar{a}}
\newcommand{\bq}{\bar{q}}

\newcommand{\tw}{\tilde{w}}
\newcommand{\RNS}{\R^N \setminus \Sigma}
\newcommand{\LLe}{{\Bbb L}_{\epsilon}}
\newcommand{\LLej}{{\Bbb L}_{\epsilon,j}}
\newcommand{\Le}{{\cal L}_{\be, \bR}}
\newcommand{\Lei}{{\cal L}_{\be, \bR}^{-1}}

\newcommand{\Laa}{\Lambda_{\be,\bR,\ba}}
\newcommand{\LaO}{\Lambda_{\be,\bR,0}}

\newcommand{\tcM}{\tilde{\cal M}}
\newcommand{\tLaa}{\tilde{\Lambda}_{\be,\bR,\ba}}
\newcommand{\tLaO}{\tilde{\Lambda}_{\be,\bR,0}}

\newcommand{\bu}{\bar{u}_{\be}}

\newcommand{\bue}{\bar{u}_{\be} (\bR,\ba ,\cdot)}
\newcommand{\buex}{\bar{u}_{\be} (\bR,\ba ,x)}

\newcommand{\ui}{u_{\e_i}(R_i,a_i,\cdot -x_i)}
\newcommand{\uix}{u_{\e_i}(R_i,a_i,x-x_i)}

\newcommand{\bw}{\bar{w}_{\be}(\bR,\cdot)}
\newcommand{\bwx}{\bar{w}_{\be}(\bR,x)}

\newcommand{\Te}{{\cal T}_{\be, \bR}}
\newcommand{\Tei}{{\cal T}_{\be, \bR}^{i}}
\newcommand{\Se}{{\cal S}_{\be, \bR}}
\newcommand{\TO}{{\cal T}_{0}}
\newcommand{\TOi}{{\cal T}_{0}^{i}}
\newcommand{\SO}{{\cal S}_{0}}

\newcommand{\del}{\partial}

\newtheorem{theorem}{Theorem}
\newtheorem{lemma}{Lemma}
\newtheorem{proposition}{Proposition}
\newtheorem{corollary}{Corollary}
\newtheorem{definition}{Definition}
\newtheorem{remark}{Remark}

\begin{document}

\title{Constant scalar curvature metrics with isolated singularities}
\date{ } 
\author{Rafe Mazzeo\thanks{Supported by the NSF with a Young Investigator 
Fellowship and by grant \#DMS-9303236}\\ Stanford University \\  
mazzeo\@@math.stanford.edu 
\and Frank Pacard \\ Universit\`e de Paris XII \\ 
Frank.Pacard\@@diam1.ens-cachan.fr}  
\maketitle
\begin{abstract}
We extend the results and methods of \cite{MP} to prove the existence of 
constant positive scalar curvature metrics $g$ which are complete and
conformal to the standard metric on $S^N \setminus \Lambda$, where
$\Lambda$ is a disjoint union of submanifolds of dimensions between
$0$ and $(N-2)/2$. The existence of solutions with isolated singularities
occupies the majority of the paper; their existence was previously
established by Schoen \cite{S}, but the proof we give here, based
on the techniques of \cite{MP}, is more direct, and provides more
information about their geometry. When $\Lambda$ is
discrete we also establish that these solutions are smooth points in the 
moduli spaces of all such solutions introduced and studied in \cite{MPU1} and
\cite{MPU2}
\end{abstract}

\section{Introduction and statement of the results}

In this paper we construct solutions of the Yamabe problem on the sphere
$(S^N,g_0)$ with its standard metric which are singular at a specified 
closed set $\Lambda \subset S^N$. More specifically, we seek a new metric 
$g$ conformal to $g_0$ and complete on $\Lambda \subset S^N$, and with 
constant scalar curvature $R$. 
The problem may be translated into a differential equation as follows.
Since $g$ is conformal to $g_0$, we may write $g = u^{\frac{4}{N-2}}g_0$
where $u$ is a positive function on $M \setminus \Lambda$. The scalar
curvature functions $R(g_0)$ of $g_0$ and $R(g)$ of $g$ are related by 
the equation 
\begin{equation}
\Delta_{g_0}u - \frac{N-2}{4(N-1)}R(g_0)u + \frac{N-2}{4(N-1)}R(g)u^{\p} = 0.
\label{eq:1.1}
\end{equation}
In order that $g$ be complete on $S^N \setminus \Lambda$ it is necessary
for $u$ to tend to infinity rather strongly on approach to $\Lambda$, and
of course, we wish to solve this equation with $R$ a (prescribed) constant.

The first two terms of the operator on the left in (\ref{eq:1.1}), namely
\begin{equation}
{\cal L}_{g_0} \equiv \Delta_{g_0} - \frac{N-2}{4(N-1)}R(g_0)
\label{eq:1.5}
\end{equation}
give a second order linear elliptic differential operator known as the 
conformal Laplacian of the metric $g_0$. It satisfies the conformal 
equivariance property that if two metrics are conformally 
related, such as $g$ and $g_0$ above, then for any function $\phi$,
\begin{equation}
{\cal L}_{g_0}(u\phi) = u^{\frac{N+2}{N-2}}{\cal L}_{g}(\phi).
\label{eq:1.6}
\end{equation}
Notice that (\ref{eq:1.1}) corresponds to the special case of
(\ref{eq:1.6}) when $\phi = 1$. 

This `singular Yamabe problem' has been extensively studied 
in recent years, also in the case when the ambient manifold 
is more general than the sphere, and many existence results as well as 
obstructions to existence are known. Briefly, for a solution to exist 
on a general compact Riemannian manifold $(M,g_0)$, the size of $\Lambda$ 
and the sign of $R$ must be related to one another~:  
if a solution exists with $R < 0$, then $\dim(\Lambda) > \frac{N-2}2$,
while if a solution exists with $R \ge 0$, then $\dim(\Lambda) 
\le \frac{N-2}2$ and in addition the first eigenvalue of the conformal 
Laplacian of $g_0$ must be nonnegative. Here the dimension is to
be interpreted as Hausdorff dimension. Unfortunately, only partial
converses to these statements are known. In particular, almost all of the
existence results require $\Lambda$ to be a submanifold, or at worst
a stratified set. 

The character of the analysis required to prove existence of solutions
when $R < 0$ is quite different than in the other cases. In fact,
it is not too difficult, using barrier methods, to construct solutions
of (\ref{eq:1.1}). The difficulty is in showing that these solutions
have a singularity at $\Lambda$ strong enough for $g$ to be
complete on the complement. Similarly, when $R=0$, (\ref{eq:1.1}) becomes 
linear, albeit of a somewhat nonstandard sort. We shall not
discuss these cases any further, but refer to the work of R. McOwen 
\cite{McO}, D. Finn \cite{F}, \cite{F2} and the first 
author \cite{M}, and references therein. 

We shall concentrate here on the case $R > 0$. The first examples of solutions
were constructed by R. Schoen \cite{S} when $\Lambda$ is either a finite set 
of points (of cardinality greater than one) or a nonrectifiable 
set arising essentially as the limit set of a Kleinian group. 
Nontrivial solutions with higher dimensional singular sets were 
constructed by the first author and N. Smale \cite{MS}, the second author 
\cite{P}, and finally in some generality by both of us \cite{MP}. This 
last result states that solutions may be constructed on an arbitrary
compact manifold of nonnegative scalar curvature $M$ whenever $\Lambda$ is 
a finite disjoint union of submanifolds of dimensions between $1$ and 
$\frac{N-2}2$. This paper is meant to extend the methods of \cite{MP} to also 
allow for the construction of solutions with isolated singularities in 
the case where $M=S^N$ and $g_0$ is the standard metric on $S^N$, i.e. to
allow $\Lambda$ to be an arbitrary disjoint finite union of submanifolds
of dimensions between $0$ and $\frac{N-2}2$. Our main result is:
\begin{theorem}
Suppose that $\Lambda$ is a disjoint union of submanifolds in $S^N$
of varying dimensions. Write $\Lambda  = \Lambda' \cup \Lambda''$,
where $\Lambda'$ is the union of all submanifolds of dimension zero,
i.e. $\Lambda' = \{p_1, \ldots, p_n\}$ is a collection of points,
and $Lambda'' = \bigcup_{j=1}^n \Lambda_j$ where $\dim \Lambda_j = k_j$. 
Suppose further that $0 < k_j \le \frac{N-2}2$ for each $j$, and 
either $n = 0$ or $n \ge 2$. Then there exists a complete metric $g = 
u^{\p}g_0$ on $S^N \setminus \Lambda$ which has constant positive scalar 
curvature $R = N(N-1)$. 
\label{th:1.1}
\end{theorem}

For most of this paper we shall consider the case where $\Lambda''$ is
empty, so that $\Lambda$ is a finite collection of points
\[
\Lambda = \{p_1, \ldots, p_n\} \subset S^N.
\]
The modifications needed to treat the general case are an amalgamation 
of the techniques here and those of \cite{MP}, and shall be described in 
a brief final section of this paper.

Of course, solutions with isolated singularities were already constructed 
by R. Schoen, but his remarkable proof is long and difficult. We feel that 
(the extension of) the methods of \cite{MP} avoid many of the difficulties
he encountered and substantially clarify the analysis. In addition, various 
properties of the solutions follow immediately from the construction here, 
but are not at all obvious for Schoen's solutions; we describe some
of these below.

Let us describe the two approaches, and some of their differences, somewhat
more specifically. If we normalize the desired constant scalar curvature 
$R$ to be $N(N-1)$, then the equation we are trying to solve on the sphere is 
\begin{equation}
{\cal N}(u) \equiv \Delta u - \frac{N(N-2)}{4} u + \frac{N(N-2)}{4} 
u^{\frac{N+2}{N-2}} = 0.
\label{eq:1.2}
\end{equation}
In order to construct a solution with a given singular set $\Lambda$,  
first a one-parameter family of approximate solutions $\tilde{u}_\e$ of 
(\ref{eq:1.2}), each element of which blows up quickly enough near $\Lambda$,
is constructed. Then (\ref{eq:1.2}) is linearized
about $\tilde{u}_\e$, and the resulting linear operator $L_\e$ is analyzed.
If it can be proved that $L_\e$ is surjective on some reasonable space
of functions, at least when $\e$ is sufficiently small, then a standard 
iteration argument may be used to correct $\tilde{u}_\e$ by adding to it 
a function $v$ to obtain an exact solution to (\ref{eq:1.1}) which blows 
up sufficiently quickly at $\Lambda$. Unfortunately, $L_\e$ is not surjective 
on $L^2$; in fact, $L_\e$ is self-adjoint on $L^2$, but $0$ is in its 
continuous spectrum, so that it does not even have closed range on
this space. Schoen's tactic is to find an explicit infinite set
of functions which span an approximate nullspace $K$, such that $L_\e$ 
is invertible on $K^\perp$. He first solves the equation 
on $K^\perp$, and then gives `balancing' conditions to ensure
that the solution of this restricted problem is a solution
of the original problem. Our somewhat different tactic is to work on 
a finite dimensional extension of a certain weighted H\"older space.
On this space, $L_\e$ is actually an isomorphism, when $\Lambda$ is
discrete, and surjective in general, for $\e$ sufficiently small. 
Unfortunately, as $\e$ tends to zero, the norm of any right inverse for 
$L_\e$ blows up (Schoen encounters a similar problem). We must then analyze
precisely the rate and manner of blowup. This is different than in \cite{MP}, 
when $\Lambda'$ is empty, and 
where the (right) inverse for $L_\e$ is bounded as $\e \rightarrow 0$. 
The advantage of working in weighted spaces over using the approximate
nullspace $K$ is that in our approach there are only a finite
number of obstructions to solving the equation, and these may
be identified explicitly and geometrically. In fact, these 
obstructions are intimately connected with the definition of the
finite dimensional extension mentioned above. 

As noted earlier, another advantage of our approach is that we
may easily derive various properties of the solutions. The main
property we are interested in is nondegeneracy, which will
be defined precisely in \S 11.  This property is important in the
study of the `marked' and `unmarked' moduli spaces ${\cal M}_{\Lambda}$ 
and ${\cal M}_n$ of solutions of this problem where the singular
set $\Lambda$ is fixed or allowed to vary amongst all configurations
of $n$ points in $S^N$. There are notions of marked and unmarked 
nondegeneracy associated to each of these spaces. These moduli spaces 
were defined and studied by the first author, with D. Pollack and 
K. Uhlenbeck in \cite{MPU1}, \cite{MPU2}. It was proved that for 
$M = S^N$ they are real analytic sets. If there exists some $g \in 
{\cal M}_\Lambda$ which is marked nondegenerate, then the top dimensional 
stratum in the component of $g$ is a real analytic manifold of dimension $n$. 
Similarly, if $g \in {\cal M}_n$ is unmarked nondegenerate, then nearby $g$
this moduli space is a real analytic manifold of dimension $n(N+1)$. 
Nondegenerate solutions on $S^N$ of a very special type were 
constructed in \cite{MPU2}; for these solutions, only certain special 
configurations $\Lambda$ (in particular, only those with cardinality 
$n$ an even number and with points clustered in pairs) are allowed. 
In contrast, we prove
\begin{theorem}
For any integer $n \ge 2$ and any configuration $\Lambda$ of $n$ points
in $S^N$ there exists an element $g \in {\cal M}_n$ which has
singular set $\Lambda$ and which is unmarked nondegenerate. 
For a generic (in fact, Zariski open) set of $\Lambda$, this solution
is marked nondegenerate, and for such a metric the points $(p_1, \ldots, 
p_n)$ of $\Lambda$ and the `Delaunay necksizes' $(\e_1, \ldots, \e_n)$
constitute a full set of coordinates in ${\cal M}_n$ near $g$,
while the Delaunay parameters alone yield coordinates in 
${\cal M}_{\Lambda}$ near $g$.
\label{th:1.10}
\end{theorem}

The Delaunay necksizes will be defined in \S 2. We remark that the
admissible sets of Delaunay necksize parameters $\{\e_1, \ldots, \e_n\}$
in the construction are not arbitrary sets of small numbers, but are 
required to satisfy a `balancing conditions' (\ref{eq:4.1}). 

Although the statement of Theorem~\ref{th:1.1} is for metrics on the 
complement in the sphere of a set $\Lambda$, we shall use the 
conformal equivariance of the equation (\ref{eq:1.2}) and prove
instead the existence of solutions of this equation on the 
complement in $\R^N$ of a (finite) set $\Sigma$ which decay
at infinity. This makes the technical details somewhat easier
in that we always have a preferred coordinate system. 

The solutions we construct here are not `the same' as the ones constructed 
by R. Schoen; they are also quite different from the ones constructed 
in \cite{MPU2}, and possibly don't even lie in the same components
of the relevant moduli spaces. The ones in the latter paper are
not required to have small Delaunay parameters, while (at least some of) 
the ones in \cite{S} may be thought of as infinite strings of spheres,
connected by very small necks, and joined together at a central convex 
sphere. Our solutions also have small necks, but the central region is 
metrically both concave and very small. 

Karen Uhlenbeck has informed us that she too has established the
existence of solutions to this problem with exactly one singular
point on compact Riemannian manifolds $(M,g_0)$ with nonnegative 
conformal Laplacian.  
The construction of the approximate solutions in this case is 
substantially different than ours, and requires the positive mass theorem,
but her linear analysis is in roughly the same spirit as ours.  

The connection between this problem and the construction of complete,
noncompact surfaces in $\R^3$ with constant mean curvature (CMC) 
is well-known. Indeed, Kapouleas' initial construction of these 
surfaces \cite{K} is related to Schoen's construction of constant 
scalar curvature metrics \cite{S}, and the analysis of the moduli 
space problem in \cite{KMP} for the CMC case was directly inspired
by \cite{MPU1}.  Maintaining this tradition, we shall show, in
a work in preparation, the existence of complete, noncompact CMC
surfaces in $\R^3$ with an arbitrary number of ends (greater than two),
and with prescribed Delaunay parameters on these ends. These solutions
are nondegenerate, hence are smooth points in the relevant moduli
spaces. The geometric difference between these solutions and those
of Kapouleas are the concavity in the central region of our
solutions, versus the convexity of his. 

The outline of this paper is as follows. In the first section 
we examine the basic models of singular solutions, the 
Delaunay solutions (the reason for this name is explained 
in \cite{MPU1}), and prove some estimates on them that will be 
required later. After this, we present the construction of
approximate solutions for the problem. These functions are 
periodic near each end in a certain natural coordinate system, and 
correspondingly the linearized operators we must analyze have 
periodic coefficients. Although an analysis of such operators was 
made in \cite{MPU1}, we need rather stronger behaviour than seems 
to be immediately available from the methods there, so the next few 
sections are devoted to a somewhat novel construction of inverses 
for the linearized operators. First of all, the linearized operators
are invertible only on rather special function spaces. Then, we 
construct an inverse for the Dirichlet problems for the linearizations in the 
$\e$-neighbourhoods of the singular points, and on the complement 
of these neighbourhoods.  By analyzing the Dirichlet-to-Neumann operators 
on the union of the boundaries of these neighbourhoods, we may join these 
right inverses together to get a global right inverse. To carry this
out, we take a slightly indirect path and prove the existence of inverses
by this method only for a model problem; the existence of inverses
for the true problem is deduced by a rather intricate sequence of
perturbations. After proving sufficiently fine estimates for these inverses, 
a rather standard contraction mapping argument (estimates for which,
unfortunately, are not so trivial) is employed to complete the construction.
After this we discuss the issue of nondegeneracy of the solutions. 
As noted earlier, in most of the paper we discuss only the case
where $\Lambda$ is discrete, but in the last section we discuss
the changes needed to handle the general case.

\section{Delaunay type solutions}

In this section we recall some well known fact about the Delaunay solutions of
(\ref{eq:1.1}) on $S^N$ that will be used extensively in
the rest of the paper.  A reference for facts not proved here is 
\cite{MPU1}, cf. also \cite{S2}.  It is known that if
$w$ is any solution of (\ref{eq:1.1}) on $S^N \setminus \{p_1, p_2\}$,
then it is invariant with respect to any conformal transformation
fixing these points; if these points are antipodal, as may be
assumed, then $u$ is rotationally invariant. In either case,
the equation it satisfies may be reduced to an ODE. 
It is convenient to stereographically project the sphere $S^N$
to ${\R}^N$ from one of the singular points, say $p_1$, so
that $p_2$ is sent to $0 \in {\R}^N$. Then the solution is
transformed to a radial solution of 
\begin{equation}
{\cal N}(u) \equiv \Delta u +\frac{N(N-2)}{4}u^{\p}=0,
\label{eq:2.1}
\end{equation}
on ${\R}^N\setminus \{0\}$ which is singular at the origin.

We may reduce (\ref{eq:2.1}) by writing 
\[
u(x)=| x |^{\frac{2-N}{2}}v(-\log |x|),
\]
and using $t = -\log x$ to get
\begin{equation}
\frac{d^2v}{dt^2} -\frac{(N-2)^2}{4} v +\frac{N(N-2)}{4}v^{\p}=0.
\label{eq:2.3}
\end{equation}
This equation is nondissipative, and the Hamiltonian energy 
\begin{equation}
H(v,\dot{v})= \dot{v}^2 - \frac{(N-2)^2}{4} v^2 +\frac{(N-2)^2}{4} v^{\q}
\label{eq:2.4}
\end{equation}
is constant along solutions of (\ref{eq:2.3}). Here and in all
that follows, $\cdot$ denotes differentiation with respect to $t$.
By examining the level curves of $H$, we see that all bounded solutions 
of (\ref{eq:2.3}) lie in the region of the $(v,\dot{v})$-plane where
$H(v, \dot{v}) \leq 0$.  There are several different types of
bounded solutions; we summarize their basic properties:
\begin{proposition}
For any $H_0 \in (-\frac{N-2}{2}(\frac{N-2}{N})^{N/2}, 0)$, there 
exists a unique bounded solution of (\ref{eq:2.3}) satisfying 
$H(v, \dot{v})=H_0$, $\dot{v} (0)=0$ and $\ddot{v}(0) >0$.
This solution is periodic and for all $t \in {\R}$ we have 
$v(t) \in (0,1)$. This solution can be indexed by the parameter $\e =v(0)
\in (0,(\frac{N-2}{N})^{(N-2)/4})$, which is the smaller of the two 
values $v$ assumes when $\dot{v}=0$.
When $H_0 = -\frac{N-2}{2}(\frac{N-2}{N})^{N/2}$, there is a unique bounded 
solution of (\ref{eq:2.3}), given by 
\[
v(t) = \big(\frac{N-2}{N}\big)^{(N-2)/4}.
\]
Finally, if $v$ is a solution with $H_0=0$ then either 
$v(t)\equiv 0$ or $v(t) =(\cosh (t-t_0))^{\frac{2-N}{2}}$ for some 
$t_0 \in {\R}$.
\label{pr:2.1}
\end{proposition}

We will call these the Delaunay type solutions. Although we do not
know them explicitly, the next two propositions give sufficient information
about their behaviour as $\e$ tends to zero for our purposes. 

\begin{proposition}
Fix $\e \in (0,(\frac{N-2}{N})^{(N-2)/4})$ and let $v_\e$ be
the corresponding Delaunay solution. Then the period $T_\e$ of $v_\e$ 
tends to infinity monotonically as $\e \rightarrow 0$ and satisfies
\[
T_{\e} = - \left(\frac{4}{N-2} + o(1)\right) \log (\e).
\]
In addition, for all $t \in {\Bbb R}$
\[
\e \leq v_\e (t) \leq \e \cosh (\frac{N-2}{2} t)
\] 
\label{pr:2.2}
\end{proposition}
{\bf Proof :}
The second claim is rather simple. Since $H$ is constant along
solutions, 
\begin{equation}
H(v_\e,\dot{v}_\e)= \dot{v}_\e^2 - \frac{(N-2)^2}{4} (v_\e^2 -v_\e^{\q})=
-\frac{(N-2)^2}{4} (\e^2 -\e^{\q}).
\label{eq:2.5}
\end{equation}
Moreover, it is easy to see that $v_\e$ is increasing for $t \in 
[0, T_\e/2]$ and decreasing for $t \in [T_\e/2, T_\e]$. In particular the 
minimum of $v_\e$ is achieved at $t=0$ and equals $\e$.
Using this lower bound in (\ref{eq:2.5}) we see that
\[
\dot{v}_\e^2 = \frac{(N-2)^2}{4} \big((v_\e^2 -\e^2) - ( v_\e^{\q}-\e^{\q})
\big) \leq \frac{(N-2)^2}{4}(v_\e^2 -\e^2). 
\]
Taking the (positive) square root and integrating this differential
inequality yields the inequality
$v_\e (t) \leq \e \cosh (\frac{N-2}{2} t)$. 

Now we derive the asymptotic estimate for the period $T_\e$.  Let
$v_{\max,\e}$ denote the larger of the two solutions of 
\[
v^2 - v^{\frac{2N}{N-2}}=\e^2 - \e^{\frac{2N}{N-2}}.
\]
(Hence $v_{\max,\e}$ is the larger of the two values $v_\e$ assumes
when $\dot{v}_\e = 0$). By (\ref{eq:2.5}), we get 
\[
\int_{0}^{T_\e/2} ( v_\e^2 - v_\e^{\q} -\e^2 +\e^{\q})^{-1/2} 
\dot{v}_\e \,dt = \frac{N-2}{4} T_\e.
\]
In the interval of integration we may change the
variable of integration to $\bar{v}=v_\e (t)/\e$. This yields
\[
\int_{1}^{v_{\max,\e}/\e} (\bar{v}^2 -1 - \e^{\frac{4}{N-2}}
(\bar{v}^{\q}-1))^{-1/2} d\bar{v} = \frac{N-2}{4} T_\e.
\]
It is then not difficult to see that this integral grows as 
$-(1+ o(1))\log \e$. 
A proof of the fact that the period is a monotone function of $\e$ is 
indicated in \cite{MPU1} and also proved by elementary ODE methods
in an unpublished manuscript of Chouikha and Weissler.
\hfill $\Box$

\begin{proposition}
For any $\e \in (0,(\frac{N-2}{N})^{(N-2)/4})$ and for any $t \in \R$ 
the Delaunay solution $v_\e$ satisfies the estimates
\[
\left|v_\e (t) - \e \cosh \left(\frac{N-2}{2} t\right)\right| 
\leq c_N \e^{\frac{N+2}{N-2}} e^{\frac{N+2}{2} |t|},
\]
\[
\left| \dot{v}_\e (t) - \frac{N-2}{2}\e \sinh \left(\frac{N-2}{2} t\right)
\right| \leq c_N \e^{\frac{N+2}{N-2}} e^{\frac{N+2}{2} |t|}, 
\]
for some constant $c_N >0$ which depends only on $N$.
\label{pr:2.3}
\end{proposition}
{\bf Proof :} We start with the {\it a priori} estimate 
\[
v_\e (t) \leq \e \cosh (\frac{N-2}{2} t)\leq \e e^{\frac{N-2}{2} t},
\]
valid for $t >0$. Next, writing the equation for $v_\e$ as
\[
\frac{d^2v}{dt^2} -\frac{(N-2)^2}{4} v = -\frac{N(N-2)}{4}v^{\p},
\]
we can represent $v_\e$ as 
\[
v_\e (t) = \e \cosh (\frac{N-2}{2} t) -\frac{N(N-2)}{4} 
e^{\frac{N-2}{2}t}\int_0^{t} e^{(2-N)s}\int_0^s 
e^{\frac{N-2}{2}z}v_\e^{\p}(z)\,dz\,ds.
\]
This leads immediately to the estimate
\[
0 \leq  \e \cosh (\frac{N-2}{2} t) -v_\e (t) \leq \frac{N-2}{8} 
\e^{\frac{N+2}{N-2}} e^{\frac{N+2}{2}t}.
\]
Finally, differentiating the formula of $v_\e$ with respect to $t$, we get
\[
\left| \dot{v}_\e (t) - \frac{N-2}{2} \e \sinh \left(\frac{N-2}{2} t\right)
\right| \leq \frac{N(N-2)}{4} \left( \frac{N-2}{2} e^{\frac{N-2}{2}t}
\int_0^{t}e^{(2-N)s}\int_0^s e^{\frac{N-2}{2}z}v_\e^{\p}(z)\,dz\,ds \right.
\]
\[
\left.
+ e^{\frac{2-N}{2}t} \int_0^t e^{\frac{N-2}{2}z}
v^{\p}(z)dzds\right)\leq \frac{N^2-4}{16} 
\e^{\frac{N+2}{N-2}}e^{\frac{N+2}{2}t}.
\]
\hfill $\Box$ 

\begin{remark}
The estimate of this Lemma is only interesting in the domain where 
\[
2c_N \e^{\frac{N+2}{N-2}}e^{\frac{N+2}{2}|t|} < (N-2) \e e^{\frac{N-2}{2}|t|},
\]
i.e. where $|t| < -\frac{2}{N-2} \log \e -\tilde{c}_N$, which is close to half a period of $v_{\e}$. 
\end{remark}

We transfer these estimates and this remark back to the $x$-coordinates
to obtain
\begin{corollary}
For any $\e \in \big(0, \left(\frac{N-2}{N}\right)^{(N-2)/4}\big]$ and any
$x \in \R^N \setminus \{0\}$, the Delaunay solution $u_\e(x)$
satisfies the estimates
\[
\left|u_{\e}(x) - \frac{\e}{2}\big( 1 + |x|^{2-N} \big)\right| \le 
C_N \e^{\frac{N+2}{N-2}}|x|^{-N},
\]
\[
\left| r\del_r u_\e(x) + \frac{N-2}{2} \e  |x|^{2-N}  \right|
\le C_N \e^{\frac{N+2}{N-2}} |x|^{-N}.
\]
In particular, in the annulus $\e^{\frac{2}{N-2}} \le |x| \le 1$, 
$u_\e$ is well-approximated by $\frac{\e}{2}|x|^{2-N}$.
\label{co:estann}
\end{corollary}

There are some important variations of these solutions, leading 
to a $2N+2$-dimensional family of Delaunay type solutions. These 
variations are families of solutions $U(s)$ of ${\cal N}(U(s)) = 0$
with $U(0) = u_\e$, depending smoothly on the parameter $s$.
The derivatives of these families with respect to $s$ at $s = 0$ 
correspond to Jacobi fields, i.e. solutions of the linearization of 
${\cal N}$ about one of the $u_{\e}$, and will be described below.

We describe these families of variations in turn. The first is quite obvious. 
It is the family where the Delaunay parameter $\e$ is varied
\[
(-\e, 1 - \e) \ni \eta \longrightarrow  u_{\e + \eta}(x).
\]
The second corresponds to the fact that if $u$ is any solutions of 
${\cal N}(u) = 0$, then $R^{\frac{2-N}{2}}u(x/R)$ also solves this 
equation. Applying this to a Delaunay solution yields the family
\[
\R^+ \ni R \longrightarrow |x|^{\frac{2-N}{2}}v_{\e}(-\log |x| + \log R).
\]
The other two families of solutions correspond to translations. The simpler
of these is the usual translation 
\[
\R^N \ni b \rightarrow u_\e(x + b).
\]
The final one corresponds to translations `at infinity'. To describe this
we use the Kelvin transform, $u(x) \rightarrow |x|^{2-N} u(\frac{x}{|x|^2})$,
which preserves the property of being a weak solution of (\ref{eq:2.1}).
Start with a Delaunay solution $u_\e(x)$. Its Kelvin transform is 
\[
|x|^{\frac{2-N}{2}} v_\e (\log |x|).
\]
Translate this by some $a \in \R^N$ to get
\[
|x -a|^{\frac{2-N}{2}} v_\e (\log |x -a|), 
\]
which has its singularity at $a$ rather than $0$. Its Kelvin transform 
yields the family
\[
\R^N \ni a \rightarrow \left| x - a|x|^2 \right|^{\frac{2-N}{2}}
v_\e(-2\log|x| + \log | x - a|x|^2|).
\]
Each function in this family has singularity at $0$ again. 

In fact, the first and third variations, where the Delaunay parameter is 
changed or where the singularity is moved, are less well behaved than
the other two: the Jacobi fields corresponding to them grow too 
quickly. Thus we shall focus attention on the smaller family of solutions
\begin{equation}
u_\e (R,a,x) = \left| x-a|x|^2 \right|^{\frac{2-N}{2}} v_\e (-2\log |x| + 
\log |x-a|x|^2 |+\log R),
\label{eq:4.def}
\end{equation}
where only translations along the Delaunay axis and of the `point 
at infinity' are allowed. 

\section{Construction of the approximate solutions}

We now proceed to define a family of approximate solutions to the
problem, using the Delaunay solutions of the last section. Although
we ultimately wish to construct solutions on $S^N \setminus \Lambda$,
we shall use the conformal invariance of the problem and construct
solutions on $\R^N \setminus \Sigma$, where 
\[
\Sigma = \{x_1, \ldots , x_n\}
\]
with appropriate regularity at infinity. This is done purely to
make the notation as simple as possible. 

The approximate solutions we define here have Delaunay type singularities 
in a ball $B(x_i,\rho_i)$ around each point of $\Sigma$ and are 
harmonic outside the union of slightly larger balls. We use
not just the radial Delaunay solutions, but elements of the
family $u_{\e}(R,a,x)$ defined above. The additional parameters
in this family yield, on the linear level, the Jacobi fields 
constituting a `deficiency space' needed to obtain surjectivity of the 
linearized operator. In addition, these extra parameters are essential
in ensuring a sufficiently good agreement in the transition regions
around the boundaries of the small balls. We emphasize again
that we do not use all possible Delaunay solutions in each ball,
but just two of the four possible families. 

Let
\[
\be = (\e_1, \ldots, \e_n), \quad \bR = (R_1, \ldots, R_n), \quad
\bar{\rho} = (\rho_1, \ldots, \rho_n), \quad \mbox{and}\quad 
\ba = (a_1, \ldots, a_n)
\]
be sets of (small) Delaunay parameters, positive real numbers corresponding
to translations of the Delaunay solutions, (small) radii and vectors 
in $\R^N$, respectively. We shall impose various constraints on these
quantities. First, assume that the balls $B(x_i,\rho_i)$ do not
intersect one another. We will denote by $\Omega_{\bar{\rho}}$ the set 
${\R}^N \setminus \cup_{i=1}^n B(x_i, \rho_i)$. Next, we shall assume
that 
\begin{equation}
\rho_i = \e_i^{\frac{4}{N^2 - 4}}.
\label{eq:defrho}
\end{equation}
Finally, fix a set of positive numbers $q_1, \ldots, q_n$, and assume
that 
\begin{equation}
\e_i = \e q_i, \qquad i =1, \ldots, n, 
\label{eq:defeps}
\end{equation}
or simply $\be = \e \bq$, for some $\e > 0$. 
We shall call a set $\be$ satisfying (\ref{eq:defeps}) {\it admissible}.
Notice that (\ref{eq:defrho}) and (\ref{eq:defeps}) imply that
the quantities $\e_1, \ldots, \e_n$ are all comparable to one another,
and similarly for the $\rho_i$. Because of this, it makes sense to estimate
various quantities later by $\e$ or $\rho = \e^{4/(N^2-4)}$. 

We shall now define the family of approximate solutions $\bue$. 
In each ball $B(x_i,\rho_i)$, $\buex$ is equal to $\uix$.
On the exterior region $\Omega_{\bar{\rho}}$, $\buex$ will equal
some harmonic function $w(x)$. With this definition, it is possible to 
make $\buex$ continuous across each $\del B(x_i,\rho_i)$, but not 
${\cal C}^1$. For example, if one takes $w(x)$ to be the unique harmonic 
function which decays at infinity and whose value on $\del B(x_i,\rho_i)$ 
is given by $\uix$ then, by definition $\buex$ is continuous but not 
${\cal C}^1$. It is easier to deal with a harmonic function $w$ which is 
a sum of multiples of Green functions with poles at the $x_i$.  
The next result shows that we can mitigate this discontinuity of 
derivatives, at least to the extent needed later, by choosing all 
parameters carefully. 

\begin{proposition} 
Suppose that $\be$ is a set of admissible Delaunay parameters, with 
$\e < \e_0$. Suppose that the parameters $R_i$ and $a_i$ satisfy
the relationships
\begin{equation}
\sum_{i\neq i_0} R_i^{\frac{N-2}{2}}R_{i_0}^{\frac{N-2}{2}} 
q_i |x_{i_0}-x_i|^{2-N} = q_{i_0}, \quad i_0 = 1, \ldots, n,
\label{eq:4.1}
\end{equation}
and 
\begin{equation}
a_{i_0} = - \frac{1}{q_{i_0}}R_{i_0}^{\frac{N-2}{2}}\sum_{i\neq i_0}
q_i R_i^{\frac{N-2}{2}}|x_{i_0} -x_i|^{-N}(x_{i_0}-x_i),
\label{eq:4.3}
\end{equation}
and that the harmonic function $\bwx$ on $\Omega_{\bar{\rho}}$ is
given by 
\begin{equation}
\bwx = \sum_{i=1}^{n} \frac{\e_i}{2}R_i^{\frac{N-2}{2}}|x-x_i|^{2-N}.
\label{eq:4.defw}
\end{equation}
Then, for $i = 1, \ldots, n$ and $j = 0, \ldots, N$, we have the estimates
\[
\int_{S^{N-1}}(\ui-\bw)(x_i+ \rho_i\theta)\phi_j(\theta)\,d\theta =
O(\e\rho^2), 
\]
and
\[
\int_{S^{N-1}}\left(\frac{\del\ui}{\del r_i}-\frac{\del\bw}{\del r_i}\right)
(x_i+\rho_i\theta)\phi_j (\theta)\,d\theta =O(\e \rho), 
\]
Here $\phi_j(\theta) = \theta_j$, $j = 0, 
\ldots, N$, are eigenfunctions corresponding 
to the first two eigenvalues of the Laplacian on $S^{N-1}$, i.e.
$\lambda_0 = 0$ and $\lambda_1 = N-1$. 
\label{pr:4.1}
\end{proposition}
\begin{remark} The definition (\ref{eq:4.1}) says that $\bq$ is the
eigenvector of some matrix, with coefficients depending on the $R_i$
and $x_i$, with eigenvalue $1$. Thus we first fix $\bq$, and then only
use parameters $R_i$ for which (\ref{eq:4.1}) is valid. Appropriate
choices of $R_i$ do not exist for every $\bq$, but it is not hard to
show that there is a cone in the positive orthant of the $\bq$ space
for which solutions exist. The quantities $\bR$ and $\ba$ will be considered 
fixed and bounded (with each $R_i^{-1}$ also bounded), so we usually
neglect to mention their appearance in various constants in the estimates
below.
\end{remark}
{\bf Proof:}
To prove these estimates, we shall expand both $\uix$ and 
$\bwx$ near $\del B(x_i, \rho_i)$.
\medskip
First we consider $\uix$. From Proposition~\ref{pr:2.3}, if $|a| \rho \le 
1/8$ as we may well assume by taking $\rho$ small, 
then near $\del B(x_i,\rho_i)$ we have 
\[
\uix =\frac{\e_i}{2}\left( R_i^{\frac{N-2}{2}}|x-x_i|^{2-N} + 
R_i^{\frac{2-N}{2}}\big(1 +(N-2)(a_i\cdot (x-x_i)\big)\right)  + O(\e \rho^2).
\]
It is precisely at this step, in estimating the size of the error term, that 
we are led to fix the size of the balls using (\ref{eq:defrho}), since 
with this definition of $\rho$ we have $\e \rho^2 = \e^{\frac{N+2}{N-2}}
\rho^{-N}$. We also have
\[
\del_{r_i} \uix = -\frac{N-2}{2}\e_i \big(R_i^{\frac{N-2}{2}} 
|x-x_i|^{1-N}-R_i^{\frac{2-N}{2}}|x-x_i|^{-1}(a_i\cdot (x-x_i))\big) 
+O (\e \rho).
\]
\medskip
Next, for the expansion of $\bwx$ near $\del B(x_{i_0},\rho_{i_0})$ we have 
\[
\bwx = \frac{\e_{i_0}}{2}R_{i_0}^{\frac{N-2}{2}}|x-x_{i_0}|^{2-N} 
+ \sum_{i\neq i_0} \frac{\e_i}{2}R_i^{\frac{N-2}{2}}\big(|x_{i_0}-x_i|^{2-N} 
+
\]
\[
(2-N) |x_{i_0} -x_i|^{-N}(x_{i_0}-x_i)\cdot(x-x_{i_0})\big) +O(\e \rho^2),
\]
and
\[
\del_{r_{i_0}}\bwx = (2-N)\big(\frac{\e_{i_0}}{2}R_{i_0}^{\frac{N-2}{2}}
|x-x_{i_0}|^{1-N} + 
\]
\[
\sum_{i\neq i_0}\frac{\e_i}{2}R_i^{\frac{N-2}{2}}
|x_{i_0} -x_i|^{-N} |x-x_{i_0}|^{-1} (x_{i_0}-x_{i}) \cdot 
(x-x_{i_0})\big) + O(\e \rho).
\]
\medskip
Taking the difference of these expressions we obtain, 
near $|x-x_{i_0}| = \rho_{i_0}$, that
\[
u_{\e_{i_0}}(R_{i_0},a_{i_0}, x-x_{i_0})-\bwx = \left(\frac{\e_{i_0}}{2}R_{i_0}^{
\frac{2-N}{2}} - \sum_{i\neq i_0}\frac{\e_{i}}{2}R_{i}^{\frac{N-2}{2}}
|x_{i_0}-x_i|^{2-N}\right)
\]
\begin{equation}
+(N-2)\left(\frac{\e_{i_0}}{2} R_{i_0}^{\frac{2-N}{2}} a_{i_0} \\ 
+\sum_{i\neq i_0}
\frac{\e_{i}}{2}R_{i}^{\frac{N-2}{2}} |x_{i_0} -x_i|^{-N}(x_{i_0}-x_i)\right)
\cdot(x-x_{i_0}) + O(\e\rho^2) 
\label{eq:4.diffuw}
\end{equation}
and
\begin{equation}
\del_{r_{i_0}} \left(u_{\e_{i_0}}(R_{i_0}, a_{i_0}, x-x_{i_0})-\bwx\right) = 
(N-2)\bigg(\frac{\e_{i_0}}{2} R_{i_0}^{\frac{2-N}{2}} a_{i_0} +
\label{eq:4.diffderuw}
\end{equation}
\[
\sum_{i\neq i_0}
\frac{\e_{i}}{2}R_{i}^{\frac{N-2}{2}} |x_{i_0} -x_i|^{-N}(x_{i_0}-x_i)\bigg)
\cdot \frac{x-x_{i_0}}{|x-x_{i_0}|} + O(\e \rho).
\]
\smallskip
Now insert (\ref{eq:4.1}) and (\ref{eq:4.3}) into these expressions,
using (\ref{eq:defrho}) again, to obtain the result. 
\hfill $\Box$
\medskip

Finally, using the function $\bwx$ given in Proposition~\ref{pr:4.1},
we define the approximate solution 
\begin{equation}
\buex= \sum_{i=1}^n \chi_i (x-x_i)\uix + \bwx (1- \sum_{i=1}^n 
\chi_i (x-x_i) ).
\label{eq:apprsol}
\end{equation}
Here the $\chi_i$ are smooth, radial functions equalling one
in $|x| \le \rho_i$, vanishing in $|x| \ge 2\rho_i$, and
satisfy $|\del_r \chi_i (x)| \leq c\rho_i^{-1}$ and 
$ |\del^2_r \chi_i (x)| \leq  c \rho_i^{-2}$.  
Notice that we may define the approximate solution $\buex$ regardless of 
whether the relations (\ref{eq:4.1}) and (\ref{eq:4.3}) are satisfied. 

We conclude this section by stating an estimate which will be used 
extensively later, and which follows easily
from the proof of Proposition~\ref{pr:4.1}.

\begin{corollary} Suppose either that $\ba$ satisfies (\ref{eq:4.3})
or $\ba = 0$, and that $\bR$ satisfies (\ref{eq:4.1}). 
Then for $\rho_i \le |x-x_i| \le 1$
\[
\left| \buex^{\frac{4}{N-2}} - \uix^{\frac{4}{N-2}}
\right| \le C \e^{\frac{4}{N-2}} (\rho |x - x_i|^{N-6} + |x - x_i|^{N-5}).
\]
\label{cor:useful}
\end{corollary}
{\bf Proof:} Write $\bar{u}_{\be}$ and $u_{\e_i}$ for the functions
appearing in the estimate. Also, assume $x_i = 0$ for convenience. 
We write the quantity to be estimated as
\[
|(\chi_i u_{\e_i} + (1-\chi_i) \bar{w}_{\be})^{\frac{4}{N-2}}
- u_{\e_i}^{\frac{4}{N-2}}| = |( \bar{w}_{\be} + \chi_i ({u}_{\e_i} - 
\bar{w}_{\be}))^{\frac{4}{N-2}} - u_{\e_i}^{\frac{4}{N-2}}|.
\]
By Taylor's theorem, this is dominated by
\[
C (\bar{w}_{\be})^{\frac{4}{N-2} - 1}|u_{\e_i} - \bar{w}_{\be}|.
\]
From (\ref{eq:4.diffuw}), we see that when $\bR$ solves (\ref{eq:4.1}), 
then $|u_{\e_i} - \bar{w}_{\be}|$ is bounded by $C\e\rho^2$ when $\ba$ solves
(\ref{eq:4.3}), and by $C \e (\rho + |x|)$ when $\ba = 0$. 
Now use that $\bar{w}_{\be} \le C \e |x|^{2-N}$ to complete 
the proof. \hfill $\Box$

\section{The linearized operator on $\R^N \setminus \{0\}$}

Much of the analysis in this paper concerns the linearization
of the nonlinear operator ${\cal N}$ about one of the approximate 
solutions $\buex$. We shall approach the study of this linearization
gradually. In this section we define the linearization about
one of the Delaunay solutions $u_\e(R,a,x)$ on $\R^N \setminus \{0\}$. 
The main point here is the introduction of the Jacobi fields
$\Psi_{\e,R,a}^{j,\pm}(x)$. In succeeding sections we discuss the
refined mapping properties of the linearization, first for the
Dirichlet problem on the unit ball, then for a simpler, locally
radial model for the true linearization on $\R^N \setminus \Sigma$, and 
finally for the true linearization itself. 

Fix one of the Delaunay solutions $u_{\e,R,a}$. The linearization 
of ${\cal N}$ at $u_{\e,R,a}$ is defined by 
\begin{equation}
L_{\e,R,a} v = \left. \frac{d\,}{ds}\right|_{s = 0} {\cal N}(u_\e(R,a,\cdot) +
sv) = \Delta v + \frac{N(N+2)}{4}u_\e(R,a,\cdot)^{\frac{4}{N-2}}v.
\label{eq:deflin}
\end{equation}
More generally, this operator can also be defined as the derivative at $s=0$ 
of ${\cal N}(U(s))$, where $U(s)$ is any one-parameter family of
solutions with $U(0) = u_\e(R,a,x)$, $U'(0) = v$. Viewed this way, it
is immediate that varying the parameters in any one of the families of
Delaunay solutions leads to solutions of $L_{\e,R,a} \Psi = 0$. 
Solutions of this homogeneous problem are called Jacobi fields. 

The Jacobi fields corresponding to these various families are easy
to determine rather explicitly, cf. \cite{MPU1}. We shall be primarily
interested, at least initially, in the Jacobi fields at $a = 0$. 
Since the functions $u_\e(R,0,x)$ are radial, we may separate
variables in this case. Introduce the eigenfunctions 
$\phi_j(\theta)$ and eigenvalues $\lambda_j$ for the Laplacian on $S^{N-1}$. 
Then corresponding to $\lambda_0 = 0$ are the Jacobi fields for the
first and second types of variations:
\begin{itemize}
\item $\Psi_{\e}^{0,-} (x) = \displaystyle{\frac{\del u_\e}{\del \e}}(x)= 
|x|^{\frac{2-N}{2}}\displaystyle{\frac{\del v_\e}{\del \e}}
(-\log|x|)\equiv |x|^{\frac{2-N}{2}}\Phi_\e^{0,-}(-\log |x|)$,
\item $\Psi_{\e}^{0,+} (x) = -|x| \displaystyle{\frac{\del u_\e}{\del r}}(x) 
-\frac{N-2}{2} u_\e(x)= |x|^{\frac{2-N}{2}}\displaystyle{\frac{\del 
v_\e}{\del t}} (-\log|x|)\equiv |x|^{\frac{2-N}{2}}\Phi_\e^{0,+}(-\log |x|)$,
\end{itemize}
Here we take $\phi_0(\theta) = 1$. The third and fourth families
have Jacobi fields corresponding to the next set of eigenvalues,
$\lambda_1 = \ldots = \lambda_N = N-1$. These are:
\begin{itemize}
\item $\Psi_{\e}^{j,+} (x) = |x|^{\frac{4-N}{2}} (\displaystyle{\frac{N-2}{2}}
v_\e(-\log|x|) - \displaystyle{\frac{\del v_\e}{\del t}}(-\log |x|))\,
\phi_j(\theta)\equiv |x|^{\frac{2-N}{2}}\Phi_\e^{1,+}(-\log |x|)\,
\phi_j (\theta)$, 
\item $\Psi_{\e}^{j,-} (x) = \displaystyle{\frac{\del u_\e}{\del x_j}}(x) 
= |x|^{-\frac{N}{2}} (\frac{2-N}{2} v_\e(-\log|x|) - 
\displaystyle{\frac{\del v_\e}{\del t}}(-\log |x|))\,\phi_j(\theta)$
\\ \qquad \qquad
$\equiv |x|^{\frac{2-N}{2}}\Phi_\e^{1,-}(-\log |x|)\,\phi_j (\theta)$.
\end{itemize}

Later we shall also use the Jacobi fields corresponding to differentiating
the family $u_\e(R,a,x)$ and evaluating not necessarily at $R=1$, $a=0$.
We shall denote these Jacobi fields by
\[
\Psi_{\e,R,a}^{j,+}(x) \equiv |x|^{\frac{2-N}{2}}\Phi_{\e,R,a}^{j,+}(x),
\qquad j = 0, \cdots, N.
\]

\section{Function spaces}
We shall consider the action of the linearized operator on weighted
H\"older spaces, which we now define. Since we shall be considering 
both local and global versions of this linearization, we first
define the appropriate spaces on the ball:
\begin{definition}
For each $k \in {\Bbb N}$, $0<\al<1$ and $\sigma \in \R^+$ set 
\[
|u|_{k,\al, [\sigma, 2\sigma]}= \sup_{|x| \in [ \sigma, 2\sigma ]}
\left(\sum_{j=0}^{k} \sigma^j |\nabla^j u(x)|\right) + 
\sigma^{k+\alpha} \sup_{|x|,|y| \in [ \sigma, 2\sigma ]}
\left(\frac{ |\nabla^k u(x) -\nabla^k u (y)|}{|x-y|^\alpha }\right).
\]
Then, for any $\mu\in {\R}$, the space ${\cal C}^{k,\al}_{\mu}
(B(0,R)\setminus \{0\})$ is the collection of functions $u$
which are locally in ${\cal C}^{k,\alpha}$ and for which the norm
\[
||u||_{k,\alpha, \mu}=\sup_{\sigma \leq R/2}\sigma^{-\mu}|u|_{k, \alpha, 
[\sigma, 2\sigma]},
\]
is finite. The subspace ${\cal C}^{k,\al}_{\mu,{\cal D}}(B(0,R)
\setminus \{0\})$
consists of those functions which, in addition, vanish at the boundary 
$|x| = R$.
\label{de:3.1}
\end{definition}
In particular, the function $r^\mu$ is in ${\cal C}^{k,\alpha}_\mu$
for any $k$, $\al$ or $\mu$. 

There are analogous spaces of functions on the complement of the 
prescribed singular set in $S^N$ of the eventual solution defined
in the obvious way. Since it is easier to work in Euclidean coordinates, we 
shall usually consider instead the induced operator and problem on 
$\R^N \setminus \Sigma$, where $\Sigma = \{x_1, \ldots, n\}$. The only 
difference is that we now need to consider decay properties near 
Euclidean infinity. Consideration of these is included in the following 
definition. Divide $\R^N \setminus \Sigma$ into the union of three 
slightly overlapping open regions: $\Omega_1 \equiv \cup_{i=1}^n B(x_i, 1)$, 
a neighbourhood of $\Sigma$, $\Omega_2 \equiv \{|x| > C\}$, a neighbourhood 
of infinity, and a bounded piece $\Omega_3$. 
\begin{definition}
For any $\nu, \nu' \in {\R}$, the space ${\cal C}^{k,\alpha}_{\nu, \nu'} 
({\R}^N \setminus \Sigma)$ is defined to be the collection of $u\in 
{\cal C}^{k,\alpha}({\R}^N\setminus \Sigma)$ for which the norm
\[
||u||_{k,\al, \nu, \nu'}= ||u||_{k,\al,\nu,\Omega_1} + 
\| |x|^{-N+\nu'} u(x/|x|^2)\|_{k, \al , B(0, 1/C)}
+ ||u||_{k,\al,\Omega_3}
\]
is finite. 
\end{definition}
In this definition $|f|_{k, \al, \Omega}$ denotes the usual H\"older norm of the function $f$ over the (bounded) open set $\Omega$.

The one result about these we shall use frequently, and without
comment, is that to check if a function $u$ is an element
of some ${\cal C}^{0,\al}_{\nu}$, say, it is sufficient to
check that $d(x)^{-\nu}|u(x)| \le C$ and $d(x)^{-\nu-1}|\nabla u(x)| 
\le C$. Here $d(x)$ is a smooth positive function on the complement
of $\Sigma$, equaling $|x-x_i|$ on each $B(x_i,1)$. 

\section{The linearized operator on the unit ball}

In this section we shall study the Dirichlet problem in the unit ball 
$B(0,1)$ for the linearization about one of the radial Delaunay 
solutions $u_\e(R,0,x)$, i.e. when $a = 0$. We restrict to the radial 
Delaunay solutions because they may be studied using separation of 
variables.  Later we shall use a perturbation argument to treat the 
more general case. 

For simplicity, denote the linearization $L_{\e,R,0}$ by $L_{\e,R}$,
and when $R = 1$, simply by $L_\e$. Our goal here is to study the problem
\begin{equation}
\left\{ \begin{array}{rrl}      
L_{\e,R} w & = & f \quad \mbox{in}\quad B(0,1)\setminus \{ 0 \}\\[2mm]
     w & = & 0  \quad \mbox{on}\quad \del B(0,1).
        \end{array}
\right.
\label{eq:3.2}
\end{equation}
Specifically, we wish to find suitable spaces on which $L_{\e,R}$ is
surjective, and to find optimal estimates for the solution of (\ref{eq:3.2})
in these spaces, uniformly as $\e$ tends to zero. The first part of
this was already considered in \cite{MPU1}. There it was proved that
$L_{\e,R}$ is surjective as a map between certain weighted Sobolev
spaces, but the uniformity in $\e$ was not considered. For various
technical reasons, it is more useful here to consider the
action of $L_{\e,R}$ on the family of weighted H\"older spaces 
defined in the previous section.

One can prove, analogous to the result of \cite{MPU1}, that 
\begin{equation}
L_{\e,R}: {\cal C}^{2,\al}_{\mu,{\cal D}}(B(0,1)\setminus \{0\}) 
\longrightarrow {\cal C}^{0,\al}_{\mu-2} (B(0,1)\setminus \{0\})
\label{eq:21.1}
\end{equation}
is Fredholm for all $\mu \in \R$, $\mu \notin \{\mu_j^\pm (\e) \}$, where
$\{\mu_j^\pm (\e) \}$ is a discrete set, with each $\mu_j^\pm (\e)$ determined 
by the ordinary differential equation induced by $L_{\e,R}$ on the
$j^{\mbox{th}}$ eigenspace of the Laplacian on $S^{N-1}$. The values
of the inadmissible weights $\mu_j^\pm (\e)$ may depend on $\e$ (although, 
somewhat miraculously, it turns out that those with $j = 0, \ldots, N$ do not).
These values play the same role as the the indicial roots of a regular 
singular problem; indeed, they arise in the same manner as the rates of 
growth or decay of solutions of $L_{\e,R}w = 0$. However, they
are determined in a far less explicit manner, analytically rather
than algebraically, and so are usually impossible to compute. 
Fortunately, we know the solutions of this equation, hence the
rates of growth of these solutions, explicitly when $j = 0, \ldots, N$, 
so we can determine, as in \cite{MPU1}, that 
\begin{equation}
\Psi^{0,+}_{\e,R,a} \in {\cal C}^{2,\al}_{\mu} \quad \mbox{for all }
\mu \le \frac{2-N}{2}, \qquad \Psi^{0,-}_{\e,R,a} \in {\cal C}^{2,\al}_{
\mu} \quad \mbox{for all } \mu < \frac{2-N}{2},
\label{eq:3.wt0}
\end{equation}
and for $j = 1, \ldots, N$, 
\begin{equation}
\Psi^{j,+}_{\e,R,a} \in {\cal C}^{2,\al}_{\mu} \quad \mbox{for all }
\mu \le \frac{4-N}{2}, \qquad \Psi^{j,-}_{\e,R,a} \in {\cal C}^{2,\al}_{
\mu} \quad \mbox{for all } \mu \le \frac{-N}{2}.
\label{eq:3.wtj}
\end{equation}
From this one sees that $\mu_0^\pm = (2-N)/2$ and $\mu_j^{+} = (4-N)/2$, 
$\mu_j^- = -N/2$, $j = 1, \ldots, N$. As we pointed out above, the values
$\mu_j^{\pm}$, for $j=0, \ldots , N$, do not depend on $\e$.

It turns out that (\ref{eq:21.1}) is surjective when $\mu < (2-N)/2$, 
$\mu \neq \mu_j^-$, and injective when $\mu > (2-N)/2$. The basic
conflict is that, although this indicates that we should be working on spaces 
with weight less than $(2-N)/2$, we are using $L_{\e,R}$ to correct
the approximate solution, which only grows like $r^{(2-N)/2}$, and we
obviously do not want the correction term to blow up faster than this. 
As proved in \cite{MPU1}, this may be remedied by working in a finite
dimensional extension of ${\cal C}^{2,\al}_{\mu}$, where $(2-N)/2 < \mu 
< (4-N)/2$ and the extension is the span of the Jacobi fields 
$\Psi_{\e,R}^{0,\pm}$. The point is that although in general (\ref{eq:3.2}) 
has no solution $w \in 
{\cal C}^{2,\al}_{\mu}$ when $f \in {\cal C}^{0,\al}_{\mu-2}$ and
$\mu > (2-N)/2$, there {\it is} a solution with the correct decay but with
eigencomponent $w_0$ having the wrong Dirichlet data. To correct this,
one may add on some multiple of one or the other of the Jacobi fields 
$\Psi^{0,\pm}_{\e,R}$. We shall prove something similar, but with two 
additional considerations. The first is that unlike the analogous situation 
in \cite{MP}, any choice of right inverse has norm blowing up as $\e$ tends to 
zero. However, this blow-up happens only on the subspace of functions spanned 
by the eigenfunctions $\phi_j(\theta)$, $j = 0, \ldots, N$. In order to obtain 
reasonable estimates, we need to work in the space with weight $\mu$ with $1 
< \mu <2$. For this reason we shall also need to append to this `deficiency 
space' the Jacobi fields with index $j = 1, \ldots N$. The other consideration 
is that we wish to prove surjectivity using only the Jacobi fields $\Psi_{\e,
R}^{j,+}$, because the others are considerably more difficult to use in the 
nonlinear analysis. In fact, using the Jacobi fields $\Psi_{\e,R}^{j,-}$
would force us to reposition the singularities or change the Delaunay 
parameters, which for technical reasons it is far better to avoid. On
the level of linear analysis this is a good choice as well, for with
this restricted domain the linearization is an isomorphism rather
than just surjective. 

We come to the main result of this section. 

\begin{proposition}
Fix $R \in \R^+$, $R\neq 1$, and $\mu \in ((4-N)/2,2)$ and define the deficiency space 
\begin{equation}
W \equiv \mbox{span\,}\{\Psi_{\e,R}^{j,+}: j = 0, \ldots, N\}.
\label{eq:3.defsp}
\end{equation}
Then there exists an $\e_0 > 0$ such that for all $\e \in (0, \e_0]$, the operator
\begin{equation}
L_{\e,R}: {\cal C}^{2,\al}_{\mu,{\cal D}}(B(0,1)\setminus \{0\}) \oplus W
\longrightarrow {\cal C}^{0,\al}_{\mu-2}(B(0,1)\setminus \{0\})
\label{eq:3.defdefop}
\end{equation}
is an isomorphism. The inverse map will be denoted $L_{\e,R}^{-1}$. 
In particular, 
for $f \in {\cal C}^{0,\al}_{\mu-2} (B(0,1)\setminus\{0\})$ there exists a 
unique solution $w$ of (\ref{eq:3.2}) which has a decomposition 
\[
w(x) = G_{\e,R} (f)(x) + \frac{1}{\e} \sum_{j=0}^N K^j_{\e,R}(f)
\Psi^{j,+}_{\e}(x/R).
\]
The operator 
\[
G_{\e,R}:{\cal C}^{0,\alpha}_{\mu -2} (B(0,1)\setminus \{ 0\})\rightarrow 
{\cal C}^{2,\al}_{\mu} (B(0,1)\setminus \{ 0\})
\]
and functionals
\[
K^j_{\e,R}:{\cal C}^{0,\alpha}_{\mu -2}(B(0,1)\setminus \{ 0\}) \rightarrow \R
\]
 satisfy the following properties: 
\begin{itemize}
\item  $G_{\e,R}$ is defined for all $\mu \in ((4-N)/2, 2)$ and is bounded 
independently of $\e$ when $1 < \mu < 2$. 

\item The restriction of $G_{\e,R}$ to the space of functions $h(r,\theta)$
with eigencomponents $h_j(r)$ vanishing for $j = 0, \ldots, N$, is defined and
bounded independently of $\e$ when $-N < \mu < 2$.

\item  The functional $K^0_{\e,R}$ is defined for all $\mu > (2-N)/2$ and is bounded independently of $\e$ when $0< \mu $.

\item The functionals $K^j_{\e,R}$, $j = 1, \ldots, N$, are defined for all $\mu > (4-N)/2$ and are bounded independently of 
$\e$ when $1< \mu $. 
\end{itemize}
\label{pr:3.1}
\end{proposition}

The proof of Proposition~\ref{pr:3.1} will be broken into several steps. 
The solution is constructed by separation of variables, restricting the 
problem to each eigenspace of the Laplacian on the sphere. First, we note 
that although the parameter $R$ is important later, at this stage it is irrelevant
and may be scaled away by setting $y = x/R$. In fact, replacing $f(x)$ by 
$R^{-2}f(x/R)$ and $w(x)$ by $w(x/R)$, we see that (\ref{eq:3.2}) is equivalent to 
\begin{equation}
\left\{ \begin{array}{rrl}      
L_\e w & = & f \quad \mbox{in}\quad B(0,1/R)\setminus \{ 0 \}\\[2mm]
    w & = & 0  \quad \mbox{on}\quad \del B(0,1/R).
        \end{array}
\right.
\label{eq:3.3}
\end{equation}
When we change to the variable $t = -\log |x|$, this rescaling has the effect of 
replacing $v_{\e,R}(t) = v_\e(t + \log R)$ by $v_\e(t)$. The $K^j_{\e,R}$ will 
be denoted simply $K^j_\e$, and so on. In the $(t,\theta)$ variables, the 
H\"older spaces are converted to the ordinary (translation-invariant)
spaces, with the weight function $e^{t\mu}$. 

Write the eigenfunction decompositions of $w$ and $f$ as
\[
w(x) = |x|^{-\frac{N-2}{2}}\sum_{j=0}^{\infty} w_j(-\log |x|)\phi_j(\theta),
\qquad \mbox{and} \quad 
f(x) = |x|^{-\frac{N+2}{2}}\sum_{j=0}^{\infty}f_j (-\log |x|)\phi_j(\theta).
\]
If $w$ solves (\ref{eq:3.3}) then $w_j(t)$ solves
\[
\left\{ 
      \begin{array}{rll}
\LLej w_j (t) & = & \displaystyle{\frac{d^2w_j}{dt^2}} -
(\frac{N-2}{2})^2 w_j -\lambda_j w_j + \frac{N(N+2)}{4}v_\e^{
\frac{4}{N-2}}w_j =f_j\\[2mm]
w_j (\log R) & = & 0 .       
      \end{array}
\right.
\]
The $\LLej$ are the restrictions to the eigenspaces of $\Delta_{S^{N-1}}$
of the operator
\[
\LLe w (t, \theta) = \frac{d^2w}{dt^2} -(\frac{N-2}{2})^2 w +
\Delta_{S^{N-1}} w + \frac{N(N+2)}{4}v_\e^{\frac{4}{N-2}}w,
\]
which is simply problem (\ref{eq:2.1}) transformed from
$\R^N \setminus \{0\}$ to the cylinder $\R \times S^{N-1}$.
\medskip

We may assume, after multiplying by a suitable factor, that 
\[
||f||_{0, \alpha , \mu -2} =1.
\]
We first study the operators $\LLej$ for $j > N$, then for $j =0$ and finally
for $j = 1, \ldots, N$. 

\medskip

{\bf Step 1.} Suppose that $-N < \mu < 2$. Let $\baf$ be the projection
of $f$ onto the sum of the eigenspaces of $\Delta_{S^{N-1}}$ with $j > N$,
\[
\baf(t,\theta)= \sum_{j\geq N+1} f_j(t )\phi_j(\theta ).
\]
We must solve 
\begin{equation}
\left\{
       \begin{array}{rll}
\LLe \bar{w} (t, \theta) & = & \baf(t, \theta) \quad \mbox{in} 
\quad (\log R, +\infty )\times S^{N-1}\\[2mm] 
\bar{w} (\log R, \theta) & = & 0.  
       \end{array}
\right.
\label{eq:3.5}
\end{equation}
To do this we first show that for each $T> \log R$ there exists a unique
solution of
\[
\left\{
       \begin{array}{rll}
\LLe \bar{w}_T (t, \theta) & = & \baf \quad \mbox{in} \quad (\log R , T)
\times S^{N-1}\\[2mm] 
\bar{w}_T (\log R, \theta) & = & 0 \\[2mm] 
\bar{w}_T (T, \theta) & = & 0.
       \end{array}
\right.
%\label{eq:3.6}
\]
This solution can be obtained variationally. Consider the energy
function ${\cal E}_T$
\[
{\cal E}_T(w)=\int_{\log R}^T \int_{S^{N-1}} \left(\dot{w}^2 +\frac{(N-2)^2}{4}
w^2 +|\nabla_{S^{N-1}} w|^2 - \frac{N(N+2)}{4}v_\e^{\frac{4}{N-2}}w^2 + 
\baf w\right) dt d\theta .
\]
Since $\lambda_j\geq\lambda_{N+1}=2N$, we estimate
\[
{\cal E}_T(w)\geq \int_{\log R}^T \int_{S^{N-1}} \left(\dot{w}^2 + \left(
\frac{(N+2)^2}{4} - \frac{N(N+2)}{4}v_\e^{\frac{4}{N-2}}\right) 
w^2 + \baf w \right) dt d\theta.
\]
Since $0< v_\e<1$, this energy is convex and proper, and the existence of 
a unique minimizer for ${\cal E}_T$ is immediate. 

\begin{lemma}
For $-N < \mu < 2$, let $\delta =\frac{N-2}{2} +\mu$. Then, for all $R >0$, there exist
constants $\e_0 >0$ and $C >0$ independent of $T$, such that for all 
$ \e \in (0,\e_0]$, we have
\[
\sup_{\theta \in S^{N-1}}\sup_{t\in [\log R,T]}
e^{\delta t} | \bar{w}_T (t, \theta) | \leq C. 
\]
\label{le:3.2}
\end{lemma}
{\bf Proof :}
We prove this by contradiction. Recall that, by assumption, 
\[
\sup_{\theta \in S^{N-1}}\sup_{t\in [\log R,T]}e^{\delta t} |\baf (t,\theta)|
\leq 1.
\]
If the assertion were not true, then there would exist 
sequences of functions $\baf_i$, numbers $T_i$ and Delaunay 
parameters $\e_i$ and a sequence of solutions $\bar{w}_{T_i}$ 
such that 
\[
\lim_{i\rightarrow \infty}(\sup_{\theta \in S^{N-1}}\sup_{\log R
\le t \le T_i} e^{\delta t}|\bar{w}_{T_i}(t, \theta)|)=\infty.
\]
We let $v_i = v_{\e_i}$.  Now choose $t_i\in (\log R, T_i)$ such that
\[ 
\sup_{\theta \in S^{N-1}}e^{\delta t_i} | \bar{w}_{T_i} 
(t_i, \theta)| = \sup_{\theta \in S^{N-1}}\sup_{\log R\le t \le T_i}
e^{\delta t} |\bar{w}_{T_i}(t, \theta)| \equiv A_i,
\]
and define 
\[
\tw_i (t, \theta) = A_i^{-1}e^{\delta t_i} \bar{w}_{T_i}(t+t_i, \theta),
\]
\[
\tilde{f}_i(t,\theta) = A_i^{-1}e^{\delta t_i} \baf(t+t_i,\theta).
\]
Then, by definition,
\[
\sup_{\theta \in S^{N-1}}\sup_{\log R-t_i \le t \le T_i-t_i}e^{\delta t}|
{\tw}_i (t, \theta)|= 1,
\]
and $\tilde{f}_i$ tends to $0$ in norm.
In addition, 
\[
\frac{d^2{\tw}_i}{dt^2} -\frac{(N-2)^2}{4}{\tw}_i +\Delta_{S^{N-1}} 
{\tw}_i + \frac{N(N+2)}{4}v_i^{\frac{4}{N-2}}(t+t_i){\tw}_i = \tilde{f}_i,
\]
on $[\log R-t_i, T_i-t_i]\times S^{N-1}$. 
Passing to a subsequence, we may assume that $\log R -t_i$ converges to 
some number $\tau_1 \in {\R^-}\cup \{-\infty \}$ and also that $T_i-t_i$ 
converges to $\tau_2 \in {\R^+}\cup \{+\infty \}$. Furthermore,
we may also assume that over every compact set of $(\tau_1,\tau_2)$, the sequence
$v_i(t+t_i)$ converges to $v_{\infty}(t)$ and ${\tw}_i$ converges to ${\tw}$.
By construction, $\tw \neq 0$. Then $v_{\infty}$ is a solution of 
(\ref{eq:2.3}) and ${\tw}$ satisfies the equation
\begin{equation}
\frac{d^2{\tw}}{dt^2} -\frac{(N-2)^2}{4}{\tw}+\Delta_{S^{N-1}}{\tw}+ 
\frac{N(N+2)}{4}v_{\infty}^{\frac{4}{N-2}}{\tw} = 0
\label{eq:3.8}
\end{equation}
on $[\tau_1,\tau_2] \times S^{N-1}$. Let $H_\infty = H(v_\infty, 
\dot{v}_\infty)$. There are a few cases to consider, depending
on the values of $\tau_1$, $\tau_2$ and $H_\infty$. 
To analyze these, we require some growth estimates.
\begin{lemma}
For every $\eta>0$, there exists $\e_0 >0$ such that, for all $j \ge N+1$ and 
for all $\e \in (0, \e_0]$, every solution of  
\begin{equation}
\frac{d^2 w}{dt^2} -\left(\frac{(N-2)^2}{4} +\lambda_j - 
\frac{N(N+2)}{4}v_\e^{\frac{4}{N-2}}\right) w = 0
\label{eq:3.9}
\end{equation}
either decays to $0$ faster than $e^{-(\frac{(N+2)^2}{4} -2\eta)^{1/2}t}$  
at $\infty$ (respectively, $e^{(\frac{(N+2)^2}{4} -2\eta)^{1/2}t}$ at 
$-\infty$) or blows up faster than  $e^{(\frac{(N+2)^2}{4} -2\eta)^{1/2}t}$ 
at $\infty$ (respectively, $e^{-(\frac{(N+2)^2}{4} -2\eta)^{1/2}t}$ at $-\infty$).
\label{le:3.3}
\end{lemma}
{\bf Proof :}
Because $j \ge N+1$, the term of order zero, $\frac{(N-2)^2}{4} +
\lambda_j - \frac{N(N+2)}{4}v_\e^{\frac{4}{N-2}}$ is positive, hence
$\LLej$ satisfies the maximum principle. Let 
\[
\delta_j = \left(\frac{(N-2)^2}{4} +\lambda_j - \frac{N(N+2)}{4}\right)^{1/2}
\qquad \mbox{and} \quad \beta(\eta) = \left(\frac{(N+2)^2}{4} - \eta\right)^{1/2}.
\]
We shall show the existence of a solution bounded above by 
$e^{-\beta(2\eta)t}$ as $t$ tends to $\infty$ and bounded below
by $e^{-\beta(2\eta)t}$ as $t$ tends to $-\infty$. 
Replacing $t$ by $-t$ yields a solution with the appropriate
exponential decay at $-\infty$ and blowup at $\infty$. Since these
span the space of all solutions, this will prove the lemma. 

Since $e^{-\delta_j t}$ satisfies $\LLej e^{-\delta_j t} \le 0$, 
it is elementary that there exists
a function $y_j(t)$ such that $y_j(0) = 1$ and 
\[
y_j(t) \geq e^{-\delta_j t} \quad \mbox{for} \quad t\leq 0, \qquad 
0 \leq y_j(t) \leq e^{-\delta_j t} \quad \mbox{for} \quad t\geq 0.
\]
Next, note that 
\begin{equation}
y_j(t + T_\e) = \frac{y_j(T_\e /2)}{y_j(-T_\e /2)} y_j(t) \qquad
\mbox{for all}\ t. 
\label{eq:30.1}
\end{equation}
To see this, simply observe that $w \equiv y_j(T_\e/2) y_j(t) -
y_j(t+T_\e)y_j(-T_\e/2)$ solves (\ref{eq:3.9}), decays exponentially 
at $\infty$ and takes the value $0$ at $-T_\e/2$, hence by the
maximum principle must vanish identically.  We next claim that 
\begin{equation}
0 <  y_j(t) \le e^{-\beta(2\eta)t} \quad \mbox{for} \quad t \in [0,T_\e/2],
\qquad \mbox {and}\ y_j(t) \ge e^{-\beta(2\eta)t} \quad \mbox{for} \quad 
t \in [-T_\e/2,0].
\label{eq:30.2}
\end{equation}
Granting this for the moment, then 
\[
\frac{y_j(T_\e/2)}{y_j(-T_\e/2)} \le \frac{e^{-\beta(2\eta)T_\e/2}}{
e^{\beta(2\eta)T_\e/2}} = e^{-\beta(2\eta)T_\e}.
\]
Using this, along with (\ref{eq:30.1}) and (\ref{eq:30.2}),
the desired estimates for $y_j(t)$ for all $t$ are immediate. 
Hence, to prove the lemma, it will suffice to examine the function 
$y_j$ on the interval $[-T_\e /2, T_\e /2]$.

To prove the claim, we use the fact that 
\begin{equation}
\lim_{\e \rightarrow 0} \frac{d y_j}{dt}(0)=-\gamma_j \equiv -
\left(\frac{(N-2)^2}{4} +\lambda_j\right)^{1/2}.
\label{eq:3.11}
\end{equation}
To see this, multiply (\ref{eq:3.9}) by $e^{-\gamma_j t}$ and integrate 
from $0$ to $+\infty$. Integrating by parts and using that $y_j(0)=1$, we get
\[
\frac{d y_j}{dt}(0) +\gamma_j = \int_0^{+\infty}
 \frac{N(N+2)}{4}v_\e^{\frac{4}{N-2}}y_j(t) e^{-\gamma_j t}\,dt.
\]
Since the right side is nonnegative, we first get that $y_j'(0) \ge -\gamma_j$.
Next, from Proposition~\ref{pr:2.2}, we can bound
$v_\e(t) \leq 2\e e^{\frac{N-2}{2} t}$ for $t \in [0,T_\e/2]$
and $v_\e(t)\leq 1$ for $t \geq T_\e /2$. Hence 
\[
\int_0^{\infty}v_\e^{\frac{4}{N-2}}y_j(t) e^{-\gamma_j t}\,dt 
\leq \int_0^{T_\e /2} (\e e^{\frac{N-2}{2}t})^{\frac{4}{N-2}} 
e^{-(\delta_j+\gamma_j)t}\,dt + \int_{T_\e /2}^{+\infty} 
e^{-(\delta_j+\gamma_j) t}\,dt.
\]
Since $T_\e \sim -\frac{4}{N-2}\log \e$, the right hand side of this 
expression tends to $0$ as $\e$ tends to $0$, and the proof
of the claim is complete.

We finally prove the sharper estimates (\ref{eq:30.2}) on $[-T_\e/2,T_\e/2]$.
From Proposition~\ref{pr:2.2} again, there exist $\e_0 >0$, 
$t_0 >0$ such that for all $\e \in (0, \e_0]$ and $t \in [-T_\e/2 + t_0, 
T_\e /2 -t_0]$ 
\begin{equation}
\frac{(N-2)^2}{4} + \lambda_j - \frac{N(N+2)}{4}v_\e^{\frac{4}{N-2}}(t) \geq 
\frac{(N+2)^2}{4} -\eta = \beta(\eta)^2.
\label{eq:30.3}
\end{equation}
Using (\ref{eq:3.11}) we may also assume that 
$y_j'(0) \leq -\beta(\eta)$  
for all $\e \in (0, \e_0]$. By (\ref{eq:30.3}), we may use the
maximum principle on the interval $[0, T_\e /2 -t_0]$ to conclude that  
\[
0\leq y_j (t) \leq e^{- \beta(\eta)t},\quad 0 \le t\le T_\e/2 -t_0, \qquad 
\mbox{and} \quad y_j (t)\geq e^{-\beta(\eta)t}, \quad - T_\e /2 +t_0 \le t
\le  0.
\]
Furthermore, $y_j$ is monotone decreasing, so its maximum on any interval
is attained at the left endpoint and its minimum at the right endpoint.
Hence, for $\e$ sufficiently small, 
\[
e^{-\beta(\eta)(\frac{T_\e}{2} -t_0)} \leq e^{-\beta(2\eta)\frac{T_\e}{2}}.
\]
This completes the proof of (\ref{eq:30.2}) and hence the lemma. 
\hfill $\Box$

We now return to our proof of Lemma~\ref{le:3.2}.
\begin{itemize}
\item First consider the case where $H_{\infty}\neq 0$. 
Decompose ${\tw}$ as 
\[
{\tw}(t, \theta) = \sum_{j\geq N+1} {\tw}_j(t)\phi_j(\theta ).
\]
If $\tau_1 =-\infty$ then, taking $2\eta \leq \displaystyle{\frac{(N+2)^2}{4}}
-\delta^2$ in Lemma~\ref{le:3.3} and using the fact that
${\tw}$ is bounded as $t \rightarrow -\infty$ by $e^{-\delta t}$ with 
$\delta < \displaystyle{\frac{N+2}{2}}$, we see that for $\e$ sufficiently 
small, ${\tw}_j$ decays exponentially as $t$ goes to $-\infty$. 
Similarly, if $\tau_2 =+\infty$, then by the same argument
$\tw_j$ decays exponentially as $t \rightarrow \infty$. 
If $\tau_1 > -\infty$, then $\tw_j (\tau_1)=0$ and if $\tau_2 < 
\infty$, then $\tw_j (\tau_2)=0$.

Therefore, we can multiply the equation satisfied by ${\tw}_j$ by 
${\tw}_j$ itself and integrate by parts. There are no boundary
terms because of the exponential decay and vanishing boundary
values, so we obtain
\[
0= \int_{\tau_1}^{\tau_2} (\dot{{\tw}}_j^2 +\frac{(N-2)^2}{4} 
{\tw}_j^2 +\lambda_j {\tw}_j^2 - \frac{N(N+2)}{4}v_{\infty}^{\frac{4}{N-2}}
{\tw}_j^2 ) dt.
\]
But since $\lambda_j \geq 2N$ and $0 < v_\infty < 1$, we conclude
that $\tw \equiv 0$.

\item Next consider the case where $H_{\infty}= 0$.  This case
is almost identical to the last one, except that the exponential
decay is simpler to obtain. Decompose ${\tw}$ as before. If 
$\tau_1 =-\infty$ then by 
Proposition~\ref{pr:2.2} and (\ref{eq:3.8}) ${\tw}_j$ grows asymptotically
like $e^{\pm\gamma_j t}$, where $\gamma_j^2= \frac{(N-2)^2}{4} +\lambda_j$,
near $-\infty$. But $\tw$ is bounded by $e^{-\delta t}$ near 
$-\infty$, where $\delta < \displaystyle{\frac{N+2}{2}} \leq \gamma_j$, and so ${\tw}_j$ decays like 
$e^{\gamma_j t}$ at $-\infty$. Similarly, if $\tau_2=+\infty$,
$\tw_j$ decays exponentially as $t$ goes to $+\infty$, and as
before, $\tw_j (\tau_1)=0$ and $\tw_j (\tau_2)=0$ if either
$\tau_1$ or $\tau_2$ are finite.

The same integration by parts as before implies that $\tw \equiv 0$.
\end{itemize}
In either case, we have produced a contradiction, since by assumption 
$\tw$ is non trivial. The proof of Lemma~\ref{le:3.2} is complete.
\hfill $\Box$

Finally, we let $T\rightarrow +\infty$ to get the existence of a 
unique solution $\bar{w}$ to (\ref{eq:3.5}) which is uniformly
bounded by $Ce^{-\delta t}$ for $t \ge 0$. 

\medskip 
The second and third steps of the proof of Proposition~\ref{pr:3.1}
are rather similar to Step 1 above. The main difference is
that an extra step is needed to ensure that $w$ has vanishing
boundary value at $t=\log R$. This is where the nonuniformity
of $L_{\e,R}^{-1}$ as $\e \rightarrow 0$ arises.

\medskip

{\bf Step 2.} Let $j = 0$ and assume that $\mu > 0$. We shall now solve 
the equation
\begin{equation}
\left\{
       \begin{array}{rll}
{\Bbb L}_{\e,0}\,  w_0 (t) & = & f_0 \quad \mbox{in} \quad (\log R, \infty)\\[2mm] 
                w_0 (\log R) & = & 0. 
       \end{array}
\right.
\label{eq:3.12}
\end{equation}
Normalize $f_0$ so that $||f_0||_{0,\al,\mu-2} = 1$, and choose an extension 
$\tilde{f}_0$ of $f_0$ to $\R$ with the property that 
$||\tilde{f}_0||_{0,\al,\mu-2} \le 2$. For each $T> \log R$, 
let $w_T$ be the unique solution of 
\[
{\Bbb L}_{\e, 0} w_T = \tilde{f}_0,\qquad  w_T(0) = 0, \qquad \dot{w}_T(0) = 0.
\]
\begin{lemma}
For $\mu >0$, let $\delta =\displaystyle{\frac{N-2}{2}} +\mu$. Then for all $R >0$, there 
are constants $\e_0 >0$ and $C >0$ independent of $T$, such that for all $\e 
\in (0,\e_0]$ we have 
\[
\sup_{t\in [\log R,T]}e^{\delta t} | \bar{w}_T (t) | \leq C.
\]
\label{le:3.4}
\end{lemma}
{\bf Proof :}
First, the two independent solutions of ${\Bbb L}_{\e,0} w=0$ are
the functions $\Phi_\e^{0, \pm} = r^{\frac{N-2}{2}}\Psi_{\e}^{0, \pm}$, 
and neither these, nor any linear combinations of them decay or grow 
exponentially. In fact, $\Phi_{\e}^{0,+}$ is periodic (of period $T_\e$), 
while $\Phi_{\e}^{0,-}$ grows linearly. This last fact may be seen by 
differentiating the equality $v_\e(t +T_\e) =v_\e(t)$ with respect to $\e$ to get
\[
\Phi_{\e}^{0,-} (t+T_\e) +\frac{d T_\e}{d\e} \Phi_{\e}^{0,+} (t+T_\e)= 
\Phi_{\e}^{0,-}(t);
\]
this gives the linear growth rate since the derivative of $T_\e$ never 
vanishes. (Even without knowing that $T_\e$ is never
stationary, it is elementary that no solutions of ${\Bbb L}_{\e,0}$
can decay at infinity.)

The remainder of the proof is nearly identical to that of Lemma~\ref{le:3.2}.
If the claim were not true, then there would exist sequences
$f_{0,i}$, $T_i$, $\e_i$ and $\bar{w}_{T_i}$ such that
$A_i \equiv \sup_{t\in (-\infty, T_i]}e^{\delta t}|\bar{w}_{T_i}|$ tends
to infinity. Notice that $\bar{w}_{T_i}$ grows at most linearly as $t$ tends to 
$-\infty$, hence the supremum in the last fomula is always achieved. 
Choosing $t_i\in (\log R, T_i)$ to maximize $e^{\delta t}|\bar{w}_{T_i}(t)|$ 
(so that $e^{\delta t_i}|\bar{w}_{T_i}(t_i)| = A_i$), we rescale the functions 
and translate the independent variable by $t_i$ to obtain a solution of
\[
\frac{d^2{\tw}_i}{dt^2} -\frac{(N-2)^2}{4}{\tw}_i + \frac{N(N+2)}{4}
v_i^{\frac{4}{N-2}}{\tw}_i = \tilde{f}_{0,i},
\]
in $[-\infty,T_i - t_i]$. 
This solution satisfies  
\[
\sup_{t\in [-\infty,T_i-t_i]}e^{\delta t} |{\tw}_i (t) | = 1,
\]
while $\tilde{f}_{0,i}$ tends to zero in norm.

Passing to a subsequence, we obtain a nontrivial solution ${\tw}$ of 
the equation
\begin{equation}
\frac{d^2{\tw}}{dt^2} -\frac{(N-2)^2}{4}{\tw} + \frac{N(N+2)}{4}
v_{\infty}^{\frac{4}{N-2}}{\tw} = 0,
\label{eq:3.15}
\end{equation}
over $(-\infty,\tau_2]$. 
\begin{itemize}
\item We can immediately see that $\tau_2 =\infty$, for if not
then $\tw$ would be a non trivial solution of (\ref{eq:3.15}) 
such that  ${\tw}(\tau_2)=\dot{\tw}(\tau_2)=0$, which would
imply that $\tw \equiv 0$. 

\item We can also easily rule out the case where $H(v_\infty, \dot{v}_\infty) 
\equiv H_\infty \neq 0$. For if ${\tw}$ were a nontrivial 
solution of (\ref{eq:3.15}) then it would decay exponentially at $\infty$. 
But we have already seen that any solutions of ${\Bbb L}_{\e,0} w=0$ decays 
exponentially.

\item  The final case, where $H_{\infty}=0$, is also easy to rule out.
By Proposition~\ref{pr:2.2} and equation (\ref{eq:3.15}), $\tw$ 
grows asymptotically near $\infty$ like $e^{\pm \frac{N-2}{2} t}$. 
Since ${\tw}$ is bounded by $e^{- \delta t}$ and because 
$\delta > \frac{N-2}{2}$, we may conclude that $\tw \equiv 0$.
\end{itemize}

With the proof of Lemma~\ref{le:3.4} complete, we let $T \rightarrow \infty$
and obtain a unique solution $\bar{w}$ to (\ref{eq:3.12}) which satisfies 
\[
\sup_{t\in [\log R,+\infty)}e^{\delta t} | \bar{w} (t) | \leq C.
\]
There is no reason why this solution should satisfy the boundary condition
at $t=\log R$, and so we must add an additional term to correct the 
boundary data. We define the solution $w_0$ by
\[
w_0(t) = {\bar{w}}(t) - {\bar{w}}(\log R)(\Phi_{\e}^{0,+}(\log R))^{-1}
\Phi_{\e}^{0,+}(t).
\]
This obviously satisfies the equation and the correct boundary
conditions, so the proof will be complete if we show that 
$\frac{1}{\e}\Phi_{\e}^{0,+}(\log R)$ is bounded from below by some constant 
independent of $\e$. The following result is deduced easily from 
Proposition~\ref{pr:2.3}.
\begin{lemma}
Given $R>0$, $R\neq 1$, there exists a constant $m_0 > 0$ such that for all 
$\e \in (0, \e_0]$, we have 
\[
\Phi_{\e}^{0,+}(\log R) \geq \e m_0 .
\]
\label{le:3.5}
\end{lemma}

\begin{remark}
The definition of $w_0$ is certainly not unique, for we could
have added any appropriate combination of $\Phi_\e^{0,-}$
and $\Phi_\e^{0,+}$ to correct the boundary data. However, $\Phi_\e^{0,-}$ is
more difficult to use in the nonlinear analysis, so it is important that 
we only use $\Phi_\e^{0,+}$ here. Moreover, the use of $\Phi_\e^{0,-}$ would oblige us to change the Delaunay parameters of our approximate solution in order to get a solution of the nonlinear problem.
\end{remark}

{\bf Step 3.} Finally, suppose that $1<\mu$. We wish to solve, for 
$j=1, \ldots , N$, the problem 
\begin{equation}
\left\{
       \begin{array}{rll}
\LLej w_j (t) & = & f_j\quad \mbox{in} \quad (\log R, \infty)\\[2mm] 
     w_j (\log R) & = & 0 .
       \end{array}
\right.
\label{eq:3.16}
\end{equation}
This is done, {\it mutatis mutandis}, following the proofs in Steps
1 and 2. We describe the minor changes that need to be made.
First, define an extension $\tilde{f}_j$ of $f_j$ to all of $\R$, and 
find the unique solution $\bar{w}_T$ of 
\[
\LLej \bar{w}_T= \tilde{f}_j, \qquad \bar{w}_T(T) = \dot{\bar{w}}_T(T) = 0.
\]
The uniform bound of $e^{\delta t} |\bar{w}_T|$ on $(-\infty,T]$ is proved 
as before by contradiction. The fact that that the supremum is achieved 
follows as before from the fact that, for $t \leq \log R$, $\bar{w}_T$ 
is a linear combination of $\phi_\e^{1,\pm}$ and thus blows up as $t$ 
tends to $-\infty$ at most like a constant times $e^{-t}$. In this step 
we obtain, by rescaling and translation, a nontrivial solution $\tw$ of
\[
\frac{d^2{\tw}}{dt^2} -\frac{(N-2)^2}{4}{\tw} -(N-1){\tw} + 
\frac{N(N+2)}{4}v_{\infty}^{\frac{4}{N-2}}{\tw} = 0.
\]
To rule out the existence of a nontrivial solution to this equation,
we reason as before. In the second case, where $H_\infty \neq 0$,
we use the two explicit solutions of ${\Bbb L}_{\e,j} w=0$ given by
$\Phi_{\e}^{j,\pm}(t)$. These solutions do grow at an exponential
rate, but by hypothesis, both either blow up or do not decay 
quickly enough, and we conclude as above that this case never occurs.
When $H_\infty =0$ we note that the asymptotic behaviour of solutions 
near $\pm \infty$ is now given by $e^{\pm \frac{N}{2}t}$.

Now let $T \rightarrow +\infty$ to obtain a suitably bounded
solution of (\ref{eq:3.16}). Again we must correct the boundary data,
so we define, for $j=1, \ldots ,N$, the solution $w_j$ by
\[
w_j(t) = {\bar{w}}(t) -\frac{2}{N-2} {\bar{w}}(\log R)(\Phi_{\e}^{j,+}
(\log R))^{-1}\Phi_{\e}^{j,+}(t).
\]
\begin{lemma}
Given $R>0$, there exists a constant $m_1 > 0$ such that for $\e \in 
(0,\e_0]$ and for any $j = 1, \ldots, N$, we have $\Phi_{\e}^{j,\pm}
(\log R)  \geq m_1  \e $. 
\label{le:3.6}
\end{lemma}
{\bf Proof :}
Recall that 
\[
\Phi_{\e}^{j,+} (t) = e^{-t}\left(\frac{N-2}{2} v_\e(t) - 
\frac{\del v_\e}{\del t}(t)\right).
\]
By Proposition~\ref{pr:2.3}, we get
\[
\Phi_{\e}^{j,+} (t) = \e \frac{N-2}{2} e^{-\frac{N}{2} t} + O(\e^{
\frac{N+2}{N-2}} e^{\frac{N}{2} t}).
\]
and the result follows.
\hfill $\Box$

Again we remark that it will be important in the nonlinear analysis that we 
have obtained a solution with a $\Phi_\e^{j,+}$, but no $\Phi_\e^{j,-}$,
component, for the latter would force the location of the singularities 
of the exact solution to be different from those of the approximate
solution. 

Now that we have obtained a solution of the original equation in
Proposition~\ref{pr:3.1} and a (weighted) $L^{\infty}$ bound for 
this solution, 
estimates for the full H\"older norm follow by standard scaling
arguments. This ends the proof of the Proposition. 
\hfill $\Box$

\section{A model for the linearization about $\buex$}

In the next two sections we finally study the linearization of ${\cal N}$
about the approximate solution $\buex$ associated to some fixed
singular set $\Sigma = \{x_1, \ldots , x_n\}$, and parameters $\bar{R} = 
(R_1,\ldots ,R_n)$ and admissible $\bar{\e} = (\e_1, \ldots , \e_n)$, where 
each $R_i > 1$. This last condition can always be assumed since we may 
dilate our problem by some factor to $\kappa >0$, this will change 
the set $\Sigma$ into $\{\kappa x_1, \ldots , \kappa x_n\}$ and will 
change the parameters $R_i$ into $\kappa^{-1}R_i$.
As usual, the main point is to analyze this operator 
uniformly as $\e$ tends to zero. We shall do this in two steps, first 
studying the somewhat simpler operator
\begin{equation}
\Le w = \Delta w + \frac{N(N+2)}{4} \sum_{i=1}^{n}(\tilde{\chi}(x-x_i)
u_{\e_i}^{\frac{4}{N-2}}((x-x_i)/R_i)) w
\label{eq:5.1}
\end{equation}
where the displacements $a_i$ have all been set to zero, and then, in the 
next section, treating the true linearization as a perturbation of $\Le$. 
Here $\tilde{\chi}$ is a smooth cutoff function equalling $1$ in $B(0,1)$,
vanishing outside $B(0,2)$, and taking values in $[0,1]$ in
the annulus $B(0,2)\setminus B(0,1)$. The point of doing this is
that $\Le$ is much simpler to study because the term of order zero
is radial in each $B(x_i,1)$.  

The main result is that on suitably weighted function spaces, $\Le$ 
is invertible, with inverse blowing up as $\e$ tends to zero, but in a manner
which we can control precisely. We construct the inverse for
$\Le$ by solving the equation 
\begin{equation}
      \Le w  =  f \quad \mbox{in}\quad {\R}^N \setminus \Sigma\\[2mm]
\label{eq:5.2}
\end{equation}
in three steps: first solve the homogeneous Dirichlet problem for 
this equation outside the union of the balls $B(x_i,1)$; next, solve 
the homogeneous Dirichlet problem for this equation in each of these 
balls; finally show that the sum of these solutions can be modified
to a true solution using the Dirichlet to Neumann maps on the boundaries 
of these balls.  

The main result of this section has a statement parallel to that 
of Proposition~\ref{pr:3.1}:

\begin{proposition}
Let $\be = \e \bq$ be an admissible set of Delaunay parameters. Suppose 
that $\bR = (R_1, \ldots R_n)$, $n \ge 3$, is a collection of numbers,
with each $R_i > 1$, and satisfying (\ref{eq:4.1}). Suppose also that
\begin{equation}
1 < \nu < 2.
\label{eq:5.3}
\end{equation}
Define the deficiency space 
\[
{\cal W}_0 \equiv \mbox{span\,}\{\tilde{\chi}\Psi_{\e_i}^{j,+}(
\frac{x-x_i}{R_i}), \ j = 0, \ldots, N\}.
\]
Then for $\e$ sufficiently small, the operator
\begin{equation}
\Le: {\cal C}^{2,\al}_{\nu,2}(\R^N \setminus \Sigma) \oplus
{\cal W}_0 \longrightarrow {\cal C}^{0,\al}_{\nu-2,-2}(\R^N \setminus
\Sigma)
\label{eq:mapmodel}
\end{equation}
is an isomorphism. In particular, for each $f\in {\cal C}^{0,\al}_{\nu-2,
-2}(\R^N\setminus \Sigma)$ there exists a unique solution of 
(\ref{eq:5.2}) which has a decomposition 
\[
w(x) = {\cal G}_{{\be}, {\bR}} (f)(x) +\sum_{i=1}^n \tilde{\chi}(x-x_i) 
\sum_{j=0}^N {\cal K}^j_{{\be}, {\bR},i} (f)\frac{1}{\e_i}\Psi^{j,+}_{\e_i}
((x-x_i)/R_i),
\]
The operator 
\[
{\cal G}_{\be,\bR}:{\cal C}^{0,\alpha}_{\nu-2,-2}(\R^N \setminus\Sigma)
\longrightarrow {\cal C}^{2,\alpha}_{\nu ,2} ({\Bbb R}^N\setminus\Sigma)\]
and functionals 
\[
{\cal K}^j_{\be,\bR,i}:{\cal C}^{0,\alpha}_{\nu-2,-2}(\R^N\setminus\Sigma)
\longrightarrow \R, \qquad j = 0, \ldots, N, \quad i = 1, \ldots, n 
\]
are bounded independently of $\e$.
% and $\bR$. 
%HERE ALSO I THINK THAT WE DONT NEED UNIFORMITY w.r.t. \bR
%IN ADDITION, IF WE DO NEED THIS UNIFORMITY THEN WE HAVE TO INCLUDE IT 
%IN THE PREVIOUS SECTION. NAMELY THAT ALL THE RESULTS ARE TRUE UNIFORMLY
%AS LONG AS R IS AWAY FROM 0, AWAY FROM 1 AND BOUNDED.
\label{pr:5.1}
\end{proposition}

We shall denote the right inverse of $\Le$ constructed here by $\Lei$. 
As noted earlier, the proof will be done in several steps.
\medskip

{\bf Step 1 : The exterior problem.} Let $\Omega = \R^N\setminus\cup_{i=1}^{n}
B(x_i, 1)$. 
\begin{proposition}
Let $f\in {\cal C}^{0,\alpha}_{\nu-2,-2 }({\R}^N\setminus\Sigma)$. Then 
for $\e_0 >0$ sufficiently small, there exists a unique solution 
$w\in {\cal C}^{2,\alpha}_{\nu,2}(\Omega)$ to 
\begin{equation}
\left\{\begin{array}{rll}
           \Le w & = & \bar{f} \quad \mbox{in}\quad \Omega \\[2mm] 
               w & = & 0 \quad \mbox{on}\quad \del \Omega,
         \end{array}
\right.
\label{eq:5.4}
\end{equation}
where $\bar{f}$ is the restriction of $f$ to $\Omega$. This solution
satisfies the estimate 
\[
||w||_{2,\alpha,\nu,2} \leq c \|\bar{f} \|_{0,\alpha,\nu-2,-2},
\]
with some constant $c>0$ independent of $\e_i \in (0, \e_0]$.
The norms here are taken in $\Omega$.
\label{pr:5.2}
\end{proposition}
{\bf Proof :} By Proposition~\ref{pr:2.3} the term of order zero in $\Le$ 
has compact support in $\Omega$, and is bounded by $C\e$. 
We can transform (\ref{eq:5.4}) to a problem of the
form $\Delta w' + \e V (x,{\be}) w' = \bar{f}'$ in some compact set
$\bar{\Omega}$ with smooth boundary by taking the Kelvin transform
at some point outside of $\Omega$. The function $V(x,{\be})$ is bounded 
in $\bar{\Omega}$ independently of ${\be}$. 
$\bar{f} \in {\cal C}^{0,\alpha}_{\nu-2,-2} (\Omega)$ implies that 
$\bar{f}' \in {\cal C}^{0,\alpha}(\bar{\Omega}\setminus\{0\})$. 
The existence and uniform estimate are now standard.
\hfill $\Box$
\medskip

{\bf Step 2 : The exterior Dirichlet to Neumann map}.
We now introduce the exterior Dirichlet to Neumann map. Let $\Psi =
\{\psi_1 (\theta) , \ldots , \psi_n (\theta) \} \in \oplus_{i=1}^n 
{\cal C}^{2, \alpha} (\del B(x_i, 1))$ be a set of boundary values.
It is standard, using Proposition~\ref{pr:5.2}, that there exists
a unique solution ${\cal C}^{2, \alpha}_{\nu,2}(\Omega)$ of the 
homogeneous problem  
\begin{equation}
\left\{ \begin{array}{rlll}
           \Le w & = & 0 & \quad \mbox{in}\quad \Omega \\[2mm] 
               w & = & \psi_i & \quad \mbox{on}\quad \del B(x_i, 1) , 
\quad \forall \ i=1, \ldots ,n,
         \end{array}
\right.
\label{eq:5.5}
\end{equation}
which we will denote by $w_\Psi$. The correspondence from $\Psi$ to
$w_\Psi$ is continuous for each $\be$ and $\bR$. Define ${\Se}$ by 
\[
{\Se} (\Psi) = \left(\left.\del_{r_1}w_{\Psi}\right|_{r_1 = 1},
\ldots , \left. \del_{r_n}w_{\Psi}\right|_{r_n = 1}\right) \in 
\oplus_{i=1}^n {\cal C}^{1,\alpha}(\del B(x_i, 1)),
\]
$r_i=|x-x_i|$. This is the Dirichlet to Neumann map for the operator $\Le$ 
on $\Omega$. 
\begin{lemma}
The norms of the mappings
\[
\oplus_{i=1}^n {\cal C}^{2,\alpha}(\del B(x_i, 1)) \ni \Psi 
\longrightarrow w_\Psi \in {\cal C}^{2, \alpha}_{\nu,2} (\Omega)
\]
and
\[
\Se: \oplus_{i=1}^n {\cal C}^{2,\alpha}(\del B(x_i, 1)) \longrightarrow 
\oplus_{i=1}^n {\cal C}^{1,\alpha}(\del B(x_i, 1)),
\]
are bounded independently of $\be$ provided each $\e_i < \e_0$ and
$\e_0$ is sufficiently small. Furthermore, if all $\e_i$ tend to $0$,
${\Se}$ converges to a limiting operator ${\SO}$ which is the Dirichlet 
to Neumann map for the Laplacian, $\Delta$, on $\Omega$.
\label{le:5.1}
\end{lemma} 
{\bf Proof :}
The boundedness follows from Proposition~\ref{pr:5.2}. Convergence 
of the ${\Se}$ is a consequence of the fact that $\Le$ tends to
zero uniformly with $\be$. 
\hfill $\Box$
\medskip

{\bf Step 3 : The interior Dirichlet to Neumann map}.
We have already proved the analogue of Proposition~\ref{pr:5.2}
in Proposition~\ref{pr:3.1}, so we may pass directly to the
corresponding interior Dirichlet to Neumann map. As above, 
the unique solution of the homogeneous problem
\begin{equation}
\left\{ \begin{array}{rlll}
           \Le w & = & 0  & \quad \mbox{in}\quad B(x_i, 1) \\[2mm]
               w & = & \psi_i  & \quad \mbox{on}\quad \del B(x_i, 1), 
         \end{array}
\right.
\label{eq:5.6}
\end{equation}
with
\[
w = w_{\psi_i}\in {\cal C}^{2, \alpha}_{\nu}(B(x_i, 1)\setminus \{x_i\})
 \oplus_{j=0}^N 
\mbox{Span}\{\frac{1}{\e_i}\Psi_{\e_i}^{j,+}((\cdot -x_i)/R_i)\}
\]
defines a continuous map
\[
\Psi \in \oplus_{i=1}^n {\cal C}^{2,\alpha}(\del B(x_i,1))\longrightarrow 
(w_{\psi_1}, \ldots , w_{\psi_n}) \in \oplus_{i=1}^{n} \big(
{\cal C}^{2, \alpha}_{\nu}(B(x_i, 1)\setminus\{x_i\})
\]
\[ 
\oplus_{j=0}^N 
\mbox{Span}\{\frac{1}{\e_i}\Psi_{\e_i}^{j,+}((\cdot -x_i)/R_i)\}\big).
\]
Now define the map $\Te$ by
\[
\Te (\Psi) = (\Te^1(\psi_1), \ldots , \Te^n(\psi_n)) = 
\left(\left. \del_{r_1}w_{\psi_1}\right|_{r_1 = 1},\ldots ,\left.\del_{r_n}
w_{\psi_n}\right|_{r_n = 1}\right) \in 
\oplus_{i=1}^n {\cal C}^{1,\alpha}(\del B(x_i, 1)).
\]
\begin{lemma}
The norms of the mappings
\[
{\cal C}^{2,\alpha}(\del B(x_i,1)) \ni \psi_i \longrightarrow 
w_{\psi_i} \in {\cal C}^{2, \alpha}_{\nu}(B(x_i, 1)\setminus \{x_i\})
\oplus_{j=0}^N \mbox{Span}\{\frac{1}{\e_i}\Psi_{\e_i}^{j,+}
((\cdot -x_i)/R_i)\}
\]
and 
\[
\Te:\oplus_{i=1}^n {\cal C}^{2,\alpha}(\del B(x_i, 1)) \longrightarrow 
\oplus_{i=1}^n {\cal C}^{1,\alpha}(\del B(x_i, 1))
\]
are bounded independently of $\e_i$ provided all $\e_i < \e_0$ and
$\e_0$ is sufficiently small. The mappings $\Te$ converge, as
$\be$ tends to $0$, to a limiting operator ${\cal T}_0$ which
acts diagonally on $n$-tuples. In terms of the eigenfunctions $\phi_j$ for 
the Laplacian on $\del B(x_i,1)$, the $i^{\mbox{th}}$ component 
$\TOi$ is determined by
\[
\TOi (\phi_j)= \left(\frac{2-N}{2} +\left(\frac{(N-2)^2}{4} +\lambda_j
\right)^{1/2}\right)\phi_j \quad \mbox{for all}\quad j \geq 1,
\]
\[
\TOi (\phi_0)=-\frac{N-2}{2}\left(1 + \frac{\cosh ((N-2)\log R_i/2)}{
\sinh ((N-2)\log R_i /2)}\right)\phi_0 = (2-N)\left(\frac{R_i^{N-2}}{
R_i^{N-2} - 1}\right) \phi_0.
\] 
\label{le:5.2}
\end{lemma}
{\bf Proof :}
We only need demonstrate the convergence of the $\Te$ and the specific
form of the limit ${\cal T}_0$. This will follow from the proof of 
Lemma~\ref{le:3.3}.  

By Lemma~\ref{le:3.3}, for each $j\geq N+1$, there exists a unique
solution $y_j(t)$ of the equation
\begin{equation}
\frac{d^2 y }{dt^2} - \left(\frac{(N-2)^2}{4} + \lambda_j - 
\frac{N(N+2)}{4}v_{\e_i, R_i}^{\frac{4}{N-2}}\right)y = 0.
\label{eq:5.7}
\end{equation}
which satisfies $y_j(0)=1$ and which decays exponentially as $t 
\rightarrow +\infty$. In this equation $v_{\e_i,R_i} = v_{\e_i}(t+\log R_i)$. 
Let $\gamma_j =(\frac{(N-2)^2}{4} +\lambda_j)^{1/2}$ and $\delta_j =
(\gamma_j^2 - \frac{N(N+2)}{4})^{1/2}$. As in the proof of that Lemma, 
we show that 
\[
|y_j(t)| \leq e^{-\delta_j t }  \quad \mbox{for all} \quad t\geq 0.
\]
Integrating the product of (\ref{eq:5.7}) with $e^{-\gamma_j t}$ from $0$ 
to $+\infty$ and integrating by parts, we find that
\begin{equation}
\frac{d y_j}{dt}(0) +\gamma_j = \int_0^{+\infty} \frac{N(N+2)}{4}\,
v_{\e_i,R_i}^{\frac{4}{N-2}}\,y_j(t)\, e^{-\gamma_j t}\,dt.
\label{eq:5.8}
\end{equation}
Since $|v_{\e_i,R_i}(t)|\leq 1$, we obtain
\[
|\frac{d y_j}{dt}(0) +\gamma_j | \leq \int_0^{+\infty}\frac{N(N+2)}{4}\, 
e^{-(\gamma_j +\delta_j) t}\,dt \leq \frac{N(N+2)}{4} 
(\gamma_j +\delta_j)^{-1}.
\]
Moreover, using (\ref{eq:5.8}) we see that, for fixed $j \geq N+1$,
\[
\lim_{\e_i \rightarrow 0} \frac{d y_j}{dt} (0) =- \gamma_j.
\]
Now, the solution of (\ref{eq:5.6}) with $\psi_i = \phi_j$ is given by 
$w(x) = |x-x_i|^{\frac{2-N}{2}}y_j(-\log |x-x_i|) \phi_j (\theta)$, and its Neumann data 
is $\big( \frac{2-N}{2} - \frac{d y_j}{dt}(0)\big)  \phi_j (\theta)$. The relevant solution
when $\e_i = 0$ is $|x-x_i|^{\frac{2-N}{2} +\gamma_j} \phi_j (\theta)$. 
Hence if $\psi(\theta) = \sum_{j=N+1}^{+\infty}a_j\phi_j(\theta)$, we see that
\[
||\Tei (\psi)-\TOi (\psi)||_{L^2}^2 =\sum_{j=N+1}^{+\infty}|a_j|^2 
\left(\left(\frac{2-N}{2}-\frac{d y_j}{dt}(0)\right)-\left(\frac{2-N}{2} +
\gamma_j\right)\right)^2
\]
and this expression tends to $0$ with $\e_i$.

When $j = 0, \ldots , N$, we can proceed by direct computation. Indeed, 
when $j = 1, \ldots, N$ and $\psi_i$ is the eigenfunction $\phi_j$, then 
\[
\left(\Psi_{\e_i}^{j,+}(\log R_i)\right)^{-1} \Psi_{\e_i}^{j,+}
(-\log|x-x_i| + \log R_i) \phi_j(\theta)
\]
is the solution of (\ref{eq:5.6}). An explicit computation, using
Proposition~\ref{pr:2.3}, shows that
\[
\lim_{\be \rightarrow 0} {\Tei}(\phi_j) = \phi_j,
\]
as expected. Similarly, when $j = 0$, the explicit solution of 
(\ref{eq:5.6}) is given by 
\[
\left(\Psi_{\e_i}^{0,+}(\log R_i)\right)^{-1}\Psi_{\e_i}^{0,+}(
-\log|x-x_i| + \log R_i).
\] 
Using Proposition~\ref{pr:2.3} again, we get
\[
\lim_{\be \rightarrow 0} {\Tei} (\phi_0) =-\frac{N-2}{2}\left(
1 + \frac{\cosh ((N-2)\log R_i /2)}{\sinh ((N-2)\log R_i /2)}\right).
\]
\hfill $\Box$
\medskip

{\bf Step 4: Invertibility of the difference of Dirichlet to Neumann maps}.
In order to glue the interior and exterior solutions together, we must
add correction terms to these solutions. These correction terms are
solutions of the homogeneous problem, and are chosen so that the Cauchy
data from the inside and outside match up. To find these, we must show
that the difference $\Se - \Te$ is invertible when all $\e_i$
are  sufficiently small.
\begin{proposition}
There exists $\e_0 >0$ such that for $\e <\e_0$, 
\[
\Se- \Te:\oplus_{i=1}^{n}{\cal C}^{2,\alpha}(\del B(x_i, 1))
\longrightarrow \oplus_{i=1}^{n}{\cal C}^{1,\alpha}(\del B(x_i,1))
\]
is invertible. The norm of its inverse is bounded by some constant $C >0$ 
independent of $\e$.
\label{pr:5.3}
\end{proposition}
{\bf Proof~:} By the $L^2$ operator norm convergence of $\Se - \Te$ to
$\SO - \TO$, it suffices to show that $\SO - \TO$ is invertible. 
Now $\SO-\TO$ is a self-adjoint first order pseudodifferential
operator. Both $\SO$ and $\TO$ are elliptic, with principal
symbols $|\xi|$ and $-|\xi|$, respectively, hence the 
difference is also elliptic and semibounded. This means that 
$\SO - \TO$ has discrete spectrum, and thus we need only prove that it 
is injective. The invertibility in H\"older spaces then follows by 
standard regularity theory. 

We argue by contradiction. Assume that $\SO-\TO$ is not injective. Then 
there exists some $\Psi_{0}\in\oplus_{i=1}^{n}{\cal C}^{\infty}(\del B(x_i,1))$ 
for which $(\SO-\TO)(\Psi_0)= 0$. We may extend the Dirichlet data
$\Psi_0$ to a harmonic function on the exterior region, decaying
at infinity, and also one 
on the interior region. $\SO$ is the Dirichlet to Neumann map for
the ordinary Laplacian on the exterior region, while $\TO$ is the
Dirichlet to Neumann map for the Laplacian on the union of balls
only for the spherical eigencomponents with index $j \ge 1$. 
Thus $\Psi_0$ extends to harmonic functions $w'$ and $w''$ on the
interior and exterior regions, respectively. These functions have
the same Dirichlet data, namely $\Psi_0$, and by assumption on
this Dirichlet data, all eigencomponents with $j\ge 1$ of their 
Neumann data agree as well. To find a harmonic function with the 
$j = 0$ eigencomponent of the Neumann data also matching at these spheres, 
it suffices to add on appropriate
multiples of the functions $|x - x_i|^{2-N} - 1$ on each
$B(x_i,1)$. We obtain in this way a harmonic function $w$ on
$\R^N \setminus \Sigma$ which decays at infinity, and with at
worst a radial singular term at each $x_i$. The only possibility
is that this function has the form 
\[
w(x) = \sum_{i=1}^n p_i |x-x_i|^{2-N}.
\]

We now find the radial part, $w^i_0(x)$, of the harmonic function $w(x)$ 
in each ball $B(x_i,1)$, and use this to compute $T_0^i \Psi_0$ directly.
Comparing this with the previous answer will lead to a contradiction.
In $B(x_{i_0},1)$, this radial part is a sum of the singular part
and a constant:
\[
w^{i_0}_0(x) = p_{i_0}|x - x_{i_0}|^{2-N} + \sum_{i \neq i_0}
\frac{1}{\omega_{N-1}} \int_{S^{N-1}}p_i |\theta + x_{i_0} - x_i|^{2-N}\,
d\theta.
\]
Here $\omega_{N-1}$ is the volume of $S^{N-1}$. The terms in the sum 
may be computed by noticing that for any fixed vector $v$, the function
\[
\frac{1}{\omega_{N-1}}\int_{S^{N-1}} |x + v|^{2-N}\, d\theta
\]
is harmonic, regular and radial in $B(0,|v|)$, hence constant. Thus,
provided $|v| > 1$, its 
value when $|x| = 1$ is the same as its value at $x = 0$, and is just
$|v|^{2-N}$. We conclude that 
\[
w^{i_0}_0(x) = p_{i_0}|x - x_{i_0}|^{2-N} + \sum_{i \neq i_0}p_i |x_i -
x_{i_0}|^{2-N}.
\]
The normal derivative of this function at $B(x_{i_0},1)$ is just 
$(2-N)p_{i_0}$, but on the other hand, it must agree
with $T_0^{i_0}(w^{i_0}_0)$. This leads to the equality
\[
(2-N)p_{i_0} = (2-N)\frac{R_{i_0}^{N-2}}{R_{i_0}^{N-2} - 1}
\left( p_{i_0} + \sum_{i \neq i_0} p_i |x_i - x_{i_0}|^{2-N} \right),
\]
or simply
\[
p_{i_0} + R_{i_0}^{N-2} \sum_{i \neq i_0} p_i |x_i - x_{i_0}|^{2-N} = 0.
\]

Now recall that the parameters $R_i$ satisfy the relation (\ref{eq:4.1}), 
\[
\sum_{i\neq i_0}(R_i)^{\frac{N-2}{2}}(R_{i_0})^{\frac{N-2}{2}}q_i
|x_{i_0}-x_i|^{2-N} =q_{i_0}.
\]
A small calculation shows that the collection of
numbers $\tilde{p}_i \equiv \displaystyle{\frac{p_i}{q_i}}
(R_i)^{\frac{2-N}{2}}$ give a solution of the system 
\begin{equation}
\sum_{i\neq i_0}(R_i)^{\frac{N-2}{2}}(R_{i_0})^{\frac{N-2}{2}}q_i 
q_{i_0}|x_{i_0}-x_i|^{2-N}(\tilde{p}_i+\tilde{p}_{i_0})=0, \qquad 
i_0 =1, \ldots ,n.
\label{eq:5.9}
\end{equation}
Let $s_{i,j}=(R_i)^{\frac{N-2}{2}}(R_{j})^{\frac{N-2}{2}}q_i 
q_j|x_j-x_i|^{2-N}$. Multiply (\ref{eq:5.9}) by $\tilde{p}_{i_0}$ and 
sum over $i_0$ to obtain
\[
\sum_{i_0=1}^n(\sum_{i\neq i_0}s_{i,i_0}(\tilde{p}_{i}+\tilde{p}_{i_0}))
\tilde{p}_{i_0}  = 
\sum_{i<i_0} s_{i,i_0}(\tilde{p}_i +\tilde{p}_{i_0})^2=0.
\]
Since all the $s_{i,i_0}$ are positive, we get that $\tilde{p}_i=
-\tilde{p}_{i_0}$. But since $n \geq 3$, this forces the vanishing
of all $\tilde{p}_i$. This is a contradiction, and the proof is finished.
\hfill $\Box$.
\medskip

{\bf Step 5: The correction term}. 
Let $f \in{{\cal C}}^{0,\alpha}_{\nu-2,-2}(\R^N\setminus\Sigma)$. By
Proposition~\ref{pr:3.1} and Proposition~\ref{pr:5.2}, there exist functions 
$w_{ext}$ and $w_{in,i}$, $i=1, \ldots , n$, such that 
\[
\left\{ \begin{array}{rllll} 
              \Le w_{ext}& = & f&  \quad \mbox{in}\quad \Omega\\[2mm] 
                  w_{ext}& = & 0 & \quad \mbox{on} \quad \del \Omega ,
        \end{array}
\right. 
\]
and, for $i=1, \ldots , n$,
\[
\left\{ \begin{array}{rllll}
                 \Le w_{int,i}& =& f \quad \mbox{in}\quad B(x_i, 1) \\[2mm] 
                     w_{int,i}& =& 0 \quad \mbox{on} \quad \del B(x_i, 1) .
        \end{array}
\right. 
\]
From Proposition~\ref{pr:3.1}, there is a decomposition inside each 
$B(x_i,1)$ 
\[
w_{int,i}(x)=G_{\e_i,R_i}(f_i)(x)+\sum_{j=0}^n K_{\e_i,R_i}^{j}(f_i)(x)
\frac{1}{\e_i}\Psi_{\e_i}^{j,+}(\frac{x-x_i}{R_i}), 
\] 
where $f_i = f\upharpoonleft _{B(x_i,1)}$. In addition, for some $c >0$
\[
||G_{\e_i,R_i}(f_i)||_{2,\alpha,\nu}+\sum_{j=0}^N |K_{\e_i,R_i}^{j}(f_i)|
\leq  c ||f||_{0,\alpha,\nu-2,-2},  
\]
and
\[
||w_{ext}||_{2,\alpha ,\nu, 2} \leq c ||f||_{0,\alpha,\nu-2, -2}.
\]

We now seek a function $w_{ker}$ in $\R^N\setminus\{\cup_{i=1}^{n}\del 
B(x_i,1)\cup\Sigma\}$ so that 
\[
w(x) \equiv\left\{ \begin{array}{rll}    
                 w_{ext}(x) +w_{ker}(x) &\quad\mbox{in}\quad\Omega \\[2mm] 
                 w_{int,i}(x) +w_{ker}(x) & \quad\mbox{in}\quad B(x_i,1),
\quad i=1, \ldots , n,
                   \end{array}
           \right.
\]
will satisfy $\Le w = f$ on $\R^N\setminus\Sigma$. 

For this to be true we need to choose $w_{ker}$ to be continuous
across each $\del B(x_i,1)$, to satisfy $\Le w_{ker}= 0$ away from
these boundaries (and from $\Sigma$), and such that the jump
of $\del_{r_i}w_{ker}$ across $\del B(x_i,1)$ equals $\del_{r_i}w_{int,i}-
\del_{r_i} w_{ext}$. This is equivalent to finding a solution $\Psi$ of 
\[
(\Se -\Te)(\Psi)= \{(\del_{r_1}w_{int,1}-\del_{r_1}w_{ext}),\ldots,
(\del_{r_k}w_{int,k}-\del_{r_k}w_{ext})\}\in\oplus_{i=1}^n{\cal C}^{1,
\alpha}(\del B(x_i, 1)).
\]
This is possible by Proposition~\ref{pr:5.3}; this same result gives the 
uniform bound on the norm of $\Psi$. All estimates for the norm of $w$ 
follow from Lemma~\ref{le:5.1} and Lemma~\ref{le:5.2}.
\hfill $\Box$

\section{ The true linearization }
We now undertake the main task of analyzing the linearization about $\buex$. 
Recall that our ultimate goal is to find a solution of the nonlinear equation
\begin{equation}
{\cal N}(w)\equiv\Delta(\bue +w)+\frac{N(N-2)}{4}(\bue +
w)^{\frac{N+2}{N-2}}=0,
\label{eq:8.1}
\end{equation}
for some choice of parameters $\ba, \bR$ (recall that $\be$ is regarded 
as fixed), and for some $w$ which is dominated by $\buex$ near $\Sigma$ . Let 
\begin{equation} 
\Laa = \Delta + \frac{N(N+2)}{4}(\bue)^{\frac{4}{N-2}}
\label{eq:8.2}
\end{equation}
denote the linearization of this equation around $w=0$. As is usual
for a Jacobi operator, one can obtain solutions of $\Laa w = 0$ as
derivatives of families of solutions of (\ref{eq:8.1}). Of course, we do not
know any families of solutions yet, but we {\it can} obtain approximate 
solutions of $\Laa w = 0$ from families of approximate solutions of 
(\ref{eq:8.1}), and we know many explicit families of approximate solutions. 
In particular, the functions
\[
\Psi_{\be,\bR,\ba, i}^{0,+}(x)=\del_{R_i} \buex , \quad \mbox{and}
\quad \Psi_{\be,\bR,\ba,i}^{ j,+}(x)=\del_{a^j_i} \buex, \quad j = 1, 
\ldots, N,\quad i = 1, \ldots, n, 
\]
defined by differentiating $\bar{u}_{\be}(\bR,\ba,x)$ with respect to
$R_i$ and $a_i^j$, respectively, have the property that 
\[
\Laa \Psi_{\be,\bR,\ba, i}^{j,+}(x)
\]
vanishes near each $x_i$. This may be seen by differentiating (\ref{eq:8.1}) 
with respect to either $R_i$ or $a_i^j$.

In this section we shall consider the problem $\Laa w = f$, for
$f \in {\cal C}^{0,\al}_{\nu-2,-2}(\RNS)$. The main linear result
in this paper is that we can find a solution with optimal bounds
in the space
\begin{equation}
{\cal M}(\ba) = {\cal C}^{2,\al}_{\nu,2}(\RNS) \oplus \mbox{Span\,} \{\Psi_{\be,\bR,\ba, i}^{j,+} \},
\label{eq:8.3}
\end{equation}
where the norm of an element $w(x) = v(x) + \sum_i S_i \Psi_{\be,\bR,\ba, i}^{0,+} + \sum_{i}\sum_{j=1}^{N} \al_i^j \Psi_{\be,\bR,\ba, i}^{j,+} $ is given by 
\[
||w||_{{\cal M}(\ba)}=\rho^{\nu}||v||_{2,\al,\nu,2}+\e\sum_{i}|S_i|
+\e\rho\sum_{i,j}|\al^j_i|.
\]

As before, the deficiency space spanned by the $\Psi_{\be,\bR,\ba,i}^{j,+}$ 
is necessary to obtain surjectivity of the linearized operator. The
rather byzantine choice of norms on ${\cal M}(\ba)$ is necessitated
by the rather intricate estimates below, and the need for uniformity
as $\e$ tends to zero. The following proposition states the expected
bijectivity of $\Laa$, but it also gives a rather sharper estimate
on the solution in terms of the linear functionals ${\cal K}_{\be,\bR,i}^{j}$
introduced in Proposition~\ref{pr:5.1}. The values ${\cal K}_{\be,\bR,i}^j(f)$
are closely related to the coefficients $S_i$ and $\al_i^j$ for
the solution $w \in {\cal M}(\ba)$ of $\Laa w = f$, but they correspond
to coefficients of the simpler, radial solutions $\Psi_{\e_i}^{j,+}$.
It may seem rather unnatural to continue to use these radial solutions
here, but once again, these sharper estimates seem necessary later.

For notational convenience, we shall often denote these approximate solutions by 
\begin{equation}
\Psi_{\be,\bR,\ba,i}^{0,+}(x) = \mu_i (\be, \bR, \ba,x) , \qquad \Psi_{\be,\bR,\ba,i}^{j,+}(x)
= \gamma_i^j(\be, \bR, \ba,x),
\label{eq:convlabel}
\end{equation}
and we shall identify  $w(x) = v(x) + \sum S_i \mu_i (\be,\bR,\ba ,x) + 
\sum \al_i^j \gamma_i^j (\be,\bR,\ba, x)$ with $(v,\bar{S}, \bar\alpha)$.

\begin{proposition} 
Let $\be$ be an admissible set of Delaunay parameters, and suppose
that $R_i >1$ and $a_i \in\R^N$ are a collection of numbers and
vectors satisfying (\ref{eq:4.1}) and (\ref{eq:4.3}). Then 
\[
\Laa: {\cal M}(\ba) \longrightarrow {\cal C}^{0,\al}_{\nu-2,-2}(\RNS)
\]
is surjective, provided all $\e_i$ are less than some sufficiently
small number $\e_0$. Furthermore, for $f\in {\cal C}^{0,\al}_{\nu-2,-2}
(\RNS)$, there exists $w = (v,\bS,\bal)$, corresponding to
\[
w(x)=v(x) + \sum_{i = 1}^N S_i \mu_i(\be, \bR, \ba, x) + \sum_{i=1}^n\sum_{j=1}^N \al_i^j \gamma_i^j(\be, \bR, \ba, x),
\]
solution of the equation $\Laa w=f$, which satisfies 
\[
\rho^{\nu}||v||_{2,\al,\nu,2}\leq c (\rho^{\nu}||f||_{0,\al,\nu-2,-2}+ 
\sup_{i}|{\cal K}_{\be, \bR,i}^{0}(f)| +\rho \sup_{i,j}
|{\cal K}_{\be, \bR,i}^{j}(f)|).
\]
\[
\e |\bS|\leq c(\rho^{\nu+1}||f||_{0,\al,\nu-2,-2}+
\sup_{i}|{\cal K}_{\be,\bR,i}^{0}(f)| +\rho^2 \sup_{i,j}
|{\cal K}_{\be, \bR,i}^{j}(f)|),
\]
and
\[
\e\rho|\bal|\leq c(\rho^{\nu +1}||f||_{0,\al,\nu-2,-2}+
\rho \sup_{i} |{\cal K}_{\be, \bR,i}^{0}(f)| + \rho 
\sup_{i,j}|{\cal K}_{\be,\bR,i}^{j}(f)|). 
\]
\label{pr:7.1}
\end{proposition}
{\bf Proof :} 
Define the space
\[
{\cal M} = {\cal C}^{2,\al}_{\nu,2}(\RNS) \oplus \mbox{Span\,}
\{\tilde{\chi}(x-x_i)\tilde{\mu}_i, 
\tilde{\chi}(x-x_i) \tilde{\gamma}_i^j\},
\]
where now $\tilde{\chi}$ is the cutoff equalling one in each $B(x_i,1)$
and vanishing outside the union of the $B(x_i,2)$ and
\[
\tilde{\mu}_i (\e_i, R_i, x)= \del_{R_i}u_{\e_i}(R_i,0,x-x_i), \qquad
\tilde{\gamma}_i^j (\e_i, R_i,x) = \del_{a_i^j}u_{\e_i}(R_i,0,x-x_i).
\]
An element $w = v + \sum S_i \tilde{\chi}(x-x_i)\tilde{\mu}_i + 
\sum \al_i^j\tilde{\chi}(x-x_i)\tilde{\gamma}_i^j$ has the norm
\[
||w||_{{\cal M}} \equiv ||(v,\bS,\bal)||_{{\cal M}} = 
||v||_{2, \alpha, \nu, 2} + \e \sum_i |S_i| +\e\sum_{i,j} |\al_i^j|.
\]
Proposition~\ref{pr:5.1} may then be rephrased as the assertion that
\begin{equation}
\Le: {\cal M} \longrightarrow {\cal C}^{0,\al}_{\nu-2,-2}(\RNS)
\label{eq:10.1}
\end{equation}
is an isomorphism, and that the inverse is bounded independently of $\be$.

We shall prove this proposition by reducing it to Proposition~\ref{pr:5.1}
by a sequence of perturbations.

\bigskip
{\bf Step 1}. We first claim that 
\begin{equation}
\LaO \equiv \Delta +\frac{N(N+2)}{4}(\bar{u}_{\be}(\bR,0,x))^{\frac{4}{N-2}}: 
{\cal M} \longrightarrow {\cal C}^{0,\al}_{\nu-2,-2}(\RNS)
\label{eq:10.2}
\end{equation}
is also an isomorphism, with inverse bounded independently of $\be$. 
To see this, we estimate the norm of the difference $A = \LaO - \Le$.
First, the difference
\[
|(\bar{u}_{\be}(\bR,0,x))^{\frac{4}{N-2}}- \sum_i\tilde{\chi}(x-x_i)
(u_{\e_i}(R_i,0,x-x_i))^{\frac{4}{N-2}}|
\]
vanishes in each $|x-x_i|\le \rho_i$ and is estimated by $c\e^{\frac{4}{N-2}}
\mbox{dist}(x,\Sigma)^{-4}$ when $\mbox{dist}(x,\Sigma)\geq 1$. We get the
estimate 
\[
\e^{\frac{4}{N-2}}(|x-x_i|^{N-6}\rho + |x-x_i|^{N-5})
\]
for this difference in the annulus $\rho_i \le |x-x_i| \le 1$ 
from Corollary~\ref{cor:useful}.

Finally, from Proposition~\ref{pr:2.3} we derive the bounds 
\begin{equation}
|\tilde{\mu}_i(\e_i, R_i,x)| \leq C\e|x-x_i|^{2-N} \quad \mbox{and}\quad
|\tilde{\gamma}_i^j(\e_i, R_i,x) |\leq C \e |x-x_i|,
\label{eq:10.25}
\end{equation}
which hold in particular for all $x$ satisfying $\mbox{dist}(x,x_i)\in
[\rho_i, 2]$. Thus, for $w \in {\cal M}$, we get
\[
||(\LaO-\Le)w||_{0,\al,\nu-2,-2} = ||(\bar{u}_{\be}(\bR,0,x)^{
\frac{4}{N-2}}-\sum_{i=1}^n \tilde{\chi}(x-x_i)
u_{\e_i}(R_i,0,x-x_i)^{\frac{4}{N-2}})w||_{0,\al,\nu-2,-2} 
\]
\begin{equation}
\leq c \rho^{N+2}( ||v||_{2,\al,\nu,2}+\e\rho^{-\nu-1}|\bS|+
\e(1+\rho^{N-2-\nu})|\bal|).
\label{eq:7.11}
\end{equation}
By (\ref{eq:defrho}), we can estimate the right side by an
arbitrarily small multiple $\eta$ of $||w||_{{\cal M}}$, provided $\e$
is sufficiently small, i.e. $||Aw||_{0,\al,\nu-2,-2} \le \eta 
||w||_{{\cal M}}$. Letting $\Le^{-1}$ denote the inverse of the
map (\ref{eq:10.1}), then
\begin{equation}
\LaO^{-1} \equiv (\Le + A)^{-1} = \Le^{-1}( I + A\Le^{-1})^{-1}
\label{eq:10.3}
\end{equation}
is a right inverse of the map (\ref{eq:10.2}). Since $\Le$ is
invertible, and since $A$ is small in operator norm,
$\LaO$ is also invertible, hence (\ref{eq:10.3}) is the
unique inverse. Since $A$ and $\Le^{-1}$ are bounded independently of 
$\be$, so is $\LaO^{-1}$. 

We now obtain sharper estimates for the solution of $\LaO w = f$. 
\begin{lemma} Suppose that $f$ is any function in ${\cal C}^{0,\al}_{
\nu-2,-2}$, and let $w = \LaO^{-1}f = (v,\bS,\bal)$ and $w_0 = 
\Le^{-1}f = (v_0, \bS_0, \bal_0)$. Then
\[
||v||_{2,\al,\nu,2} \le C \big(||v_0||_{2,\al,\nu,2} + \rho^{N+1-\nu}\e
|\bS_0| + \rho^{N+2}(1 + \rho^{N-2-\nu})\e |\bal_0|\big),
\]
\[
\e |\bS| \le C \big( \rho^{N+2}||v_0||_{2,\al,\nu,2} + \e |\bS_0|
+ \rho^{N+2}(1 + \rho^{N-2-\nu})\e |\bal_0|\big),
\]
and
\[
\e |\bal| \le C \big( \rho^{N+2}||v_0||_{2,\al,\nu,2} +
\rho^{N+1-\nu} \e |\bS_0| + \e |\bal_0| \big).
\]
\label{le:diffLaL}
\end{lemma}
\begin{remark} We point out explicitly again that the coefficients 
$\bS_0$ and $\bal_0$ here are the values of the functionals 
${\cal K}^j_{\be,\bR,i}(f)$ introduced in Proposition~\ref{pr:5.1},
corresponding to the radial model operator $\Le$. Specifically, 
\[
\e_i (S_0)_i ={\cal K}_{\be,\bR,i}^{0}(f), \quad \e_i (\al_{0})_i^j =
{\cal K}_{\be,\bR,i}^{j}(f).
\]
Although it might seem more natural to view the coefficients $\bS$ and 
$\bal$ coming from the solution of $\LaO w = f$  (and later, $\Laa w = f$) 
in terms of them, it turns out to be simpler in the long run to 
view these as primary. 
\end{remark}
{\bf Proof:} Rewrite (\ref{eq:10.3}) as 
\[
\LaO^{-1} = \Le^{-1}(I + A\Le^{-1})^{-1}(I + A \Le^{-1} - A \Le^{-1}) = 
\Le^{-1}(I - (I + A\Le^{-1})^{-1}A\Le^{-1}).
\]
Applying this to $f$ yields 
\[
w = w_0 - \Le^{-1}(I + A\Le^{-1})^{-1} A w_0.
\]
The estimates of the lemma follow by considering the different components 
of the two sides of this equation, and using the boundedness of 
$\Le^{-1}(I + A\Le^{-1})^{-1}$, and (\ref{eq:7.11}). \hfill $\Box$

Using these estimates, and the fact that 
$||w_0||_{{\cal M}} \le C ||f||_{0,\al,\nu-2,-2}$, we conclude
\[
||v||_{2,\al,\nu,2}\leq c(||f||_{0,\al,\nu-2,-2}+\rho^{N+1-\nu}\sup_{i}|
{\cal K}_{\be,\bR,i}^{0}(f)| + \rho^{N+2}(1+\rho^{N-2-\nu})\sup_{i} 
\sup_{j = 1, \ldots, N}|{\cal K}_{\be, \bR,i}^{j}(f)|),
\]
\[
\e |\bS| \leq c(\rho^{N+2}||f||_{0,\al,\nu-2,-2}
+\sup_{i}|{\cal K}_{\be,\bR,i}^{0}(f)| + \rho^{N+2}(1+\rho^{N-2-\nu})
\sup_{i}\sup_{j=1,\ldots ,N} |{\cal K}_{\be, \bR,i}^{j}(f)|)  
\]
and
\[
\e |\bal|\leq c(\rho^{N+2}||f||_{0,\alpha, \nu-2,-2} 
+\rho^{N+1-\nu} \sup_{i} |{\cal K}_{\be, \bR,i}^{0}(f)|+\sup_{i}\sup_{j=1, 
\ldots ,N}|{\cal K}_{\be, \bR,i}^{j}(f)|).
\] 

\medskip

{\bf Step 2}. We now perturb further. We claim that 
$\LaO$ is invertible from ${\cal M}(0)$ to ${\cal C}^{0,\al}_{\nu-2,-2}$,
with inverse uniformly bounded in $\be$. In fact, because of the
uniform invertibility of (\ref{eq:10.2}), it suffices to show that
the inclusion
\[
{\cal M} \longrightarrow {\cal M}(0)
\]
is uniformly bounded. Equivalently, we must show that
\[
||w||_{{\cal M}(0)} \le C ||w||_{{\cal M}}
\]
for some $C$ independent of $\e$. 
To calculate the norm of $w$ in ${\cal M}_0$, we must write this function
in the form $w = \bar{v} + \sum S_i \mu_i (\be, \bR, 0, \cdot) + \sum \al_i^j \gamma_i^j (\be, \bR, 0, \cdot)$.
Since $w = v + \tilde{\chi}(\sum S_i \tilde{\mu}_i (\e_i, R_i,\cdot ) + \sum \al_i^j \tilde{\gamma}_i^j (\e_i, R_i, \cdot))$ we get 
\[
\bar{v} = v + \sum S_i ( \tilde{\chi}\tilde{\mu}_i(\e_i, R_i, \cdot) - \mu_i (\be, \bR,0, \cdot) )
+ \sum \al_i^j (\tilde{\chi}\tilde{\gamma}_i^j (\e_i, R_i, \cdot) - \gamma_i^j  (\be, \bR, \ba, \cdot)).
\]
Now observe that 
\[
\tilde{\mu}_i - \mu_i = (1 - \chi_i)\del_{R_i}(u_{\e_i}(R_i,0,x-x_i)- 
\bar{w}_{\be}),\qquad \tilde{\gamma}_i^j - \gamma_i^j = 
(1-\chi_i)\del_{a_i^j}(u_{\e_i}(R_i,0,x-x_i)- \bar{w}_{\be})
\]
in each $B(x_i,1)$. From (\ref{eq:4.diffuw})
\begin{equation}
||\tilde{\chi}\tilde{\mu_i} - \mu_i||_{2,\al,\nu,2} \le C\rho^{-\nu}\e,
\qquad \mbox{and} \quad ||\tilde{\chi}\tilde{\gamma}_i^j - 
\gamma_i^j||_{2,\al,\nu,2} \le C \rho^{1-\nu}\e.
\label{eq:12.80}
\end{equation}
Hence
\[
||\sum S_i (\tilde{\chi}\tilde{\mu}_i - \mu_i)||_{2,\al,\nu,2} +
||\sum \al_i^j (\tilde{\chi}\tilde{\gamma}_i^j - \gamma_i^j)||_{2,\al,\nu,2}
\le C \rho^{-\nu} (\e |\bS| + \e \rho |\bal|),
\]
and thus
\[
||w||_{{\cal M}(0)} = \rho^{\nu}||\bar{v}||_{2,\al,\nu,2} 
+ \e |\bS| + \e \rho |\bal| \le C ( ||v||_{2,\al,\nu,2} + 
\e |\bS| + \e \rho |\al|) = C ||w||_{{\cal M}},
\]
as desired. 

We also derive the precise estimates for the inverse of 
\begin{equation}
\LaO: {\cal M}(0) \longrightarrow {\cal C}^{0,\al}_{\nu-2,-2}(\RNS)
\label{eq:10.4}
\end{equation}
from the ones stated at the end of Step 1. In fact, replace $v$ by 
$\bar{v}$ there and use (\ref{eq:12.80}) and the sharp estimates for 
$\e|\bS|$ and $\e\rho |\bal|$ to get 
\[
\rho^{\nu}||\bar{v}||_{2,\al,\nu,2} \leq c(\rho^{\nu}||f||_{0,\al,\nu-2,-2}+ 
\sup_{i}|{\cal K}_{\be,\bR,i}^{0}(f)| +\rho\sup_{i}\sup_{j=1,\ldots ,N}
|{\cal K}_{\be, \bR,i}^{j}(f)|),
\]
\[
\e|\bS| \leq  c(\rho^{N+2} ||f||_{0,\al,\nu-2,-2}+ \sup_{i}|{\cal K}_{\be,
\bR,i}^{0}(f)| +\rho^{N+2}(1+ \rho^{N-2-\nu})\sup_{i}\sup_{j=1, 
\ldots ,N} |{\cal K}_{\be, \bR,i}^{j}(f)|),
\]
and
\[
\e \rho |\bal| \leq  c(\rho^{N+3}||f||_{0,\al,\nu-2,-2}+\rho^{N+2-\nu}
\sup_{i} |{\cal K}_{\be, \bR,i}^{0}(f)|) + \rho\sup_{i}
\sup_{j=1, \ldots ,N} |{\cal K}_{\be,\bR,i}^{j}(f)|.
\]
\medskip 

{\bf Step 3.} We perform a final perturbation of (\ref{eq:10.4})
to show that 
\begin{equation}
\Laa: {\cal M}(\ba) \longrightarrow {\cal C}^{0,\al}_{\nu-2,-2}(\RNS)
\label{eq:10.5}
\end{equation}
is an isomorphism, with inverse uniformly bounded as $\be$ tends to zero.
In order to apply the same sort of perturbation argument, we would
like $\Laa$ and $\LaO$ to be acting on the same space of functions.
To this end, define
\[
\tcM = \{(v,\bS,\ba) \in {\cal C}^{2,\al}_{\nu,2}(\RNS) \oplus
(\R^+)^n \oplus (\R^N)^n \},
\]
with norm given by
\[
||(v,\bS,\bal)||_{\tcM} = \rho^{\nu}||v||_{2,\al,\nu,2}+ \e |\bS| +\rho
\e |\bal|.
\]
Then $\tcM$ is equivalent to any of the spaces ${\cal M}(\ba)$ via
the map 
\[
\iota_{\ba}: (v,\bS,\bal) \longrightarrow v(x) + \sum S_i \mu_i(\be, \bR, \ba,x) +
\sum \al_i^j \gamma_i^j(\be, \bR, \ba,x).
\]
We denote
\[
\tLaa = \Laa \circ \iota_{\ba}: \tcM \longrightarrow {\cal C}^{0,\al}_{
\nu-2,-2}(\RNS).
\]

The uniformly bounded inverse $\tLaO^{-1}$ is a first approximation
to the inverse of $\tLaa$. In fact,
\[
\tLaO^{-1}\tLaa = I + \tLaO^{-1}(\tLaa - \tLaO).
\]
We shall show that the final term $\tLaO^{-1}(\tLaa - \tLaO)$
has small norm as a map on $\tcM$. Because of the uniform boundedness
of $\tLaO^{-1}$, it should suffice to show that the difference
\[
\tLaa - \tLaO: \tcM \longrightarrow {\cal C}^{0,\al}_{\nu-2,-2}(\RNS)
\]
has small norm. This is almost true, although at one stage of the 
argument we require slightly more information about the inverse. 
However, let us concentrate on showing that most terms in this
difference are small.

Let $w = (v,\bS,\ba) \in \tcM$, and for convenience write 
\[
U_{\ba}(x) = \sum S_i \mu_i(\be, \bR, \ba,x) + \sum \al_i^j \gamma_i^j(\be, \bR, \ba,x).
\]
Then 
\begin{equation}
(\tLaa - \tLaO)(w) = (\buex^{\frac{4}{N-2}} - 
\bar{u}_{\be}(\bR,0,x)^{\frac{4}{N-2}}) v
+ (\Laa U_{\ba}(x) - \LaO U_{0}(x)).
\label{eq:10.6}
\end{equation}
We estimate the two terms on the right in turn.

The first is simpler. Going back to the definition of $\buex$, we 
expand in $B(x_i, \rho_i)$ to get 
\[
(\buex)^{\frac{4}{N-2}}- (\bar{u}_{\be}(\bR,0,x))^{\frac{4}{N-2}} 
=
\]
\[
2 a_i \cdot (x-x_i)|x-x_i|^{-2} v_{\e_i}^{\frac{6-N}{N-2}}
(-\log |x-x_i| +\log R_i)( v_{\e_i} (-\log |x-x_i| +\log R_i)
\]
\begin{equation}
-\frac{2}{N-2}\dot{v}_{\e_i}(-\log |x-x_i| +\log R_i))
+O(1)v_{\e_i}^{\frac{4}{N-2}} (-\log |x-x_i| +\log R_i).
\label{eq:7.2}
\end{equation}
Here we have used the fact that for all $T>0$, there exists some
$c_T >0$, independent of $\e \le \e_0$, such that for all $t_0 \in [-T,T]$
\[
v_{\e}(t+t_0) +|\dot{v}_{\e}(t+t_0)|+|\ddot{v}_{\e}(t+t_0)|\leq c_T v_{\e}(t).
\]

We consider this difference in three separate regions about each $x_i$.
In the innermost one, where $\mbox{dist}(x,x_i)\leq \rho^{\frac{N+2}{2}}
= \e^{\frac{2}{N-2}}$, since $v_\e$ and $\dot{v}_\e$ are bounded we get
\[
|(\buex)^{\frac{4}{N-2}}-(\bar{u}_{\be}(\bR,0,x))^{\frac{4}{N-2}}
|\leq C|x-x_i|^{-1}.
\]
For the middle region, where $\mbox{dist} (x,x_i)\in [\e^{\frac{2}{N-2}},
\rho_i]$, we further expand (\ref{eq:7.2}), using Proposition~\ref{pr:2.3},
to get
\begin{equation}
|(\buex)^{\frac{4}{N-2}} - (\bar{u}_{\be}(\bR,0,x))^{\frac{4}{N-2}}
|\leq C(\e^{\frac{4}{N-2}}|x-x_i|^{-2}+\e^{\frac{8}{N-2}}|x-x_i|^{-5}).
\label{eq:7.21}
\end{equation}
A similar estimate also holds if $\mbox{dist}(x,x_i)\in [\rho_i, 2\rho_i]$.
Finally, when $\mbox{dist}(x,x_i)\geq 2\rho_i$ for all $i$, then
\[
(\buex)^{\frac{4}{N-2}} -(\bar{u}_{\be}(\bR,0,x))^{\frac{4}{N-2}}=0,
\]
since $\bu$ does not depend on $\ba$ if $\mbox{dist}(x,x_i) \geq 2\rho_i$.
The estimates for the full H\"older norms are similar.

Putting these estimates together, we get
\begin{equation}
||(\buex ^{\frac{4}{N-2}} - \bar{u}_{\be}(\bR,0,x)^{\frac{4}{N-2}})
v||_{0,\al,\nu-2,-2}\leq c\e^{\frac{2}{N-2}}||v||_{2,\al,\nu,2}.
\label{eq:10.10}
\end{equation}

The second term on the right in (\ref{eq:10.6}) is more complicated
to estimate. First write
\[
\Laa(U_{\ba}) -\LaO(U_0) = \Delta (U_{\ba} - U_0) 
\]
\begin{equation}
+ {\frac{N(N+2)}{4}}\left((\buex)^{\frac{4}{N-2}}U_{\ba} - 
(\bar{u}_{\be}(\bR,0,x))^{\frac{4}{N-2}}U_0\right).
\label{eq:10.11}
\end{equation} 
Recall that $U_{\ba}$ satisfies
\[
\Delta U_{\ba} +\frac{N(N+2)}{4}(\buex)^{\frac{4}{N-2}}U_{\ba} =0
\]
when $\mbox{dist}(x,x_i) \leq \rho_i$. Hence in this ball, (\ref{eq:10.11}) 
vanishes identically. It also vanishes when $\mbox{dist}(x,x_i)\geq 2\rho_i$,
because $\buex$ does not depend on $\ba$ and so the difference cancels.

We come at last to the transition annulus, where $\mbox{dist}(x,x_i)\in 
[\rho_i,2\rho_i]$. Recalling that $U_{\ba}$ is the differential of
$\buex$ and that $\buex = \chi_i (x-x_i) \uix + (1-\chi_i )\bar{w}_{\be}(x)$,
we write
\[
U_{\ba} = \chi_i U_{a_i}^i + (1-\chi_i) W_{\be}
\]
where $U_{a_i}^i$ and $W_{\be}$ are the differentials of 
$u_{\e_i}$ and $\bar{w}_{\be}$, respectively. In this annulus we
can estimate
\begin{equation}
|U_{a_i}^i| \le C (\rho^{2-N} \e |S_i| + \e \rho |\al|), \qquad
\mbox{and} \quad |W_{\be}| \le C \e \rho^{2-N} |S_i|.
\label{eq:10.30}
\end{equation}
Now write (\ref{eq:10.11}) as a sum of three terms, $I + II + III$, where
\[
I = (\Delta \chi_i) (U_{a_i}^i - U_0^i) + 2 \nabla \chi_i \cdot \nabla (
U_{a_i}^i - U_0^i),
\]
\[
II = \chi_i \Delta (U_{a_i}^i - U_0^i) +  \chi_i \frac{N(N+2)}{4}\big( \buex^{\frac{4}{N-2}}U_{a_i}^i - 
\bar{u}_{\be}(\bR,0,x)^{\frac{4}{N-2}} U_0^i \big),
\]
and
\[
III = (1 - \chi_i) \frac{N(N+2)}{4} W_{\be} \big( \buex^{\frac{4}{N-2}}
- \bar{u}_{\be}(\bR,0,x)^{\frac{4}{N-2}} \big).
\]
In the last of these, we have used that $W_{\be}$ does not depend on $\ba$.

This last term is the simplest to estimate. In fact, by 
Corollary~\ref{cor:useful}, 
\[
|\bar{u}_{\be}(\bR,\ba,x)^{\frac{4}{N-2}}-\bar{u}_{\be}(R,0,x)^{\frac{4}{
N-2}}| \le \e^{\frac{4}{N-2}}\rho^{N-5}.
\]
Hence 
\begin{equation}
||III||_{0,\al,\nu-2,-2} \le C \rho^{2-\nu}\e 
\rho^{2-N} \e^{\frac{4}{N-2}} \rho^{N-5} |\bS|
\le C (\rho^{N + 1-\nu}) \e |\bS|,
\label{eq:10.31}
\end{equation}
using (\ref{eq:defrho}). 

To handle the middle term, recall that
\[
\Delta U_{a_i}^i + \frac{N(N+2)}{4}\uix^{\frac{4}{N-2}} U_{a_i}^i = 0.
\]
Using this, we rewrite $II$ as
\[
\chi \frac{N(N+2)}{4} \big( ( \buex^{\frac{4}{N-2}} - \bar{u}_{\e_i}
(R_i,a_i,x)^{\frac{4}{N-2}})U_{a_i}^i - (\bar{u}_{\be}(\bR,0,x)^{
\frac{4}{N-2}} - u_{\e_i}(R_i,0,x-x_i)^{\frac{4}{N-2}}) U_0^i \big).
\]
Since
\[
|\buex^{\frac{4}{N-2}} - \uix^{\frac{4}{N-2}}|
\le C \e^{\frac{4}{N-2}}\rho^{N-5},
\]
we get
\begin{equation}
||II||_{0,\al,\nu-2,-2} \le C \rho^{2-\nu}( \rho^{2-N} \e |\bS|
+ \e \rho |\bal|) \e^{\frac{4}{N-2}} \rho^{N-5} \le \eta (\e |\bS|
+ \e \rho |\bal|),
\label{eq:10.32}
\end{equation}
where $\eta$ tends to zero with $\e$. 

The first term, $I$, is the most difficult to estimate suitably. 
First observe that 
\[
U_{a_i}^i - U_0^i = S_i \del_{R_i}(\uix - u_{\e_i}(R_i,0,x-x_i))
\]
\[
+ \sum_j \al_i^j \del_{a_i^j}(\uix - u_{\e_i}(R_i,0,x-x_i)).
\]
From the proof of  Proposition~\ref{pr:4.1} we see that
\begin{equation}
\uix - u_{\e_i}(R_i,0,x-x_i) = \e_i \frac{N-2}{2}
R_i^{\frac{2-N}{2}}a_i \cdot (x - x_i) + O(\e|x - x_i|^2).
\label{eq:10.33}
\end{equation}
Noting that the derivative with respect to $a$ of the leading 
term of this expression is independent of $a$, hence cancels in the
difference, we get
\[
|U_{a_i}^i - U_0^i| \le C \rho(\e|\bS|+\e\rho |\bal|), \qquad \mbox{and}
\quad |\nabla(W_{\ba} - W_0)| \le C(\e|\bS| + \e \rho |\bal|).
\] 
Thus we see that both terms in $I$ are of size $ \rho^{-1} (\e |\bS|
+ \e \rho |\bal|)$, and in particular, they are not small relative
to the norm of $\tcM$. Fortunately, the term we actually want to 
estimate, $\LaO^{-1}(I)$, is still small because $I$ has support
in the union of annuli $B(x_{i}, 2\rho_{i})\setminus B(x_{i}, \rho_{i})$,
which has small volume. The following result will be proved in the 
next section.

\begin{lemma}
Suppose $\nu \in (1,2)$. Let $h$ be the term $I$. Expand it in terms of the 
eigenfunctions on $S^{N-1}$: 
\[
h(x_i +r \theta )=\sum_{j=0}^{\infty}h_j(r) \phi_j(\theta),
\]
and let 
\[
h' = \sum_{j=1}^N h_j(r) \phi_j(\theta) \qquad \mbox{and} \qquad
h'' =\sum_{j=N+1}^{\infty}h_j(r) \phi_j(\theta).
\]
Then the solution $\hat{w} = \Le^{-1}h = (\hat{v},\hat{S},\hat{\al})
\equiv (\hat{v}, {\cal K}^0(h), {\cal K}'(h))$ satisfies
\[
||v_0||_{2, \alpha, \nu, 2} \leq C \big( \rho^{2-\nu}\sup |h_0| + 
\rho^{2-\nu} \sup |h'| + ||h''||_{0,\al,\nu-2}\big).
\] 
In addition, for any $\mu\in (-N, 2)$, we have 
\[
\e |\hat{S}_i| \equiv |{\cal K}^{0}_{\be,\bR, i} (h)| \leq C(\rho^2
\sup |h_0(r)| + \rho^{N+1} \sup |h'(r)| + ||h''||_{0,\al,\mu-2}),
\]
and, for $j=1,\ldots, N$, 
\[
\e |\hat{\al}_i^j| \equiv |{\cal K}^{j}_{\be,\bR,i}(h)| 
\leq C(\rho^2 \sup |h_0(r)|+\rho \sup |h'(r)|+||h''||_{0,\al,\mu-2,-2}).
\]
If $i' \neq i$, then for $j = 1, \ldots, N$, 
\[
|\hat{\al}_{i'}^j| \equiv |{\cal K}^j_{\be, \bR,i'}(h)|\leq C(\rho^2 
\sup |h_0(r)| + \rho^{N+1} \sup |h'(r)| + ||h''||_{0,\al,\mu-2,-2}).
\]
\label{le:6.2}
\end{lemma}

The point of this result, of course, is that because $h$ is supported
in a region with small volume, one can obtain better estimates for 
$\Le^{-1}h$ (and then also for $\LaO^{-1}h$) than the obvious ones. 
We apply this as follows. Using (\ref{eq:10.33}), we see
that the terms in its eigenfunction expansion satisfy
\[
|h_0| \le \e (|\bS| + |\bal|), \qquad |h'| \le C\rho^{-1}( \e|\bS|
+ \e \rho |\bal|),
\]
and
\[
||h''||_{0,\al,\mu-2,-2} \le C \e \rho^{2-\mu} (|\bS| + |\bal|),
\]
for any $\mu \in (-N,2)$. We also see that
\[
\rho^\nu||h||_{0,\al,\nu-2,-2} \le C (\rho^{2} \sup |h_0|
+ \rho^2 \sup |h'| + \rho^\nu ||h''||_{0,\al,\nu-2,-2}) \le
C\rho ( \e |\bS| + \e \rho |\bal|).
\]

Now apply the Lemma to find a solution $\Le^{-1}h = \hat{w} = (\hat{v},
\hat{S},\hat{\al})$. It satisfies 
\[
\rho^{\nu}||\hat{v}||_{2,\al,\nu,2} \le C( \rho^{2}(\e(|\bS| + |\bal|)
+ \rho^{2}(\rho^{-1}(\e|\bS| + \e \rho |\bal|) + \rho^\nu( \e 
\rho^{2-\nu} (|\bS| + |\bal|)), 
\]
\[
\e |\hat{S}| \le C \big(\rho^2\e(|\bS| + |\bal|) + \rho^{N}(\e|\bS|+
\e\rho|\bal|) + \rho^{2 - \mu} \e (|\bS| + |\bal|)\big),
\]
and
\[
\e \rho |\hat{\al}|\le C\big(\rho^3 \e(|\bS|+|\bal|)+\rho^2 (\rho^{-1}(\e|\bS|
+\e \rho|\bal|)) + \e \rho^{3 -\mu}(|\bS| + |\bal|)\big).
\]
In the second and third expressions, we have estimated $h''$ in the 
${\cal C}^{0,\al}_{\mu-2,-2}$ norm, for any $\mu \in (-N,2)$.
By taking $\mu = 0$, for example, we see that we can bound all three of these
expressions by
\[
C \rho (\e |\bS| + \e \rho |\bal|).
\]

The solution $\tilde{w} = \LaO^{-1}h$ that we are actually interested in
can be estimated in terms of these quantities, using Lemma~\ref{le:diffLaL}.
Combined with the preceding, we see that
\[
||\tilde{w}||_{\tcM} \le \rho (\e |\bS| + \e \rho |\bal|).
\]

We combine this estimate with the estimates (\ref{eq:10.10}), (\ref{eq:10.31})
and (\ref{eq:10.32}), and use the uniform boundedness of $\LaO^{-1}$
for these terms to finally conclude that the operator norm of 
$\LaO^{-1}(\tLaa - \tLaO)$ on $\tcM$ can 
be made as small as desired by choosing $\e$ small. This will show that 
\[
\tLaO^{-1}\tLaa = I + \tLaO^{-1}(\tLaa - \tLaO)
\]
is invertible, hence
\[
G = (I + \tLaO^{-1}(\tLaa - \tLaO) )^{-1} \tLaO^{-1}
\]
is a left inverse for $\tLaa$. Clearly it is uniformly bounded in $\e$.
Direct computation shows that
\[
\tLaa G=(\tLaO +(\tLaa -\tLaO))G =\tLaO (I+\tLaO^{-1}(\tLaa -\tLaO))G = I
\]
as well. Thus $\tLaa$ has a two-sided inverse, which is therefore unique,
and finally, $\Laa$ does as well.

\medskip
{\bf Step 4.} We conclude the proof of the theorem by showing that
the sharper estimates for the solution of $\Laa w = f$ are valid.
We proceed as in Lemma~\ref{le:diffLaL}. Thus, let $w_0 = 
(v_0, \bS_0, \bal_0) =\LaO^{-1}(f)$ and $w = (v,\bS,\bal) = 
\Laa^{-1}f$. Then writing $B = \Laa - \LaO$, we derive
\[
\Laa^{-1} = (I + \LaO^{-1}B)^{-1}\LaO^{-1} = (I + \LaO^{-1}B - 
\LaO^{-1}B) (I + \LaO^{-1}B)^{-1}\LaO^{-1} = 
\]
\[
\LaO^{-1} -
\LaO^{-1}B(I + \LaO^{-1}B)^{-1}\LaO^{-1}.
\]
Apply this to $f$ to get
\begin{equation}
w = w_0 - \LaO^{-1}B w_1,
\label{eq:12.1}
\end{equation}
where we have set $w_1 = (I + \LaO^{-1}B)^{-1}w_0$. Clearly we have
\begin{equation}
||w_1||_{\tcM} \le C ||w_0||_{\tcM}.
\label{eq:12.3}
\end{equation}
We divide $w_1$ into a sum $v_1 + \tilde{w}_1$, where we have just
lumped the second and third terms in the usual decomposition into
$\tilde{w}_1$. Thus if we let 
\[
-\LaO^{-1}Bv_1 = (\hat{v}, \hat{S}, \hat{\al}), \qquad \mbox{and}\quad
- \LaO^{-1}B\tilde{w}_1 = (\tilde{v},\tilde{S}, \tilde{\al}),
\]
then we have
\[
(v,\bS,\bal) = (v_0,\bS_0,\bal_0) + (\hat{v},\hat{S},\hat{\al})
+ (\tilde{v},\tilde{S},\tilde{\al}).
\]

We estimate these new terms in turn. First, from (\ref{eq:10.10}),
\[
||\LaO^{-1}Bv_1||_{\tcM} = 
\rho^{\nu}||\hat{v}||_{2,\al,\nu,2} + \e |\hat{S}| + \e\rho|\hat{\al}| 
\le C \rho^{\frac{N+2}{2}}||v_1||_{2,\al,\nu,2} \le C \rho^{\frac{N+2}{2}
-\nu} ||w_1||_{\tcM}. 
\]
To estimate this, use (\ref{eq:12.3}) and the sharp estimates for
$w_0$ from the end of Step 2 to conclude that this term is bounded by
\[
\rho^{\frac{N+2}{2} - \nu}( \rho^{\nu}||f||_{0,\al,\nu-2,-2} +
\sup_i |{\cal K}^0_{\be,\bR,i}(f)| +
\rho \sup_i \sup_j |{\cal K}^j_{\be,\bR,i}(f)| ).
\]
For the other term, use Lemma~\ref{le:6.2} to estimate 
$||\LaO^{-1}B\tilde{w}_1||_{\tcM}$ first in terms of $w_1$ and then
$w_0$. Finally, substitute the sharp estimates from the
end of Step 2 for $w_0$ in terms of $f$. We omit the details,
which are straightforward. In the end, we get altogether
\[
\rho^{\nu}||v||_{2,\al,\nu,2}\leq c(\rho^{\nu}||f||_{0,\al,\nu-2,-2}+
\sup_{i}|{\cal K}_{\be, \bR,i}^{0}(f)| +\rho\sup_{i}\sup_{j=1,\ldots ,N}
|{\cal K}_{\be, \bR,i}^{j}(f)|),
\]
\[
\e|\bS| \leq  c( \rho^{\nu+1}||f||_{0,\al,\nu-2,-2}+\sup_{i}
|{\cal K}_{\be,\bR,i}^{0}(f)|) +\rho^2 \sup_{i}\sup_{j=1, \ldots ,N}
|{\cal K}_{\be, \bR,i}^{j}(f)|
\]
and 
\[
\e\rho |\bal|\leq  c(\rho^{\nu+1}||f||_{0,\al,\nu-2,-2}+
\rho \sup_{i} |{\cal K}_{\be,\bR,i}^{0}(f)|) +\rho\sup_{i}\sup_{j=1, 
\ldots ,N}|{\cal K}_{\be,\bR,i}^{j}(f)|).
\]
The proof is complete.
\hfill $\Box$

\section{Estimates on the error term}

In this section we shall derive estimates on the size of the
error term in the approximate solution. 
\begin{lemma}
Suppose that the parameters $R_i >1$ and $a_i\in {\R}^N$ satisfy 
the relations (\ref{eq:4.1}) and (\ref{eq:4.3}) of Proposition~\ref{pr:4.1}.
Let $\bu = \bu (\bR, \ba, \cdot)$, and define $\zeta \equiv \Delta 
\bu +\frac{N(N-2)}{4} \bu^\p$; then set $w = \Le^{-1}\zeta$. Then, 
for $\nu \in (1,2)$, and for some constant $C$ independent of $\e$, 
\[
||w||_{\tcM} \le C \e \rho^2.
\]
\label{le:6.1}
\end{lemma}
{\bf Proof :} 
Write $\Le^{-1}\zeta = w = (v,\bS,\bal)$ as usual. The estimate
for $v$ is rather straightforward; however, the estimates for
$\bS$ and $\bal$ are more delicate, since they rely on the particular
structure of $\zeta$, and so we need to follow the construction of 
$\Le^{-1}$ in Proposition~\ref{pr:5.1} step by step. 

{\bf Step 1.} From the definition of $\bu$, we see that in 
${\R}^N\setminus \cup_{i=1}^n B(x_i, 2\rho_i)$
\[
\zeta (x)=\frac{N(N-2)}{4} \bu^\p(x),
\]
and so
\[
|\zeta (x)|\leq C \e^{\frac{N+2}{N-2}} \mbox{dist}(x,\Sigma)^{-N-2}. 
\]
In particular, in $\Omega = \R^N \setminus \cup_{i = 1}^n B(x_i,1)$, 
$|\zeta| \le C \e^{\frac{N+2}{N-2}}$. From Proposition~\ref{pr:5.2}, 
there exists a solution of
\begin{equation}
           \left\{ \begin{array}{rllll} 
               \Le w_{ext}&=&\zeta & \quad \mbox{in}\quad \Omega\\
                   w_{ext}&=&0  &\quad \mbox{on} \quad \del \Omega .
                     \end{array}
\right.
\label{eq:6.1}
\end{equation}
which satisfies
\[
||w_{ext}||_{2, \alpha , \mu, 2} \leq C\e^{\frac{N+2}{N-2}}.
\]

{\bf Step 2.} Our task now is to estimate $\zeta$ in each $B(x_i, 1)$. 
In this step, we consider its eigenfunction expansion,
\[
\zeta (x) =\sum_{j=0}^{+\infty}\zeta_j^i (r)\phi_j(\theta ),
\]
where $r = |x - x_i|$, and estimate the terms in three separate regions. 

From the definition of $\bu$, we see that inside each $B(x_i, \rho_i)$, 
$\zeta= \Delta\bu +\frac{N(N-2)}{4} \bu^\p=0$, and so all
$\zeta_j^i = 0$. 

Next, $\zeta= \frac{N(N-2)}{4} \bu^\p$ in $B(x_i, 1)\setminus 
B(x_i, 2\rho_i)$. Thus 
\[
|\zeta_0^i(r) |\leq c \e^{\frac{N+2}{N-2}} r^{-(N+2)}, \quad 
\sum_{j=1}^N |\zeta_j^i(r)|\leq c \e^{\frac{N+2}{N-2}} r^{-3},
\]
and
\begin{equation}
| \sum_{j=N+1}^{+\infty}\zeta_j^i(r)\phi_j(\theta )|
\leq C\e^{\frac{N+2}{N-2}} r^{-2}. 
\label{eq:6.2}
\end{equation}

Finally, in $B(x_i, 2\rho_i)\setminus B(x_i, \rho_i)$, we have
\[
\zeta (x)= 2 (\del_{r_i}\chi_i)(\del_{r_i}\uix -\del_{r_i}\bwx) + 
(\Delta \chi_i) (\uix - \bwx )
\]
\[
-\frac{N(N-2)}{4}\big(\chi_i(\uix)^{\p} - (\chi_i\uix + (1-\chi_i )\bwx )^\p
\big).
\]
From this we get the estimates
\begin{equation}
\sum_{j=0}^N |\zeta_j^i(r)| + |\sum_{j\geq N+1} \zeta_{j}^i(r)\phi_j(
\theta )| \leq C\e. 
\label{eq:6.3}
\end{equation}

{\bf Step 3.} By Proposition~\ref{pr:3.1}, we get, for $i=1, \ldots , n$,
solutions of the problem
\[
              \left\{ \begin{array}{rllll}
                   \Le w_{int,i}&=&\zeta & \quad \mbox{in}\quad B(x_{i}, 1) \\
                       w_{int,i}&=&0 &\quad \mbox{on} \quad \del B(x_{i},1).
                      \end{array}
              \right.
\]
These satisfy 
\[
w_{int, i} = G_{\e_i, R_i}(\zeta_{|B(x_i,1)}) +
\sum_{j=0}^N {K}_{\e_i,R_i}^{j}(\zeta_{|B(x_i,1)})\frac{1}{\e_i}
\Psi_{\e_i}^{j,+}((x-x_i)/R_i).
\]
We now estimate the different terms in this formula.

First consider the solution $w_0^i$ in the $j=0$ eigencomponent. 
We can represent $w_0^i$ explicitly in terms of $\zeta_0^i$ using
the variation of constants formula. Write 
\[
\psi_{\e_i, R_i}^{0,+}(r)=r^{\frac{2-N}{2}}\frac{\del v_{\e_i}}{\del t}
(-\log (r/R_i))=r^{\frac{2-N}{2}}\Phi_{\e_i}^{0,+}(-\log (r/R_i)).
\]
Using Proposition~\ref{pr:2.3}, and the fact that $R_i \ge 1 + c_0$,
we see that, up to a constant, $\psi_{\e_i,R_i}^{0,+} \sim \e r^{2-N}$
for $\rho \le r \le 1$. In particular, this function never vanishes.
This means that we can write 
\[
w_0^i (x)= {K}^{0}_{\e_i,R_i}(\zeta_{|B(x_i,1)}(x))
\frac{1}{\e_i}\Psi_{\e_i}^{0,+}((x-x_i)/R_i)+ \tilde{w}_0^i(x-x_i).
\]
where
\[
\tilde{w}_0^i(x) =\psi_{\e_i, R_i}^{0,+}(r)\int_{r}^{\rho_i}
(\psi_{\e_i, R_i}^{0,+}(s))^{-2} s^{1-N} \int_s^{\rho_i}
\psi_{\e_i, R_i}^{0,+}(t) t^{N-1} \zeta_0^i (t) \, dt\,ds, 
\]
and the first term in the decomposition of $w_0^i$ is intended
to correct the boundary value. From the bounds above on $\zeta_0^i$, and 
using the approximation above for $\psi$, we find that
\begin{equation}
|\tilde{w}_0^i|\leq C\e \rho^2 \quad \mbox{for all}\quad 
\rho_i \leq |x-x_i| \leq 1.  
\label{eq:6.4}
\end{equation}
Therefore, since $\e_i^{-1}\Psi_{\e_i}^{0,+}$ is bounded above and
below on $\del B(x_i,1)$, we also get
\[
|{K}^{0}_{\e_i, R_i}(\zeta_{|B(x_i,1)})| \leq c \e \rho^2.
\]

We also require an estimate for the normal derivative of $w_0^i$
on $\del B(x_i,1)$. To get this, we use the slightly different
representation
\[
w_0^i(r) = \psi_{\e_i,R_i}^{0,+}(r) \int_r^1 (\psi_{\e_i,R_i}^{0,+}(s))^{-2}
s^{1-N} \int_{s}^{\rho_i} \psi_{\e_i,R_i}^{0,+}(t)t^{N-1}\zeta_0^i(t)\,dt\,ds.
\]
Differentiating this gives
\[
|\del_{r}w_0^i(x)| \leq  (\psi_{\e_i, R_i}^{0,+}(1))^{-1} 
\int^1_{\rho_i} \psi_{\e_i, R_i}^{0,+}(t) t^{N-1} |\zeta_0^i (t) | dt\quad
\mbox{ for all }\,x \in \del B(x_i,1),
\]
and so
\begin{equation}
\sup_{\del B(x_i,1)} |\del_{r_i}w_0^i(x)|\leq c \e \rho^{2}.
\label{eq:6.5}
\end{equation}

For $j=1, \ldots , N$, we may perform a similar analysis. Letting 
$w_j^i$ denote the component of the solution in the $j^{\mbox{th}}$
eigencomponent, then as before
\[
w_j^i (r)= {K}^{j}_{\e_i,R_i}(\zeta_{|B(x_i,1)})
\frac{1}{\e_i}\Psi_{\e_i}^{j,+} ((x-x_i)/R_i) + \tilde{w}_j^i(r)
\]
where, denoting by
\[
\psi^{1,+}_{\e_i,R_i} (r) = r^{\frac{4-N}{2}} (\frac{N-2}{2} 
v_{\e_i}(-\log (r/R_i)) - \displaystyle{\frac{\del v_{\e_i}}{\del t}}(
-\log (r/R_i)))= r^{\frac{2-N}{2}} \Phi_{\e_i}^{1,+} (-\log (r/R_i)),
\]
we represent
\[
\tilde{w}^i_j(r) =\psi^{1,+}_{\e_i,R_i}(r)\int_{r}^{\rho_i} 
(\psi^{1,+}_{\e_i,R_i}(s))^{-2} s^{1-N} \int_s^{\rho_i} \psi^{1,+}_{\e_i,
R_i}(t)  \zeta_j^i (t) t^{N-1}\,dt\,ds.
\]
The estimates 
\begin{equation}
|\tilde{w}_j^i|\leq C\e \rho |x-x_i| \qquad \mbox{and}\quad 
|{K}^{j}_{\e_i, R_i}(\zeta_{|B(x_i,1)})| \leq C \e \rho
\label{eq:6.6}
\end{equation}
follow as above. 

The estimate for the normal derivative follows from the 
alternate representation
\[
w_j^i(x) =\psi_{\e_i, R_i}^{1,+}(r)\int_{r}^1 (\psi_{\e_i,R_i}^{1,+}(s))^{-2}
s^{1-N} \int_s^{\rho_i} \psi_{\e_i, R_i}^{1,+}(t) t^{N-1} \zeta_j^i (t)\,
dt\,ds,
\]
which yields 
\[
|\del_{r_i}w_j^i(x)| \leq (\psi_{\e_i, R_i}^{1,+}(1))^{-1}  \int^1_{\rho_i} \psi_{\e_i, R_i}^{1,+}(t) t^{N-1} |\zeta_j^i (t)| dt .
\]
The estimate 
\begin{equation}
|\del_{r_i}w^i_j(x)|\leq c \e \rho^{N+1}, \quad \forall x \in \del B(x_i,1),
\label{eq:6.7}
\end{equation}
follows immediately. 

Let $\tilde{\zeta}^{i}$ denote the sum over all higher eigencomponents
of $\zeta$. From (\ref{eq:6.2}) and (\ref{eq:6.3}), and from
Proposition~\ref{pr:3.1}, for each $\mu \in (-N,2)$, 
\begin{equation}
||G_{\e_i,R_i}(\tilde{\zeta}^i)||_{2,\alpha,\mu} 
\leq C \e \rho^{2-\mu},
\label{eq:6.8}
\end{equation}
where $C$ is independent of $\e$.  In particular, using this with
$\mu = \nu$, and also using (\ref{eq:6.4}) and (\ref{eq:6.6}), we get
\[
||G_{\e_i, R_i}(\zeta_{|B(x_i,1)})||_{2,\alpha, \nu}  \leq C \e 
\rho^{2-\nu}. 
\]

Finally, take $\mu =0$ in (\ref{eq:6.8}) to get the estimate 
\begin{equation}
||\del_{r_j} G_{\e_i,R_i}(\tilde{\zeta^i})||_{{\cal C}^{1, \alpha}(\del 
B(x_i, 1))}  \leq c \e \rho^2
\label{eq:6.9}
\end{equation}
on the boundary of $B(x_i, 1)$.
\medskip

{\bf Step 4}. The final step is to consider the size of the
correction term $w_{ker}$, which is chosen to make the 
Dirichlet and Neumann data match. Its size is regulated by
the size of the jump of the Neumann data of the interior
and exterior solutions. Specifically, its Dirichlet data
$\Psi = (\psi_1, \ldots, \psi_n) \in \oplus_{i=1}^n {\cal C}^{1,\al}(\del 
B(x_i, 1))$ is chosen to satisfy 
\[
({\Se} - {\Te})(\Psi)= \{\left.(\del_{r_1} w_{int,1}-\del_{r_1} w_{ext})
\right|_{r_1 = 1}, \ldots , \left.(\del_{r_n} w_{int,n}-\del_{r_n} w_{ext})
\right|_{r_n = 1} \}.
\]
Thus, we can estimate the size of $\Psi$ from (\ref{eq:6.5}), 
(\ref{eq:6.7}) and (\ref{eq:6.9}), and the fact that the norm of the right 
side is bounded by $C \e \rho^2$. Hence $||w_{ker}||_{2,\al,\nu,2}
\le C\e\rho^2$ as well. The estimate $||w||_{\tcM} \le C\e\rho^2$
then follows easily. \hfill $\Box$ 
\medskip

{\bf Proof of Lemma \ref{le:6.2}:} The details of this proof are very
similar to those of the proof above. In fact, proceeding through
the steps of the previous proof, we find first that $w_{ext} = 0$.
The solutions $w_{int,i}$ may be found by the same procedure.
The integral representations of $w_j^i$ show that for $j = 0$ and
$j = 1, \ldots, N$, the ${\cal C}^{2,\al}_{\nu,2}$ norms of these terms 
may be bounded by $\rho^{2-\nu} \sup |h_0|$ and $\rho^{2-\nu}|\sup h'|$, 
respectively. The remaining details are left to the reader. \hfill $\Box$
 
\section{The nonlinear fixed point argument}

We are now in a position to find a solution of the nonlinear problem.
As usual, we fix the singular set $\Sigma=\{ x_1, \ldots , x_n\}$ and 
positive parameters $q_1, \ldots , q_n >0$ associated to a set of
admissible Delaunay parameters $\be$. We also fix an approximate
solution $\buex$ associated to the data $\bR$ and $\ba$ 
satisfying (\ref{eq:4.1}) and (\ref{eq:4.3}). The perturbation
of $\buex$ to an exact solution will involve not only a decaying term, but
also a slight adjustment of the parameters $\bR$ and $\ba$. Thus,
recalling that an element $w$ in $\tilde{\cal M}$ has components 
$(v, \bS, \bal )$, and changing our previous notation slightly, 
we wish to solve the equation ${\cal N}(w)=0$ where by definition
\begin{equation}
 {\cal N}(w)\equiv\Delta(\bar{u}_{\be}(\bR +\bS , \ba +\bal, \cdot) +v)
+\frac{N(N-2)}{4}(\bar{u}_{\be}(\bR +\bS , \ba +\bal, \cdot) 
+v)^{\frac{N+2}{N-2}}.
\label{eq:9.1}
\end{equation}

Recalling the linearization $\tLaa$ of ${\cal N}$ at $w = 0$, we can 
rewrite (\ref{eq:9.1}), using a Taylor expansion, as
\[
{\cal N}(w) ={\cal N}(0) +\tLaa w 
+ \int_0^1 (D{\cal N}_{|(tv,\bS,\bal)}-D{\cal N}_{|(0,\bS,\bal)} ) 
(v,0,0)\, dt 
\]
\[
+ (D{\cal N}_{|(0,\bS,\bal)}-D{\cal N}_{|(0,0,0)} ) (v,0,0)  
+ \int_0^1 (D{\cal N}_{|(0, t\bS, t\bal)} -D{\cal N}_{|(0,0,0)}) 
(0,\bS, \bal)\,dt .
\]
We have used a somewhat more elaborate expression than usual for the 
remainder term; this is necessary in the estimates below.
Therefore, the equation ${\cal N}(w)=0$ can be written as
\[
w =-\tLaa^{-1} {\cal N}(0) 
-\tLaa^{-1} \left( \int_0^1 (D{\cal N}_{|(tv,\bS,\bal)}-D{\cal N}_{
|(0,\bS,\bal)})(v,0,0)\,dt \right.
\]
\[
+ (D{\cal N}_{|(0,\bS,\bal)}-D{\cal N}_{|(0,0,0)} ) (v,0,0)   
\left. + \int_0^1 (D{\cal N}_{|(0, t\bS, t\bal)} -D{\cal N}_{|(0,0,0)})
(0,\bS, \bal)\, dt \right).
\]

We define the mapping ${\cal P}: \tilde{\cal M} \longrightarrow \tilde{\cal M}$
by 
\[
{\cal P}(w) = -\tLaa^{-1} {\cal N}(0) -\tLaa^{-1} \left( \int_0^1 
(D{\cal N}_{|(tv,\bS,\bal)}-D{\cal N}_{|(0,\bS,\bal)} ) (v,0,0)\,dt \right.
\]
\[
+ (D{\cal N}_{|(0,\bS,\bal)}-D{\cal N}_{|(0,0,0)} ) (v,0,0) 
\left. + \int_0^1 (D{\cal N}_{|(0, t\bS, t\bal)} -D{\cal N}_{|(0,0,0)}) 
(0,\bS, \bal)\,dt \right).
\]
We shall prove that for $\e$ small enough, ${\cal P}$ is a contraction 
on a small ball in $\tilde{\cal M}$ of radius $C_0 \e \rho^2$, where 
$C_0$ is twice the constant defined in Lemma~\ref{le:6.1}. This will give
the existence of a solution of  ${\cal N}(w)=0$ in the space  $\tilde{\cal M}$.

Let $w_1$, $w_2$  be two elements in this ball. We begin by estimating
the image by $\tLaa^{-1}$ of
\[
I(w_1, w_2) \equiv (D{\cal N}_{|(v_2,\bS_2,\bal_2)}-D{\cal N}_{|(0,
\bS_2,\bal_2)} ) (v_2,0,0) - (D{\cal N}_{|(v_1,\bS_1,\bal_1)}-
D{\cal N}_{|(0,\bS_1,\bal_1)} ) (v_1,0,0)  .
\]
\begin{lemma}
There exists a constant $c > 0$ (depending on $C_0$) such that 
\[
||\tLaa^{-1} I(w_1,w_2)||_{\tilde{\cal M}} \leq c 
\rho^{2-\nu} ||w_2 -w_1||_{\tilde{\cal M}} .
\]
\end{lemma}
{\bf Proof :} From the structure of ${\cal N}(w)$ we see that 
\[   
(D{\cal N}_{|(v,\bS,\bal)}-D{\cal N}_{|(0,\bS,\bal)} ) (v,0,0)  
\]
\[
= \frac{N(N+2)}{4} 
((\bar{u}_{\be} (\bR +\bS, \ba +\bal, \cdot) +v)^{\frac{4}{N-2}}
-(\bar{u}_{\be} (\bR +\bS, \ba +\bal, \cdot))^{\frac{4}{N-2}}) v.
\]
We write
\[
I(w_1,w_2) (x)= \frac{N(N+2)}{4} (I_1 +I_2 +I_3),
\]
where
\[
I_1= ((\bar{u}_{\be} (\bR +\bS_2, \ba +\bal_2, x) +v_2(x))^{\frac{4}{N-2}}
-(\bar{u}_{\be} (\bR +\bS_2, \ba +\bal_2, x))^{\frac{4}{N-2}}) (v_2 -v_1)(x)
\]
\[
= \frac{4}{N-2} v_2(x) (v_2 -v_1)(x) \int_0^1(\bar{u}_{\be}(\bR +\bS_2,\ba 
+\bal_2, x) +s v_2(x))^{\frac{6-N}{N-2}}\,ds,
\]
\[
I_2 = ((\bar{u}_{\be} (\bR +\bS_2, \ba +\bal_2, x) +v_2)^{\frac{4}{N-2}}
-(\bar{u}_{\be} (\bR +\bS_1, \ba +\bal_1, x) +v_1)^{\frac{4}{N-2}}) v_1(x)
\]
\[
=\frac{4}{N-2}v_1(x)\int_0^1\left((\bar{u}_{\be} (\bR +s\bS_2+(1-s)\bS_1,
\ba + s\bal_2+(1-s)\bal_1, x) +sv_2(x) +(1-s)v_1(x)\right)^{\frac{6-N}{N-2}}
\]
\[
\times \big(\sum_i \mu_i (\bR +s\bS_2+(1-s)\bS_1, \ba +s\bal_2+(1-s)\bal_1, x)
(\bS_{2,i}-\bS_{1,i}) 
\]
\[
+\sum_{i,j}\gamma_{i}^j (\bR +s\bS_2+(1-s)\bS_1, \ba +s\bal_2+(1-s)\bal_1, x)(
\bal^j_{2,i}-\bal^j_{1,i}) +s(v_2(x) -v_1(x))\big)\,ds 
\]
and finally
\[
I_3 = - ((\bar{u}_{\be} (\bR +\bS_2, \ba +\bal_2, x))^{\frac{4}{N-2}}
-(\bar{u}_{\be} (\bR +\bS_1, \ba +\bal_1, x))^{\frac{4}{N-2}}) v_1(x)
\]
\[
=- \frac{4}{N-2}v_1(x)\int_0^1\left((\bar{u}_{\be} (\bR +s\bS_2+(1-s)\bS_1,
\ba +s\bal_2+(1-s)\bal_1, x)\right)^{\frac{6-N}{N-2}}
\]
\[
\big(\sum_i \mu_i(\bR +s\bS_2+(1-s)\bS_1,\ba +s\bal_2+(1-s)\bal_1, x)(
\bS_{2,i}- \bS_{1,i}) 
\]
\[
+\sum_{i,j} \gamma_{i}^j (\bR +s\bS_2+(1-s)\bS_1, \ba +s\bal_2+(1-s)\bal_1, x)(
\bal^j_{2,i}-\bal^j_{1,i})\big)\,ds .
\]
\begin{itemize}
\item
First consider the region where $\mbox{dist}(x,\Sigma) \leq 
\e^{\frac{2}{N-2}}$. Here we use the fact that for some constant $c>0$ 
depending only on $N$, we have the estimates
\[
|\mu_i (\bR +\bS, \ba +\bal) |\leq c\,\bar{u}_{\be}
(\bR+\bS, \ba+\bal) \quad \mbox{and}\quad |\gamma_i^j(\bR +\bS, \ba +\bal)| \leq c \,
\mbox{dist} (x, \Sigma) \bar{u}_{\be} (\bR+\bS, \ba+\bal),
\]
where the constant $c>0$ depends on $C_0$. These yield
\[
|I_1| \leq c\, \e \rho^{2-\nu} \sup (1, \e^{\frac{6-N}{N-2}}) 
||v_2-v_1||_{2, \alpha, \nu,2}  \mbox{dist}(x,\Sigma)^{\frac{N-6}{2} +2 \nu}
\]
\[
|I_2| \leq  c\, \e  \rho^{2-\nu} ( \mbox{dist}(x,\Sigma)^{\nu -2}
(|\bS_2 -\bS_1| +\mbox{dist}(x,\Sigma) |\bal_2 -\bal_1 | ) 
\]
\[
+\sup (1, \e^{\frac{6-N}{N-2}}) \mbox{dist}(x,\Sigma)^{\frac{N-6}{2} +2\nu } 
||v_2 -v_1||_{2, \alpha, \nu,2})
\]
and
\[
|I_3| \leq c\,\e \rho^{2-\nu}(\mbox{dist}(x,\Sigma)^{\nu -2}(|\bS_2 -
\bS_1| +\mbox{dist}(x,\Sigma) |\bal_2 -\bal_1|).
\]
\item
Next, assume that $\mbox{dist}(x,\Sigma) \in [ \e^{\frac{2}{N-2}},1]$.
In this case we use the estimates
\[
|\mu_i (\bR +\bS, \ba +\bal) |\leq c\, \e \mbox{dist}(x,\Sigma)^{2-N}
 \quad \mbox{and}\quad |\gamma_i^j(\bR +\bS, \ba +\bal)| 
\]
\[
\leq c \, \e
 \mbox{dist} (x, \Sigma)  1+ \e^{\frac{4}{N-2}} \mbox{dist}(x,\Sigma)^{-N}),
\]
where the constant $c>0$ depends on $C_0$.
Then as above, we derive 
\[
|I_1| \leq c\,\e^{\frac{4}{N-2}} \rho^{2-\nu}||v_2-v_1||_{2, \al,\nu,2}
\mbox{dist}(x,\Sigma)^{N-6 +2 \nu},
\]
\[
|I_2| \leq c\, \e \rho^{2-\nu}(\e^{\frac{4}{N-2}}\mbox{dist}(x,
\Sigma)^{\nu -4} (|\bS_2 -\bS_1|+ 
(\mbox{dist}(x,\Sigma)^{N-1} +\e^{\frac{4}{N-2}} \mbox{dist}(x,\Sigma)^{-1})|\bal_2 -\bal_1|)
\]
\[ 
+\e^{\frac{6-N}{N-2}} \mbox{dist}(x,\Sigma)^{N-6 +2\nu}||v_2 -v_1||_{2,\al,
\nu,2})
\]
and
\[
|I_3| \leq c\,\e^{\frac{N+2}{N-2}} \rho^{2-\nu}(\mbox{dist}(x,\Sigma)^{
\nu -4} (|\bS_2 -\bS_1|+(\mbox{dist}(x,\Sigma)^{N-1}+\e^{\frac{4}{N-2}} \mbox{dist}(x,\Sigma)^{-1}) |\bal_2 -\bal_1 | ).
\]
\item
Finally, if $\mbox{dist}(x,\Sigma) \geq 1$, we obtain
\[
|I_1| \leq c\,\e^{\frac{4}{N-2}}C_0\rho^{2-\nu}||v_2-v_1||_{2,\al,\nu,2}
\mbox{dist}(x,\Sigma)^{-N-2},
\]
\[
|I_2| \leq c\,\e  \rho^{2-\nu}(\e^{\frac{4}{N-2}}
|\bS_2 -\bS_1|+\e^{\frac{6-N}{N-2}} ||v_2 -v_1||_{2, \al,\nu,2}) \mbox{dist} (x, \Sigma)^{-N-2}
\]
and
\[
|I_3| \leq c\,\e^{\frac{N+2}{N-2}}C_0\rho^{2-\nu}|\bS_2 -\bS_1| \mbox{dist}(x,\Sigma)^{
-N-2} .
\]
\end{itemize}
Combining these gives
\[
||I(w_1,w_2)||_{0,\al,\nu-2,-2} \leq  c\, \rho^{2-\nu}(\e^{\frac{4}{
N-2}}||v_2 -v_1||_{2,\al,\nu,2} +\e (|\bS_2-\bS_1| +
\e^{\frac{2}{N-2}}|\bal_2 -\bal_1|)).
\]
In particular 
\[
||I ( w_1, w_2)||_{0, \alpha, \nu-2, -2} 
\leq  c\, \rho^{2-\nu}(\rho^\nu||v_2 -v_1||_{2,\al,\nu,2} 
+\e |\bS_2-\bS_1| +\e \rho|\bal_2 -\bal_1|)
\]
\[
= c \, \rho^{2-\nu} ||w_2 -w_1||_{\tilde{\cal M}}.
\]
The estimate of the Lemma now follows from the boundedness of
$\Laa^{-1}$, as proved in Proposition~\ref{pr:7.1}.
\hfill $\Box$
\medskip

Next shall estimate the image by $\tLaa^{-1}$ of
\begin{equation}     
II(w_1, w_2) = (D{\cal N}_{|(0,\bS_2,\bal_2)}-D{\cal N}_{|(0,0,0)})
(v_2,0,0)  - (D{\cal N}_{|(0,\bS_1,\bal_1)}-D{\cal N}_{|(0,0,0)})(v_1,0,0).
\label{eq:9.3}
\end{equation}
\begin{lemma}
There exists some constant $c>0$ (depending on $C_0$) such that 
\[
||\tLaa^{-1} II(w_1,w_2)||_{\tilde{\cal M}} \leq c 
\rho^{2-\nu} ||w_2 -w_1||_{\tilde{\cal M}}.
\]
\end{lemma}
{\bf Proof :}
Recall that
\[
(D{\cal N}_{|(0,\bS,\bal)}-D{\cal N}_{|(0,0,0)} ) (v,0,0)  
= \frac{N(N+2)}{4} 
((\bar{u}_{\be} (\bR +\bS, \ba +\bal, \cdot ))^{\frac{4}{N-2}}
-(\bar{u}_{\be} (\bR, \ba, \cdot))^{\frac{4}{N-2}}) v.
\]
Therefore, we may decompose $II$ as 
\[
II = \frac{N(N+2)}{4} (II_1 +II_2),
\]
where
\[
II_1= ((\bar{u}_{\be} (\bR +\bS_2, \ba +\bal_2, x))^{\frac{4}{N-2}}
-(\bar{u}_{\be} (\bR, \ba, x))^{\frac{4}{N-2}}) (v_2 -v_1)(x)
\]
\[
=\frac{4}{N-2}(v_2 -v_1)(x)\int_0^1 ((\bar{u}_{\be}(\bR +s\bS_2,\ba +s\bal_2, 
x))^{\frac{6-N}{N-2}}
\]
\[
(\sum_i \mu_i (\bR +s\bS_2, \ba +s\bal_2, x)\bS_{2,i} +\sum_{i,j} 
\gamma_{i}^j (\bR +s\bS_2, \ba +s\bal_2, x)\bal^j_{2,i})\,ds) 
\]
and 
\[
II_2 = ((\bar{u}_{\be} (\bR +\bS_2, \ba +\bal_2, x))^{\frac{4}{N-2}}
-(\bar{u}_{\be} (\bR +\bS_1, \ba +\bal_1, x) )^{\frac{4}{N-2}}) v_1(x)
\]
\[
 =\frac{4}{N-2} v_1 (x) \int_0^1 ((\bar{u}_{\be} (\bR +s\bS_2+(1-s)\bS_1, \ba 
+s\bal_2+(1-s)\bal_1, x))^{\frac{6-N}{N-2}}
\]
\[
(\sum_i \mu_i (\bR +s\bS_2+(1-s)\bS_1, \ba +s\bal_2+(1-s)\bal_1, x)
(\bS_{2,i}-\bS_{1,i}) 
\]
\[
+\sum_{i,j} \gamma_{i}^j (\bR +s\bS_2+(1-s)\bS_1, \ba +s\bal_2+(1-s)\bal_1, 
x)(\bal^j_{2,i}-\bal^j_{1,i}))\,ds .
\]
\begin{itemize}
\item
As before, first consider the region where $\mbox{dist}(x,\Sigma) \leq 
\e^{\frac{2}{N-2}}$. As in the previous Lemma, we use the fact that there 
exists a constant $c>0$ only depending on $N$ such that
\[
|\mu_i (\bR +\bS, \ba +\bal) |\leq c\, \bar{u}_{\be} (\bR+\bS,\ba+\bal) 
\quad \mbox{and}\quad |\gamma_i^j| \leq c\,  \mbox{dist} (x, \Sigma) 
\bar{u}_{\be} (\bR+\bS, \ba+\bal).
\]
So, we get the estimate
\[
|II_1| \leq   c\, \mbox{dist}(x,\Sigma)^{\nu -2} (C_0 \rho^2 +\mbox{dist}(x,
\Sigma) C_0 \rho ) ||v_2 -v_1||_{0, \alpha, \nu,2},
\]
\[
|II_2| \leq   c\,   \e \rho^{2-\nu} \mbox{dist}(x,\Sigma)^{\nu -2} 
(|\bS_2 -\bS_1| +\mbox{dist}(x,\Sigma) |\bal_2 -\bal_1 | ). 
\]
\item
Next, assume that $\mbox{dist}(x,\Sigma) \in [ \e^{\frac{2}{N-2}},1]$, then 
as above, we can derive the estimate
\[
|II_1| \leq   c\, \e^{\frac{4}{N-2}} \mbox{dist}(x,\Sigma)^{\nu -4} ( \rho^2 
+(\mbox{dist}(x,\Sigma)^{N-1} +\e^{\frac{4}{N-2}}\mbox{dist}(x,\Sigma)^{-1}) 
\rho ) ||v_2 -v_1||_{0, \alpha, \nu,2} ,
\]
\[
|II_2| \leq   c \,  \rho^{2-\nu} \e^{\frac{N+2}{N-2}} \mbox{dist}(x,
\Sigma)^{\nu -4} (|\bS_2 -\bS_1| +(\mbox{dist}(x,\Sigma)^{N-1}+
\e^{\frac{4}{N-2}}\mbox{dist}(x,\Sigma)^{-1}) |\bal_2 -\bal_1 | ) .
\]
\item
And finally, if $\mbox{dist}(x,\Sigma) \geq 1$, we obtain
\[
|II_1| \leq   c \e^{\frac{4}{N-2}} \mbox{dist}(x,\Sigma)^{-N-2} \rho^2 
||v_2 -v_1||_{0, \alpha, \nu,2},
\]
\[
|II_2| \leq   c  C_0  \rho^{2-\nu} \e^{\frac{N+2}{N-2}} \mbox{dist}(x,
\Sigma)^{-N-2} |\bS_2 -\bS_1| .
\]
\end{itemize}
This allows us to conclude that
\[
 || II(w_1,w_2) ||_{0, \alpha, \nu-2, -2}\leq
c C_0 ( \rho^2 ||v_2-v_1||_{2, \alpha, \nu,2}
+ \rho^{2-\nu} \e (|\bS_2 -\bS_1| +\e^{\frac{2}{N-2}} |\bal_2 -\bal_1 | )),
\]
which leads to 
\[
|| II(w_1,w_2) ||_{0, \alpha, \nu-2, -2}\leq
c C_0 \rho^{2-\nu}  ||w_2-w_1||_{\tilde{\cal M}}.
\]
Then the estimate of the Lemma follows from the result of 
Proposition~\ref{pr:7.1}.
 \hfill $\Box$

We now turn to the estimate of the image by $\tLaa^{-1} $ of
\[
III (w_1, w_2) \equiv  (D{\cal N}_{|(0, \bS_2, \bal_2)} -D{\cal N}_{
|(0,0,0)}) (0,\bS_2, \bal_2)
-(D{\cal N}_{|(0, \bS_1, \bal_1)} -D{\cal N}_{|(0,0,0)}) (0,\bS_1, \bal_1)
\]
As before we decompose $III =III_1 +III_2$ where
\[
III_1 \equiv  (D{\cal N}_{|(0, \bS_2, \bal_2)} -D{\cal N}_{|(0,0,0)})
(0,\bS_2-\bS_1, \bal_2-\bal_1)
\]
\[
III_2 = (D{\cal N}_{|(0, \bS_2, \bal_2)}  -D{\cal N}_{|(0, \bS_1, \bal_1)} )
(0,\bS_1, \bal_1)
\]
To estimate either of these terms, it is convenient to consider
them as special cases of the more general quantity
\[
\overline{III}(\bS,\bS',\bS'', \bal,\bal',\bal'') = (D{\cal N}_{|(0, 
\bS'+\bS, \bal'+\bal)}  -D{\cal N}_{|(0, \bS, \bal)} ) (0,\bS'', \bal'')
\]
We first prove the Lemma~:
\begin{lemma}
There exists some constant $c>0$ (depending on $C_0$) such that
\[
||\tLaa^{-1} \overline{III}(\bS,\bS',\bS'', \bal,\bal',\bal'') 
||_{\tilde{\cal M}} \leq 
c \e (|\bS'||\bS''| +\rho^2 |\bal'||\bal''| + \rho |\bS'||\bal''| + \rho
|\bal'||\bS''|).
\]
\end{lemma}
{\bf Proof :}
We may write 
\[
\overline{III}(\bS,\bS',\bS'', \bal,\bal',\bal'') 
=\int_{0}^{1}D^2 {\cal N}_{|(0, \bS +s\bS', \bal +s\bal')} ((0, \bS', \bal'),
(0, \bS'', \bal''))\,ds,
\]
and shall concentrate on estimating the integrand
\begin{equation}
\overline{iii}(\bS,\bS',\bS'', \bal,\bal',\bal'') \equiv D^2 {\cal N}_{
|(0, \bS, \bal)}((0, \bS', \bal'),(0, \bS'', \bal'')).
\label{eq:9.5}
\end{equation}
We take advantage from the observation that, for all $\bR \in {\Bbb R}^n$ 
and for all $\ba \in ({\Bbb R}^N)^n$
\[
\Delta \bar{u}_{\be} (\bR,\ba,x) +\frac{N(N-2)}{4}(\bar{u}_{\be}
(\bR,\ba,x))^{\frac{N+2}{N-2}} =0,
\]
if $\mbox{dist}(x,x_i) \leq \rho_i$.
Therefore, if $\mbox{dist} (x,x_i)\leq \rho_i$, (\ref{eq:9.5}) is 
identically $0$.
Moreover, if $\mbox{dist} (x,x_i)\geq 2\rho_i$, we observe that 
$\bar{u}_{\be} (\bR,\ba,x)=\bar{w}_{\be}(\bR, x) $ is harmonic and does not 
depend on $\ba$. Therefore, in this case,  $\overline{iii}$ reduces to
\[
\overline{iii} = \frac{N(N-2)}{4}\sum_{i,j} \del^2_{R_i,R_j}((\bar{w}_{\be} (\bR + \bS)^\p)S_{i}'S_{j}'',
\]
And, still assuming that $\mbox{dist} (x,x_i)\geq 2\rho_i$, we find the estimate
\[
|\overline{iii}| \leq c \e^{\frac{N+2}{N-2}} \mbox{dist}(x, \Sigma)^{-N-2} |\bS'| |\bS''|.
\]
Finally, if  $\mbox{dist}(x,x_i) \in [\rho_i,2\rho_i]$, then by 
definition 
\[
\buex=\sum_i \chi_i(x-x_i)\uix +\bwx (1-\sum_i\chi_i(x-x_i)),
\] 
so we can write 
\[
\Delta \buex +\frac{N(N-2)}{4}(\buex)^{\frac{N+2}{N-2}}
\]
\[
= \frac{N(N-2)}{4} (-\chi_i (\ui)^{\p}+ (\chi_i \ui +(1-\chi_i)\bw)^{\p}) 
\]
\[
+\Delta (\chi_i ) (\ui -\bw)+2\nabla (\chi_i)  \nabla (\ui-\bw).
\]

Differentiating twice with respect to both $\bR$ and $\ba$, it is easy to 
derive the estimate
\[
| D^2_{|(\bR+\bS,\ba+\bal)} (-\chi_i (u_{\e_i})^{\p}+ (\chi_i u_{\e_i} +
(1-\chi_i)\bar{w}_{\be} )^{\p} ) (\bS', \bal') (\bS'', \bal'')|
\]
\[
\leq c \e^{\p} ( \rho^{-4}
 (\rho^{2-N} |\bS'| |\bS''| +\rho (|\al'||\bS''|+|\al''||\bS'|) +\rho^{2} 
|\al'||\al''|)
\]
\[
+\rho^{N-6} (\rho^{2-N} |\bS'| +\rho |\bal'| )(\rho^{2-N} |\bS''| +\rho 
|\bal''| ))
\]
\[
\leq c \e^{\p} (\rho^{-2-N} |\bS'||\bS''| +\rho^{-3} (|\al'||\bS''|+
|\al''||\bS'|) +\rho^{-2} |\al'||\al''|).
\]
Using the result of Proposition~\ref{pr:4.1} and differentiating 
twice with respect to both $\bR$ and $\ba$, we get the expansion
\[
D^2_{|(\bR+\bS, \ba+\bal) }(u_{\e_{i_0}} -\bwx)(\bS', \bal')(\bS'', \bal'')
\]
\[
=D^2_{|(\bR+\bS, \ba +\bal)} (\frac{\e_{i_0}}{2} R_{i_0}^{\frac{2-N}{2}} - 
\sum_{i\neq i_0} \frac{\e_{i}}{2} R_i^{\frac{N-2}{2}} |x_{i_0}-x_i|^{2-N})
(\bS' ,\bal')(\bS'', \bal'')
\]
\[
+\frac{N-2}{2} D_{|(\bR+\bS, \ba+\bal)}^2 ((\e_{i_0} R_{i_0}^{\frac{2-N}{2}} 
a_{i_0} +\sum_{i\neq i_0}\frac{\e_{i}}{2} R_i^{\frac{N-2}{2}} 
|x_{i_0} -x_i|^{-N}(x_{i_0}-x_i))\cdot(x-x_{i_0}))(\bS',\bal')(\bS'', \bal'')
\]
\begin{equation}
 +O(\e \rho^2)(|\bS'||\bS''| + |\bal' ||\bal''| +|\bS'||\bal''| + 
|\bal'| |\bS''|).
\label{eq:7.9}
\end{equation}
and
\[
D_{|(\bR +\bS, \ba +\bal)}^2 (\del_{r_{i_0}}(u_{\e_{i_0}}
(R_{i_0},a_{i_{0}},x-x_{i_0}) - \bwx))(\bS', \bal')(\bS'', \bal'')
\]
\[
 =\frac{N-2}{2}
 D^2_{|(\bR +\bS,\ba+\bal)} ( (\e_{i_0} R_{i_0}^{\frac{2-N}{2}} a_{i_0} 
-\sum_{i\neq i_0}\frac{\e_i}{2} R_i^{\frac{N-2}{2}} |x_{i_0} -x_i|^{-N}
(x_{i_0}-x_i))\cdot \frac{x-x_{i_0}}{|x-x_{i_0}|})(\bS',\bal')(\bS'', \bal'') 
\]
\begin{equation}
 +O(\e \rho)(|\bS'||\bS''| + |\bal' ||\bal''| +|\bS'||\bal''| + 
|\bal'| |\bS''|).
\label{eq:7.10}
\end{equation}
This expansion allows us to estimate
\[
D_{|(\bR+\bS, \ba+\bal)}^2 (\Delta (\chi_i ) (\ui -\bw)+
\]
\[
2\nabla (\chi_i)  
\nabla (\ui-\bw))(\bS',\bal')(\bS'', \bal'') .
\]
As in the setting of Lemma~\ref{le:6.2}, we shall decompose this quantity 
as $h_0 +h' +h''$ as before, then we get the estimates
\[
|h_0|\leq c \,
 \e (\rho^{-2} |\bS'||\bS''| + |\bal' ||\bal''| +|\bS'||\bal''| + |\bal'|
 |\bS''|).
\]
\[
|h'|\leq c \,
 \e (\rho^{-1} |\bS'||\bS''| + |\bal' ||\bal''| +\rho^{-1}|\bS'||\bal''| 
+\rho^{-1} |\bal'| |\bS''|).
\]
and finally
\[
|h''|\leq  c   \e ( |\bS'||\bS''| + |\bal' ||\bal''| +|\bS'||\bal''| +|\bal'| 
|\bS''|).
\]
Therefore, using the result of Proposition~\ref{pr:7.1} as well as the result 
of Lemma~\ref{le:6.2} we obtain
\[
||\tLaa^{-1} \overline{iii}(\bS,\bS',\bS'', \bal,\bal',\bal'') 
||_{\tilde{\cal M}} 
\leq c 
 \e (|\bS'||\bS''| +\rho^2 |\bal' ||\bal''| +\rho |\bS'||\bal''| + 
\rho |\bal'| |\bS''|).
\]
As a consequence we get that
\[
||\tLaa^{-1} \overline{III}(\bS,\bS',\bS'', \bal,\bal',\bal'') 
||_{\tilde{\cal M}} \leq c \e (|\bS'||\bS''| +\rho^2 |\bal' ||\bal''|
 +\rho |\bS'||\bal''| + \rho |\bal'| |\bS''|).
\]
\hfill $\Box$

As a corollary of this last lemma, recalling that $\e |\bS_j| +
\e \rho |\bal_j| \le C_0 \e \rho^2$, we get that
\[
||\tLaa^{-1} III ( w_1,w_2) ||_{\tilde{\cal M}} \leq c
 \e ( \rho^2 |\bS_2-\bS_1'| +\rho^3 |\bal_2 -\bal_1| )\leq c \rho^2 
||w_2-w_1||_{\tilde{\cal M}}.
\]

Collecting the results of these lemmata, we have established
that for some constant $c > 0$ depending on $C_0$, 
\[
||{\cal P}(w_2) -{\cal P}(w_1)||_{\tilde{\cal M}} \leq  c \rho^{2-\nu} 
||w_2-w_1||_{\tilde{\cal M}}.
\]
By taking $\e$, and hence $\rho$, sufficiently small, we have shown
that the map ${\cal P}$ is a contraction on the ball of radius
$C_0 \e \rho^2$ in $\tcM$. Hence it has a unique fixed point $w$,
and this function $w = (v,\bS,\bal)$ is a solution of the
equation (\ref{eq:9.1}). It is clear that $\bar{u}_{\be}(\bR + \bS,
\ba + \bal,\cdot) + v$ is positive near the points of $\Sigma$, hence
by the maximum principle is positive everywhere.

This completes the existence of the solution promised in 
Theorem~\ref{th:1.1}. 

\section{The nondegeneracy of the solutions}

We now show that for $\e$ sufficiently small, the solutions
we have constructed above are nondegenerate in the sense defined 
in \cite{MPU1}, \cite{MPU2}.  Actually, there are two 
closely related notions of nondegeneracy, the definitions of which
we now recall, in terms of the notations of this paper:

\begin{definition} Let $g_0$ be the standard metric on $S^N$.  A metric 
$g = u^{\frac{4}{N-2}}g_0$ of constant positive scalar curvature on 
$S^N \setminus \Lambda$, as well as the corresponding conformal factor $u$, 
is called {\it marked nondegenerate} if the linearization of the scalar 
curvature operator ${\cal N}$ for $g_0$ about the solution $u$
is injective on the function space ${\cal C}^{2,\al}_{\mu}(S^N \setminus
\Lambda)$ for all $\mu > (2-N)/2$, or equivalently, if the linearization
of ${\cal N}$ relative to the metric $g$ about the constant solution $1$
is injective on the function space  ${\cal C}^{2,\al}_{\mu'}(S^N \setminus 
\Lambda)$ (defined with respect to the metric $g$) for all $\mu' > 0$. 
The metric 
$g$, or solution $u$, is called {\it unmarked nondegenerate} if the 
linearization is injective for all $\mu > (4-N)/2$, or equivalently, for 
all $\mu' > 1$. 
\end{definition}

These two nondegeneracy conditions are precisely what is needed to ensure 
the smoothness of the marked and unmarked moduli spaces ${\cal M}_{\Lambda}$ 
and ${\cal M}_n$ at $g$. The former of these spaces is the space of all
metrics of constant positive scalar curvature on the complement of the 
finite set $\Lambda$ in the sphere $S^N$, while the latter is the set 
of all such metrics on the complement of any finite set $\tilde{\Lambda}$
of cardinality $n$. As proved in \cite{MPU1} and \cite{MPU2}, these moduli 
spaces are real  analytic sets, hence are stratified and may be written 
as the union of smooth, real analytic manifolds of varying dimensions. The 
existence of one smooth point in a given component shows that the top 
dimensional stratum in that component is of the dimension predicted by the
formal, index-theoretic, calculations, namely $n$ for the marked spaces
and $n(N+1)$ for the unmarked ones. 

In this section we shall prove
\begin{proposition} The solutions constructed here are unmarked
nondegenerate. For generic configurations $\Lambda$ they are also
marked nondegenerate. 
\end{proposition}

As a corollary, we obtain smoothness of the unmarked moduli spaces 
without any restriction (when $\e$ is small) and of the marked
moduli spaces for generic configurations. In fact, near a generic 
$\Lambda$, and when $\e$ is sufficiently small, we may use $(p_1, 
\ldots, p_n, \e_1, \ldots, \e_n) = (\bar{p},\be) \in (S^N)^n \times 
(\R^+)^n$ as coordinates on the unmarked moduli space ${\cal M}_n$, 
and $\be$ as coordinates on the marked moduli spaces.  We note 
that marked nondegenerate solutions (which are {\it a fortiori} also 
unmarked nondegenerate) were constructed in \cite{MPU2} for certain very 
special configurations $\Lambda$, which in particular contain only even 
numbers of points; the Delaunay parameters of those solutions are not 
precisely prescribed, but they need not be close to zero. Even when $n$ 
is even, it is not clear that these solutions lie in the same component 
of ${\cal M}_n$ as the ones we construct above. 

We first demonstrate  the unmarked nondegeneracy; the proof is
by contradiction. By a slight 
change of notation from the rest of the paper, we consider the solution 
$u$ on $S^N \setminus \Lambda$ to have the form $u=\bue +v$, where the 
approximate solution on the sphere (rather than on $\R^N$) is now denoted 
$\bue$, and $v$ is an element of ${\cal C}^{2,\al}_{\nu}$ with $1 < \nu < 2$,
with norm in this space bounded by $C_0 \e \rho^{2-\nu}$.
Assume that for some sequence of $\e_k$ tending to $0$, the linearized
operator
\[
{\cal L}_k = \Delta_{S^N} - \frac{N(N-2)}{4} + \frac{N(N+2)}{4} 
u_k^{\frac{4}{N-2}} 
\]
is not injective on ${\cal C}^{2,\al}_{\mu}(S^N \setminus \Lambda)$
for some $\mu > (4-N)/2$. Here we have denoted by $u_k$ the solution 
$\bar{u}_{\be_k}(\bR_k,\ba_k,\cdot) + v_k$. Then there is some element
$w_k \in {\cal C}^{2,\al}_{\mu}$ such that ${\cal L}_k w_k=0$.

First normalize $w_k$, multiplying it by a suitable constant, so that
$\sup d(y)^{-\mu}|w_k(y)| = 1$, where $d(y)$ is the distance of the
point $y$ from $\Lambda$, say in the spherical metric.  Choose a point 
$y_k \in S^N \setminus \Lambda$ realizing this supremum, i.e. such
that $d(y_k)^{-\mu}|w_k(y_k)| = 1$.  As $k$ tends to infinity,
the function $v_k$ tends to zero in ${\cal C}^{2,\al}_{\nu}$, since its 
norm is dominated by $C_0\e_k \rho_k^{2-\nu}$, and $\bar{u}_{\be_k}$
converges uniformly to zero on compact subsets of $S^N \setminus \Lambda$. 
In addition, we can assume that the sequences $\bR_k$ and $\ba_k$ converge 
to some fixed $\bR^*$ and $\ba^*$, respectively. 

Suppose first that some subsequence of the $y_k$ converges
to a point $y_0 \in S^N \setminus \Lambda$. Then we can extract a
subsequence of the $w_k$ which converge (in ${\cal C}^\infty$) to a 
limiting function $w$ on $S^N \setminus \Lambda$; $w$ must be nontrivial 
since $d(y_0)^{-\mu}|w(y_0)|=1$. Furthermore, $Lw \equiv \left(\Delta_{S^N} 
- (N-2)^2/4\right)w = 0$, and $|w(y)| \le d(y)^{\mu}$. Since $\mu > (4-N)/2$,
it is standard that $w$ is a weak solution of $Lw = 0$ on all of $S^N$, 
hence also a strong solution. But clearly $L$ has only
trivial nullspace, hence $w \equiv 0$, which is a contradiction.

If, on the other hand, some subsequence of the $y_k$ converge
to one of the points $p_{i_0} \in \Lambda$, then we must argue
somewhat differently.  Choose a function $A$ on $S^N \setminus
\{p_{i_0},q\}$, where $q \notin \Lambda$, which transforms
this twice-punctured sphere to the cylinder $\R \times S^{N-1}$
and such that the standard spherical metric $g_S$ and the
(dilated) cylindrical metric $g_C$ are related by
\[
g_S = A^{\frac{4}{N-2}} g_C, \qquad g_C = \frac{N-2}{N}(dt^2 + d\theta^2). 
\]
On $C$, the function $A$ is simply a multiple of $(\cosh t)^{\frac{2-N}{2}}$,
and on $S^N$ is of the order $d(y,p)^{\frac{2-N}{2}}$, for $p = p_{i_0}$ 
or $q$.  Noting that $g_C$ has scalar curvature  
$N(N-1)$, the same as $g_S$, the conformal Laplacians of these two metrics, 
$L_C$ and $L_S$ satisfy the usual transformation rule
\begin{equation}
L_C (A\phi) = A^{\frac{N+2}{N-2}}L_S \phi,
\label{eq:35.1}
\end{equation}
for any function $\phi$. 
The solution $u_k$ on $S^N$ corresponds to a solution $Au_k$ on $C$.
It is easy to check that the linearizations ${\cal L}_{S,k}$ of the 
scalar curvature operators on $S^N$ at $u_k$ and ${\cal L}_{C,k}$ on 
$C$ at $A u_k$ satisfy the same transformation rule (\ref{eq:35.1}) as 
$L_S$ and $L_C$, cf. \cite{MPU2}. (Unlike 
the transformation rule for the conformal Laplacian, this holds
only because we are dealing with two metrics of the same constant
scalar curvature.) 

Because of these machinations, we may replace $w_k$ by a solution
$\tilde{w}_k$ of ${\cal L}_{C,k}\tilde{w}_k = 0$. For convenience,
we call ${\cal L}_{C,k}$ simply ${\cal L}_k$ and $\tilde{w}_k$ simply
$w_k$ again, and let $y = (t,\theta)$ denote the variable
on the cylinder. Define $\mu' = \mu + (N-2)/2$, so that $\mu' > 1$.
These new functions $w_k$ satisfy
\begin{equation}
\sup d(y)^{-\mu'}|w_k| = 1,
\label{eq:35.2}
\end{equation}
where $d(y)$ is some function equalling the distance to the set of 
other singular points $\Lambda' = \Lambda \setminus \{p_{i_0}\}$
(transplanted to the cylinder) in a neighbourhood of this set, and equalling
$\mbox{sech\,}t$ outside this neighbourhood. Because of (\ref{eq:35.2}),
$w_k$ decays at either end of the cylinder. 

As before, let $y_k = (t_k,\theta_k)$ denote the point on $C$ where 
the supremum in (\ref{eq:35.2}) is attained. We already are assuming 
that $t_k \rightarrow \infty$. By translating back by $t_k$ and
multiplying by a suitable constant, we find yet another sequence
of solutions, which we again call $w_k$, attaining their maximum
at $t = 0$, and which solve the translated equation, which we
again write as ${\cal L}_k w_k = 0$. Here ${\cal L}_k$ is the
linearized scalar curvature operator relative to a metric 
$g_k = (A(\bar{u}_{\be_k} + v_k))^{4/(N-2)}g_C$ which is singular at a 
finite collection of points which are translating toward $t = -\infty$. 
Since $v_k$ is tending to zero, $g_k$ more and more nearly approximates
the Delaunay metric $V_{\e_k}^{4/(N-2)}g_C$, where $V_{\e_k}$ is
one of the ($\theta$-independent) Delaunay solutions considered in \S 2. 

It is not hard to see that the $w_k$ converge to a nontrivial
solution $w$ of the limiting equation ${\cal L} w = 0$, and
that $w$ is bounded by $e^{-\mu' t}$ for all $t$. 

There are two subcases. In the first, $V_k(0,\theta)$ tends to zero. 
Then $w$ satisfies the equation 
\[
\frac{N}{N-2}\left(\del_t^2 + \Delta_{S^{N-1}}\right)w - \frac{N(N-2)}{4}w 
= \frac{N}{N-2}\left( \del_t^2 + \Delta_{S^{N-1}} - \frac{(N-2)^2}{4}
\right) w = 0.
\]
By decomposing $w$ into eigencomponents with respect to $\Delta_{S^{N-1}}$,
we see that any eigencomponent $w_j$ is a sum of exponentials,
$w_j = a_j^+ e^{\gamma_j t} + a_j^- e^{- \gamma_j t}$. Since $w$
decays as $t \rightarrow +\infty$, $a_j^+ = 0$. But then it is
clear that no function of the form $e^{-\gamma_j t}$ can be bounded
for all $t$ by $e^{-\mu' t}$ unless $\mu' = \gamma_j$, which is not
the case, so we arrive at a contradiction.
In the second, $V_k(0,\theta)$ does not tend to zero. By translating
by a fixed finite amount, we may assume that $V_k$ tends to the
function $(\cosh t)^{(2-N)/2}$, and hence, after pulling out the
superfluous constants, that the limiting function $w$ satisfies
\[
\left(\del_t^2 + \Delta_{S^{N-1}} - \frac{(N-2)^2}{4} + \frac{N^2 - 4}{4}
\mbox{sech\,}^2 t \right) w = 0.
\]
Again separate $w$ into its eigencomponents $w_j$. Then
\[
\del_t^2 w_j - \left(\frac{(N-2)^2}{4} + \lambda_j\right) w_j +
\frac{N^2 - 4}{4} \mbox{sech\,}^2 t\, w_j = 0.
\]
For $j = 0$ the indicial roots of this equation at both $\pm \infty$
are $\pm (N-2)/2$, for $j = 1, \ldots, N$ they are $\pm N/2$, and for 
$j > N$ they all satisfy $|\gamma_j^{\pm}| \ge (N+2)/2$. 

The components $w_j$ with $j > N$ are easy to eliminate. In 
fact, these $w_j$ must decay faster than $e^{\pm (N+2)|t|/2}$
at $\pm \infty$, so we may multiply the equation satisfied by $w_j$ 
and integrate by parts to obtain
\[
\int_{-\infty}^{\infty} (\del_t w_j)^2 + \left(\lambda_j
+ \frac{(N-2)^2}{4}\right) w_j^2  - \frac{N^2 - 4}{4}\mbox{sech\,}^2 
t\, w_j^2 \,dt =0.
\]
Since $\lambda_j \ge 2N$, the integrand is nonnegative, hence $w_j = 0$. 

For the remaining cases, when $j \le N$, the indicial roots
at $\pm \infty$ are less than $(N+2)/2$ in absolute
value. On the other hand, to check unmarked nondegeneracy it
suffices to use any $\mu' > 1$ which is also less than the next omitted weight 
for the linearization about a Delaunay solution, as determined in \S 5.  
But by the results of that section, as $\e$ tends to zero, this 
next omitted weight tends to $(N+2)/2$. Thus now choose $\mu'$ in the 
range $(N/2, (N+2)/2)$; then each of the $w_j$ with $j \le N$  
decays less quickly than $e^{-\mu' t}$ as $t \rightarrow
+\infty$, which implies that these $w_j$ too must vanish.
This is a contradiction. 

This covers all cases, so we have showed that the linearization is
injective on the appropriate weighted H\"older spaces.  

This completes the proof of unmarked nondegeneracy of the solutions
when $\e$ is sufficiently small. We complete the proof of the
proposition by demonstrating the marked nondegeneracy of solutions
for generic configurations $\Lambda$. This is a simple consequence
of Sard's theorem. In fact, as discussed in \cite{MPU2}, near smooth
points of the unmarked moduli space ${\cal M}_n$, there is a
real analytic fibration $\pi: {\cal M}_n \rightarrow {\cal C}_n$
onto the configuration space of $n$ distinct points in $S^N$. 
${\cal C}_n$ is itself a real analytic manifold. Standard
differential topological arguments now imply the surjectivity
of the differential $\pi_*$ for all points in a generic fibre
of $\pi$. By a dimension count, surjectivity of this differential
is equivalent to the marked nondegeneracy of all points in
the fibre.

\section{The general case}

In this brief final section we discuss the essentially minor changes
that need to be made in order to prove the more general statement
of Theorem~\ref{th:1.1}, where the singular set is allowed to
have components of positive dimension.  As in the main body of
the paper, there are three steps. First, we must construct an 
approximate solution, or rather, a family of approximate solutions that 
become increasingly concentrated at $\Lambda$; next we must prove that 
the linearization about one of these approximate solutions is surjective 
on an appropriate function space provided the approximate solutions are 
sufficiently concentrated; finally, we perturb once again to an exact 
solution.  In each of these steps, we must somehow combine the
constructions and proofs from our previous paper \cite{MP} with
the ones here. The last step, showing that an appropriate map is
a contraction, is straightforward and we shall not comment on it
further. We now describe the first two steps.

As in the introduction, divide the components of the singular set 
$\Lambda$ into two groups, $\Lambda = \Lambda' \cup \Lambda''$, where the
first is a finite set and the latter contains all the higher dimensional 
components. The construction of a family of approximate solutions around
a component $\Lambda_j$ of dimension $k$, where $0 < k \le (N-2)/2$,
is given in detail in \cite{MP}, but briefly, it is obtained by first 
fixing a tubular neighbourhood ${\cal T}(\Lambda_j)$ 
and identifying it with a neighbourhood of the zero section in
the normal bundle $N\Lambda_j$. On each of the fibres $N_p\Lambda_j$ 
we glue in a sufficiently dilated solution of the equation
\[
\Delta u + u^{\frac{N+2}{N-2}} = 0,
\]
cut off to be supported in ${\cal T}(\Lambda_j)$. 
Because the fibres are of dimension $N-k$, this equation is subcritical,
and one can show that the radial solutions have the form $|x|^{(2-N)/2}v$, 
where $v$ tends to zero rather quickly as $|x|\rightarrow \infty$. There 
is a one-parameter family $u_\e$ of such solutions, which are all, in a 
suitable sense, dilated versions of one another.  Again, we refer
to \cite{MP} for all details.

To construct the approximate solutions when $\Lambda'$ and $\Lambda''$ are 
both nonempty, we fix Delaunay parameters $\be'$ for the points
of $\Lambda'$ and dilation parameters $\be''$ for the components
of $\Lambda''$, and set $\be = (\be',\be'')$. The elements of each of 
these subsets are mutually commensurable, i.e. $\be' = \e' \bar{q}'$ 
and $\be'' = \e'' \bar{q}''$, where $\bar{q} = (\bar{q}',\bar{q}'')$ 
is a vector with all components positive. Now construct the 
approximate solutions singular at the points of $\Lambda'$ as
in \S 3, balanced exactly as before (here we use that $|\Lambda'| \ge 2$).
This approximate solution is of size $\e'$ outside a fixed neighbourhood 
of $\Lambda'$, and so we use a partition of unity to glue it to 
the approximate solution defined in each ${\cal T}\Lambda_j$. 

The main step is to show that the linearization of the scalar curvature 
operator, $L_{\be}$, is surjective as a map 
\[
L_{\be}: {\cal C}^{2,\al}_{\nu',\nu''}(S^N \setminus \Lambda) \oplus
{\cal W} \longrightarrow {\cal C}^{0,\al}_{\nu'-2,\nu''-2}(S^N
\setminus \Lambda),
\]
at least when all the components of $\be$ are sufficiently small.
Here $1 < \nu' < 2$ is the weight parameter determining growth 
of functions in a neighbourhood of each point of $\Lambda'$ and
$\nu'' < (2-N)/2$ determines the growth near each component of $\Lambda''$. 
${\cal W}$ is the same deficiency space as before. From the results
of this paper, there exists an inverse for the linearization about the 
part of the approximate solution which is singular only at $\Lambda'$;
we denote this inverse, as well as its Schwartz kernel, by $H'_{\be'}$.
Let $H''_{\be''}$ denote the right inverse, or Schwartz kernel, for the
linearization about the part of the approximate solution which is
singular only at $\Lambda''$, as constructed in \cite{MP}. 
Now let $H_{\be}$ be a Schwartz kernel obtained by using a partition of 
unity (e.g. the same one as used to construct the full approximate
solution) to glue together these two pieces. Then it is easy to
check that
\[
L_{\be} H_{\be} = I + R,
\]
where $R$ is of size $\max \{\e', \e''\}$. Clearly, then, $I+R$ is
invertible, and a right inverse for $L_{\be}$ is given
by $H_{\be} (I + R)^{-1}$. It is straightforward to verify the
necessary mapping properties for this operator. 

As noted above, the final step, using a contraction mapping argument to 
show the existence of a solution of the problem, is standard.

This completes the proof of Theorem~\ref{th:1.1} in full generality.

\end{document}